\documentclass{article}

\usepackage{arxiv}

\usepackage[utf8]{inputenc} % allow utf-8 input
\usepackage[T1]{fontenc}    % use 8-bit T1 fonts
\usepackage{hyperref}       % hyperlinks
\usepackage{url}            % simple URL typesetting
\usepackage{booktabs}       % professional-quality tables
\usepackage{amsfonts}       % blackboard math symbols
\usepackage{nicefrac}       % compact symbols for 1/2, etc.
\usepackage{microtype}      % microtypography
\usepackage{lipsum}
\usepackage{authblk}
\usepackage[sort&compress,square,comma,numbers]{natbib}

\usepackage{enumitem}

\usepackage{float} 
\usepackage{times}
\usepackage{epsfig}
\usepackage{graphicx}
\usepackage{amsmath, amsfonts} %, amssymb}
\usepackage{bm}
\usepackage{algorithm}
\usepackage{algorithmicx}
\usepackage{comment}
\usepackage{fullpage}
\usepackage{tikz}
\usepackage{pifont}
\usepackage{amsthm}
\usepackage{subcaption}
\newtheorem{theorem}{Theorem}
\newtheorem{definition}{Definition}

\newtheorem{lemma}{Lemma}[theorem]

\newtheorem*{remark}{Remark}

\DeclareMathOperator*{\argmin}{arg\,min}
\DeclareMathOperator{\Tr}{Tr}

\def\code#1{\texttt{#1}}

\usepackage{mathtools}
\usepackage{multirow}
\usepackage[noend]{algpseudocode}
\usepackage{booktabs}
\usepackage[export]{adjustbox}
\usepackage{url}

\usepackage[textsize=scriptsize,backgroundcolor=green]{todonotes}
\providecommand{\keywords}[1]{\textbf{\textit{Index terms---}} #1}

\usepackage{xspace}

\newcommand*{\st}{\textit{s.t.}\@\xspace}
\newcommand*{\vs}{\textit{vs.}\@\xspace}
\newcommand*{\wrt}{\textit{w.r.t.}\@\xspace}
\makeatletter
\newcommand*{\etc}{%
	\@ifnextchar{.}%
	{\textit{etc}}%
	{\textit{etc.}\@\xspace}%
}
\makeatother
\algtext*{EndWhile}
\algtext*{EndFor}
\algtext*{EndIf}
\algdef{SE}[DOWHILE]{Do}{doWhile}{\algorithmicdo}[1]{\algorithmicwhile\ #1}%

\makeatletter
\def\BState{\State\hskip-\ALG@thistlm}
\makeatother

\title{Range-Net: A High Precision Streaming SVD for Big Data Applications}

\author[1]{\textbf{Gurpreet Singh} \textsuperscript{\dag}}
\author[2]{\textbf{Soumyajit Gupta} \textsuperscript{\dag}}
\author[3]{\textbf{Matthew Lease}}
\author[4]{\textbf{Clint Dawson}}
\affil[2]{Department of Computer Science}
\affil[3]{School of Information}
\affil[4]{Oden Institute for Computational Engineering and Sciences}
\affil[1]{The University of Texas at Austin}
\affil[ ]{\texttt{\{gurpreet, smjtgupta, ml\}@utexas.edu}, \texttt{clint.dawson@oden.utexas.edu}}

\begin{document}

\maketitle

{\let\thefootnote\relax\footnote{{\dag contributed equally to this work.}}}% under the supervision of \ddag.}}}

\begin{abstract}

In a Big Data setting computing the dominant SVD factors is restrictive due to the main memory requirements. Recently introduced streaming Randomized SVD schemes work under the restrictive assumption that the singular value spectrum of the data has exponential decay. This is seldom true for any practical data. Although these methods are claimed to be applicable to scientific computations due to associated tail-energy error bounds, the approximation errors in the singular vectors and values are high when the aforementioned assumption does not hold. Furthermore from a practical perspective, oversampling can still be memory intensive or worse can exceed the feature dimension of the data. To address these issues, we present Range-Net as an alternative to randomized SVD that satisfies the tail-energy lower bound given by Eckart-Young-Mirsky (EYM) theorem. Range-Net is a deterministic two-stage neural optimization approach with random initialization, where the main memory requirement depends explicitly on the feature dimension and desired rank, independent of the sample dimension. The data samples are read in a streaming setting with the network minimization problem converging to the desired rank-r approximation. Range-Net is fully interpretable where all the network outputs and weights have a specific meaning. We provide theoretical guarantees that Range-Net extracted SVD factors satisfy EYM tail-energy lower bound at machine precision. Our numerical experiments on real data at various scales confirms this bound. A comparison against the state of the art streaming Randomized SVD shows that Range-Net accuracy is better by six orders of magnitude while being memory efficient.

\end{abstract}

\section{Introduction}

Singular Value Decomposition (SVD) is pivotal to exploratory data analysis in identifying an invariant structure under a minimalistic representation (assumptions on the structure) %. The core objective is to identify principal directions that are invariant under transformation eluding 
containing the span of resolvable information in the dataset. Finding a low rank structure is a fundamental task in applications including Image Compression \cite{de2015data}, Image Recovery \cite{brand2002incremental}, Background Removal \cite{wang2018robust}, Recommendation Systems \cite{zhang2005using} and as a pre-processing step for Clustering \cite{drineas2004clustering} and Classification \cite{jing2017multi}. With the advent of digital sensors and modern day data acquisition technologies, the sheer amount of data now requires that we revisit the solution scheme with reduced memory consumption as the target. In this work, we reformulate SVD with special emphasis on the main memory requirement that precludes it's use for big data applications.

\begin{figure}[h]
    \centering
    \includegraphics[width=\linewidth]{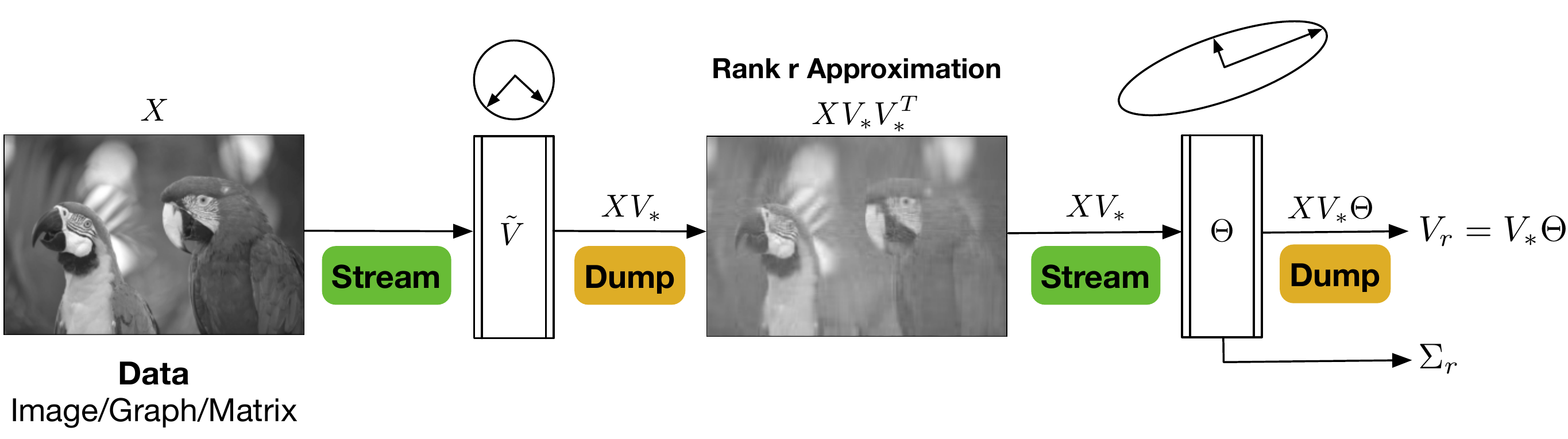}
    \caption{An overview of the low-memory, two-stage Range-Net SVD for Big Data Applications. Stage 1 identifies the span of the desired rank-r approximation. Stage 2 rotates this span to align with the singular vectors while extracting the singular values of the data. The input data can be streamed from either a server or secondary memory.}
    \label{fig:overview}
    \vspace{-5mm}
\end{figure}

It is well known that natural data matrices have a decaying spectrum wherein saving the data in memory in its original form is either redundant or not required from an application point of view. However, any assumption on the decay rate can only be made if the singular value decomposition is known \textit{a priori}, which is seldom the case in exploratory data analysis. Visually assessing a rank-$r$ approximation for image processing applications might seem correct qualitatively but are still prone to large errors due to limited human vision acuity. This is further exacerbated when the application at hand is associated with scientific computations wherein the anomalies or unaccounted phenomena are still being explored from large scale datasets. The reader is preemptively referred to \textbf{Fig. \ref{fig:sandy_sing}} where the high frequency features related to turbulence cannot be disregarded. Furthermore, for classification and clustering problems where feature dimension reduction is desirable it is imperative that a low-rank approximation of a dataset contains most $(\geq 90\%)$ of the original information content without altering the subspace information. 

In the following, we first present the problem statement describing a low rank approximation of any data matrix $X$ following Eckart-Young-Mirsky theorem \cite{eckart1936approximation,mirsky1960symmetric}. This is followed by an overview of the current limitations of low-rank SVD solvers drawing upon the problem statement, user requirements, and practical considerations in big data applications. We then present Range-Net (shown in Fig. \ref{fig:overview}) as a high-precision, low memory requirement, streaming low-rank SVD solver that achieves the EYM tail energy bound at machine precision with associated theoretical guarantees in \textbf{Appendix \ref{app:thm}}. A number of numerical experiments are presented considering practical and synthetic datasets, confirming our theoretical bounds, for detailed verification and benchmarking.

\subsection{Problem Statement}\label{sec:prob}

As per the standard notation in literature, let us denote the raw data matrix as $X \in \mathbb{R}^{m \times n}$ of rank $f \leq \min(m,n)$ and its approximation as $X_{r} \in \mathbb{R}^{m \times n}$. The singular value decomposition of $X = U\Sigma V^{T}$, where $U \in \mathbb{R}^{m \times f} = [u_1, \cdots, u_f]$ and $V \in \mathbb{R}^{n \times f} = [v_1, \cdots, v_f]$ are its left and right singular vectors respectively, and $\Sigma \in \mathbb{R}^{f \times f} = diag(\sigma_1,\cdots,\sigma_f)$ are its singular values. The rank $r \, (r \leq f)$ truncation of $X$ is then $X_r = U_r \Sigma_r V_r^T$, where $\Sigma_{r}$ is a diagonal matrix of the largest $r$ singular values of $X$, and $U_{r} = U_{[1:r]}$ and $V_{r}=V_{[1:r]}$ are the corresponding left and right singular vectors. In other words, $X = U\Sigma V^T=U_r\Sigma_rV^T_r+U_{f \textbackslash r}\Sigma_{f \textbackslash r}V^T_{f \textbackslash r} = X_{r} + X_{f \textbackslash r}$. Here, $U_{f \textbackslash r},V_{f \textbackslash r}$ are the trailing $f-r$  left and right singular vectors, respectively. The problem statement is then: Given $X \in \mathbb{R}^{m \times n}$ find $\hat{X}$ such that,
\begin{equation}
    \argmin \limits_{\substack{\hat{X} \in \mathbb{R}^{m \times n}\\\mathrm{rank}(\hat{X}) \leq r}} \|X - \hat{X}\|_{F} \label{eq:eym1}
\end{equation}
following EYM theorem \cite{eckart1936approximation,mirsky1960symmetric} equipped with a Frobenius norm. An informal equivalent is given by \textbf{Theorem \ref{thm:orig_eym}}. In effect, the minimizer $\hat{X}_{*}$ of the above problem gives us the rank-$r$ approximation of $X$ such that $X_{r} = \hat{X}_{*}$. In this work we utilize the minimizer of the above problem to extract the top rank-$r$ SVD factors of $X$ without loading the entire data matrix into the main memory.

\noindent \textbf{Target:} The above problem statement now sets our target as finding the minimizer $\hat{X}_{*}$ of the aforementioned tail energy equipped with a Frobenius norm. Please note that the minimizer naturally gives the lower bound on this tail energy in addition to being a rank-$r$ approximation.

\subsection{Limitations of Current Approaches}

Randomized low-rank SVD solvers have recently gained traction due to their low main memory requirements for big data. The first one-pass algorithm based on random projections appeared in \cite{woolfe2008fast}. Since then, a number of improvements have been proposed by \cite{halko2011finding,woodruff2014sketching,boutsidis2016optimal,upadhyay2016fast,tropp2017fixed}. The most recent and the current state of the art SketchySVD \cite{tropp2019streaming} is shown to provide consistent approximations compared to the predecessors with tighter upper bounds on tail energy. However, these current best streaming randomized SVD \cite{upadhyay2016fast,tropp2017fixed,tropp2019streaming} approaches, although claimed to perform for scientific computing, work under severely restrictive assumptions that the singular value spectrum of the data has an exponential decay. This is seldom true for any practical dataset, leading to large errors in the estimated factors. 
\begin{remark}
From a practical point of view, the user is interested in low-approximation errors in the estimated singular factors (vectors and values) and not in the tail-energy approximation errors (relative or absolute).
\end{remark}

An interesting point to note here is that even when the relative tail-energy approximation errors are small, the approximation errors in the singular factors can still be large subjective to the dataset at hand. Consequently, comparing the relative tail energies between randomized and conventional SVD solvers for different rank approximations is misleading. The reader is referred to \textbf{Fig. \ref{fig:sketchy_err}} where a synthetic dataset with non-exponential decay of singular values is considered to compare the error in singular values and this relative tail energy between SketchySVD and conventional SVD approaches. Additionally, upper-bounding the tail energy approximation error is of little consequence when the solution algorithm is not iterative and can therefore never achieve the target lower bound given by \textbf{Eq. \ref{eq:eym1}}. 

We now summarize the limitations in the current best randomized SVD approaches that motivated us to reformulate the problem in the spirit of EYM theorem:

\textbf{Tall and Skinny Matrices:} For a rank-$r$ approximation of $X \in \mathbb{R}^{m\times n}$ of rank-$f$, Randomized SVD methods rely upon rank-$k$ ($k>r$) sketches of $X$. However, these methods are useful only when $k \geq f$ but for practical datasets $f \approx min(m,n)$ and therefore the memory requirement can still be overbearing. The reader is referred to \textbf{Section \ref{sec:mnist}} for a practical example.

\textbf{Exponential Decay of Singular Values:} Assuming exponential decay implies the rank of $X$ itself is such that $f \ll \min(m,n)$. For real world applications, the data matrices are almost full rank $f=\mathcal{O}(\min(m,n))$, where a rank $r$ truncation is chosen such that the desired dominant features are accounted for. \textbf{Section \ref{sec:practical}} shows a synthetic case with non-exponential decay of singular values where sketching accrues substantial errors. We again point out that assumptions on the decay rate cannot be made preemptively on any practical dataset without performing SVD on it. Further, it is difficult to assume that the decay rate will follow a strict functional form: mixture of linear, exponential and others (see \textbf{Fig. \ref{fig:spectrum}}).

\textbf{Upper Bound on Tail Energy:} The problem statement in \textbf{Eq. \ref{eq:eym1}} suggests finding a minimum with the minimizer providing a lower bound on the tail-energy. Even if the solution scheme is upper bounded \cite{halko2011finding,upadhyay2016fast,tropp2017fixed,tropp2019streaming}, the minimizer $\hat{X}_{*}$ in Eq. \ref{eq:eym1} or equivalently achieving the lower bound is necessary.  

\textbf{Approximation Errors:} A low-rank SVD solver that does not iteratively compute the projection (left or right) while solving \textbf{Eq. \ref{eq:eym1}} cannot extract SVD factors with low errors even with multiple passes over the data matrix or multiple runs. The probability of arbitrarily selecting left and right projection matrices (over-sampled or otherwise) such that the tail energy is minimized is almost negligible. As shown later in \textbf{Theorem \ref{thm:eym_mod}}, any subspace projection (left or right) of the original data matrix that does not correspond to the minimizer in \textbf{Eq. \ref{eq:eym1}} increases the tail energy and therefore results in incorrect low-rank SVD factors (singular values and vectors). This implies that sketching matrices can only satisfy EYM lower bound \textit{if and only if} they span the rank-$r$ subspace of $X$ and hence should not be chosen arbitrarily. Please note that power iterations on sketches cannot reduce the sketching errors (once introduced) to arbitrarily small values.

\textbf{Memory Requirement:} Randomized SVD algorithms also require an optimal choice of hyper-parameters (sketch sizes \etc) that are subjective to the dataset being processed. In a practical, limited memory scenario this entails tuning the hyper-parameters to achieve optimal trade-off between memory requirement, compute time and an approximation error that does not violate the upper bound on the tail energy. The reader is referred to \textbf{Table \ref{tab:comp}} for space complexity comparison between different state of the art, streaming, randomized SVD algorithms, conventional SVD and our proposed Range-Net. Further, power iterations are either not feasible when the data matrix $X$ is itself too large to be loaded on to the main memory or the reduction in approximation errors stagnates after a few iterations. The reader is suggested to attempt the numerical experiment presented in \textbf{Section \ref{sec:practical}} with more than 5 power iterations to ratify this statement.

\noindent These limitations lead us to the following questions:
\begin{enumerate}[leftmargin=*]
    \item Can we design a low-rank, streaming, SVD solver that does not compromise accuracy for reduced memory consumption?
    \item Can we minimize the solver memory load for a desired low-rank approximation?
    \item Can we modularize low rank approximation and decomposition as two stages to increase user choices?
    \item Can we verify the results without performing a full SVD or relying upon an upper bound on tail-energy that is certainly not the minimizer of Eq. \eqref{eq:eym1}?
\end{enumerate}

Motivated by these limitations and pitfalls of the current best practices we present \textbf{Range-Net} as a low-rank SVD solver that explicitly relies upon tail-energy minimization following EYM theorem to achieve machine precision results with no assumptions on the decay rate of the singular values. Our approach also draws upon seminal work by Sagner \cite{sanger1994two}, where a  simple neural network approximates the relationship between the input and output of a linear transmission channel. This work inspired our low-weight network design where the invariants, right singular vectors and values, are extracted as network weights.

\subsection{Main Contributions}

These aforementioned limitations have not been comprehensively addressed previously, and is therefore the focus of our work in an effort to achieve practical stream processing algorithms for big-data applications. Our main contributions are as follows:

\textbf{Data and Representation Driven Neural SVD:} The representation driven network loss terms ensures that the data matrix $X$ is decomposed into the desired SVD factors such that $X = U\Sigma V^{T}$. In the absence of the representation enforcing loss term, the minimizer of Eq. \ref{eq:eym1} results in an arbitrary decomposition such that $X = ABC$ different from SVD factors. For example, for a rank-$r$ approximation of a rank-$f$ data matrix $X$ ($r \geq f$) the user can remove the representation loss from the stage-1 minimization problem to arrive at non-orthonormal vectors following \textbf{Lemma \ref{lem:1.2}} that closely satisfy (GPU precision) the EYM tail energy bound of 0.  

\textbf{A Deterministic Approach with GPU Bit-precision Results:} The network is initialized with weights drawn from a random distribution while the iterative minimization is deterministic. Although not advised for big data matrices, a full gradient descent converges to the same minimizer as a stochastic gradient descent. The streaming order of the samples is of no consequence and the user is free to choose the order in which the samples are processed in a batch-wise manner (indexed or randomized). Further, the sole difference between a full or batch-wise gradient descent is that the latter provides a low-cost approximation to the full gradient computation. We suggest the user attempt all possibilities to gain confidence in the workings of the network architecture.

\textbf{First Streaming Architecture with Exact Low Memory Cost:} Range-net requires an exact memory specification based upon the desired rank-$r$ and data dimensions $X \in \mathbb{R}^{m\times n}$ given by $r(n+r)$ and not $\mathcal{O}(r(n+r))$ independent of the sample dimension $m$. Additionally, this is the first streaming \footnote{\textbf{Note:} Range-Net does not require the entire dataset to be present in the main memory at any point during the iterative minimization.} algorithm that does require the user to wait until the streaming step is complete, contrary to sketching based randomized streaming algorithms. The stage-1 minimization problem can be started as soon as the first batch of data arrives. Further, stage-1 and stage-2 minimization problems of Range-net need not be kept separate and can be combined for a fully online, low-rank streaming SVD solver without loss of generality. This will be presented in our future work and we limit ourselves to a two-stage framework to avoid any confusion.

\textbf{Layer-wise Fully Interpretable:} Range-Net is a low-weight, fully interpretable, dense neural network where all the network weights and outputs have a precise definition and the choice of network activations is strict. The network weights are placeholders for the right (or left) orthonormal vectors upon convergence of the network minimization problems. The user can explicitly plug a ground truth solution to verify the network design and directly arrive at the tail energy bound. For example, if the ground truth low rank singular vectors are known, the user can supply them as the network weights for the stage-1 problems to retrives the singular values directly from stage-2.

% \textbf{Extendable to Eigen and PCA:} As a consequence of being a low-rank SVD solver, the architecture can be trivially extended for Eigen value and vector decomposition and principal component analysis without any loss of generality.

% \section{Motivation}
% \sg{move down to end of lit review}
% In the recent years, we are experiencing an exponential increase in the availability of high quality (high SNR) time evolving satellite image data for natural and artificial processes. Here, identifying long range, time-periodic structures will directly lead to improvements in several areas including but not limited to environmental impact assessment, Geo-sciences, planetary scale phenomena identification, and resource management. SVD is our closest ally and the most prominent tool in identifying the dominant features of a dataset when addressing this challenge. In an attempt to further this pursuit, our high precision neural SVD solver draws upon two seminal works in the field of large-scale matrix decomposition:

\section{Related Works}

Sketching based algorithms \cite{mahoney2011randomized} rely upon informed projections, where the goal is to form a sketch of the original matrix, that is small memory-wise while preserving important properties of the original space. Popular methods include Gaussian Projection \cite{johnson1984extensions} and Subsampled Randomized Hadamard Transform \cite{woodruff2014sketching}. Both of these methods involve loading the entire data into memory, forming the projection matrix and then computing a sketch. Random Sampling \cite{drineas2012fast} and Count Sketch \cite{clarkson2017low} circumvents this memory consumption by directly computing a sketch without loading the data.

The core idea behind randomized matrix decomposition is to make one or two passes over the data and compute efficient sketches. They can be broadly categorized into four main branches: 1) Sampling based methods (Subset Selection \cite{boutsidis2014near} and Randomized CUR \cite{drineas2006fast}), 2) Random Projection based QR \cite{halko2011finding}, 3) Randomized SVD \cite{halko2011finding} and 4) Power iteration methods \cite{musco2015stronger}. The sketches can represent any combination of row space, column space or the space generated by the intersection of rows and columns (core space). However, all of these methods require loading the entire data in memory. Readers are referred to \cite{kishore2017literature,ye2019fast} for an expanded survey.

Conventional SVD although deterministic and accurate, becomes expensive when the data size increases and requires $r$ passes over the data (see \textbf{Table \ref{tab:comp}}). The two branches of interest to us are the randomized SVD and Power iteration Methods for extracting SVD factors. Randomized SVD algorithms \cite{halko2011finding} are generally a two stage process:
\begin{enumerate}[leftmargin=*]
    \item \textbf{Randomized Sketching} uses random sampling to obtain a reduced matri(x/ces) which covers any combination of the row, column and core space of the data.
    \item \textbf{Deterministic Post-processing} performs conventional SVD on the reduced system from Randomized Sketching stage.
\end{enumerate}
These approaches make only one pass over the data assuming that the singular value spectrum decays rapidly. On the other hand, a Power iteration based approach \cite{musco2015stronger} requires multiple passes over the data and are used when the singular value spectrum decays slowly. This class of algorithm constructs a Krylov matrix inspired by Block Lanczos \cite{golub1977block} method to obtain a polynomial series expansion of the sketch. Although power iteration based algorithms achieve lower tail-energy errors, they cannot be used in big-data applications when $X$ itself is too large to be retained in the main memory. Here, constructing a Krylov matrix with higher order terms such as ($AA^T$ or $A^{T}A$) is not feasible\footnote{Readers are referred to \textbf{Fig. \ref{fig:err_trend_synth}} and \textbf{Fig. \ref{fig:err_trend_parrot}} to see a performance comparison between power iteration schemes and Range-Net for a synthetic and mid-sized real datasets that can be loaded into the main memory of our compute machine.}.

\begin{table*}[ht]
    \centering
     \resizebox{\columnwidth}{!}{
    \begin{tabular}{c|cccc|cc}
    \toprule
        Method & \citet{halko2011finding} & \citet{upadhyay2016fast} & \citet{tropp2017practical} & \citet{tropp2019streaming} & Range-Net & Conventional SVD \\ \midrule
        Space Complexity & $\mathcal{O}(k(m+n))$ & $\mathcal{O}(k(m+n)+s^2)$ & $\mathcal{O}(km+nl)$ &  $\mathcal{O}(k(m+n)+s^2)$ & $r(n+r)$ & $n(m+2n)$ \\
        \# Passes & 1 & 1 & 1 & 1 & $\leq 5$ & r \\ 
        Type & Randomized & Randomized & Randomized & Randomized & Deterministic & Deterministic \\ \bottomrule
    \end{tabular}}
    \caption{Current Best Randomized Streaming SVD Methods. $k,l,s$ are overestimated sketch sizes for a rank-$r$ estimate \st $k,l,s>r$. Note that Range-Net has an exact memory requirement (as conventional SVD), unlike order bounded Randomized methods.}
    \label{tab:comp}
\end{table*}

Due to main-memory restrictions on remote compute machines, streaming \cite{clarkson2009numerical} algorithms are becoming popular. For low-rank SVD approximations these involve streaming the data and updating low-memory sketches covering the row, column and core spaces. Existing randomized SVD capable of streaming include \cite{halko2011finding,upadhyay2016fast,tropp2017fixed,tropp2019streaming}, each with different sketch sizes and upper bounds on approximation errors (see \textbf{Table \ref{tab:comp}}). Amongst these, SketchySVD \cite{tropp2019streaming} is the state of the art streaming randomized SVD, with sketch sizes comparable to it's predecessors and tighter upper bounds on the tail energy and lower errors than the previous methods. As a two stage approach, SketchySVD (\textbf{Alg. \ref{alg:sketchy}}) constructs an overestimated rank sketch of the data matrix relying upon row, column and core projections. A QR decomposition on the row and column sketches then gives an estimate of the rank-$k$ subspace. This is followed by a conventional SVD on the core matrix to extract it's singular values and vectors. Finally, the singular vectors are returned after projecting them back to the original row and column space. Note that the time complexity of SketchySVD is $\mathcal{O}(k^2(m+n))$ with memory cost $\mathcal{O}(k(m+n)+s^2)$ with oversamling parameters $k=4r+1$ and $s=2k+1$.

\subsection{Need for Range-Net} \label{sec:practical}

To illustrate the limitations of current streaming randomized SVD approaches, we consider a synthetic dataset to numerically demonstrate the consequences of preemptive assumptions on the decay rate of the singular value spectrum. The numerical results section later presents a number of these singular value spectra for different practical datasets (\textbf{Fig. \ref{fig:spectrum}}) to demonstrate that the decay rates are subjective to the problem at hand. Let us consider a matrix $X$ with slow decay in the singular value.
\begin{align*}
    X = \mathrm{diag}(\underbrace{450,449,\cdots,2,1}_{f=450},\underbrace{0,\cdots,0}_{n-f})
\end{align*}
Here $X \in \mathbb{R}^{m \times n}$ is a strictly diagonal matrix with $m=n=500$ with rank $f=450$, where the singular value spectrum decays linearly. We suggest that the reader also attempt the case where all the diagonal entries are strictly ones and zeros under a high rank setting. Over multiple runs we identified the following requirements for SketchySVD to return SVD factors with relatively lower approximation errors:
\begin{enumerate}[leftmargin=*]
    \item Decay rate of singular values of a dataset must be exponential.
    \item For a rank-$f$ matrix, the desired rank $r$ must be chosen such that the oversampled rank $k$ is strictly greater than $f$ $(k \geq f)$ to achieve lower errors at scale compared to other runs.
\end{enumerate}

\begin{figure}[th]
    \centering
    \begin{subfigure}{.32\linewidth}
      \centering
      \includegraphics[width=\linewidth]{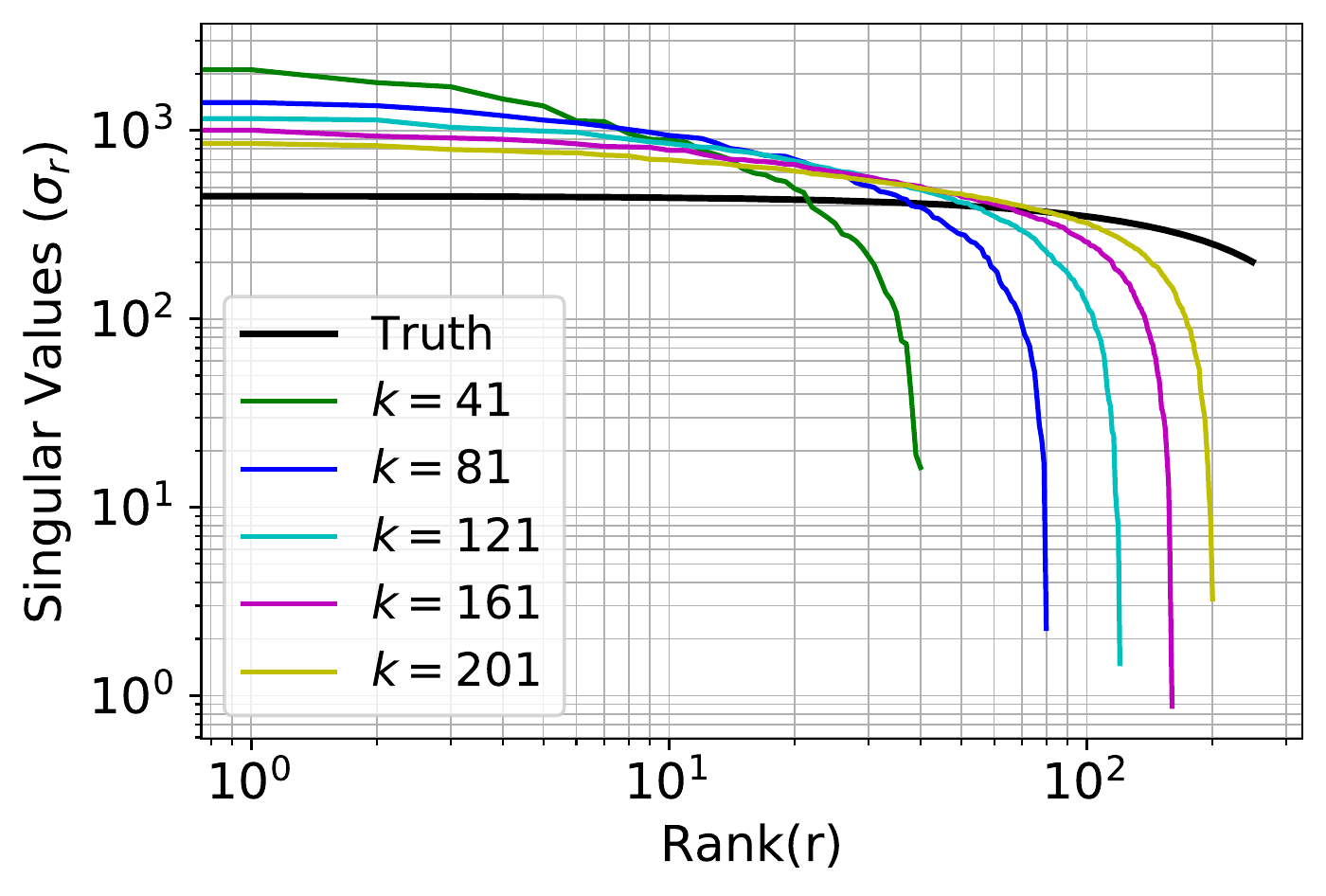}  
      \caption{Estimated Singular Values}
    \end{subfigure}
    \begin{subfigure}{.32\linewidth}
      \centering
      \includegraphics[width=\linewidth]{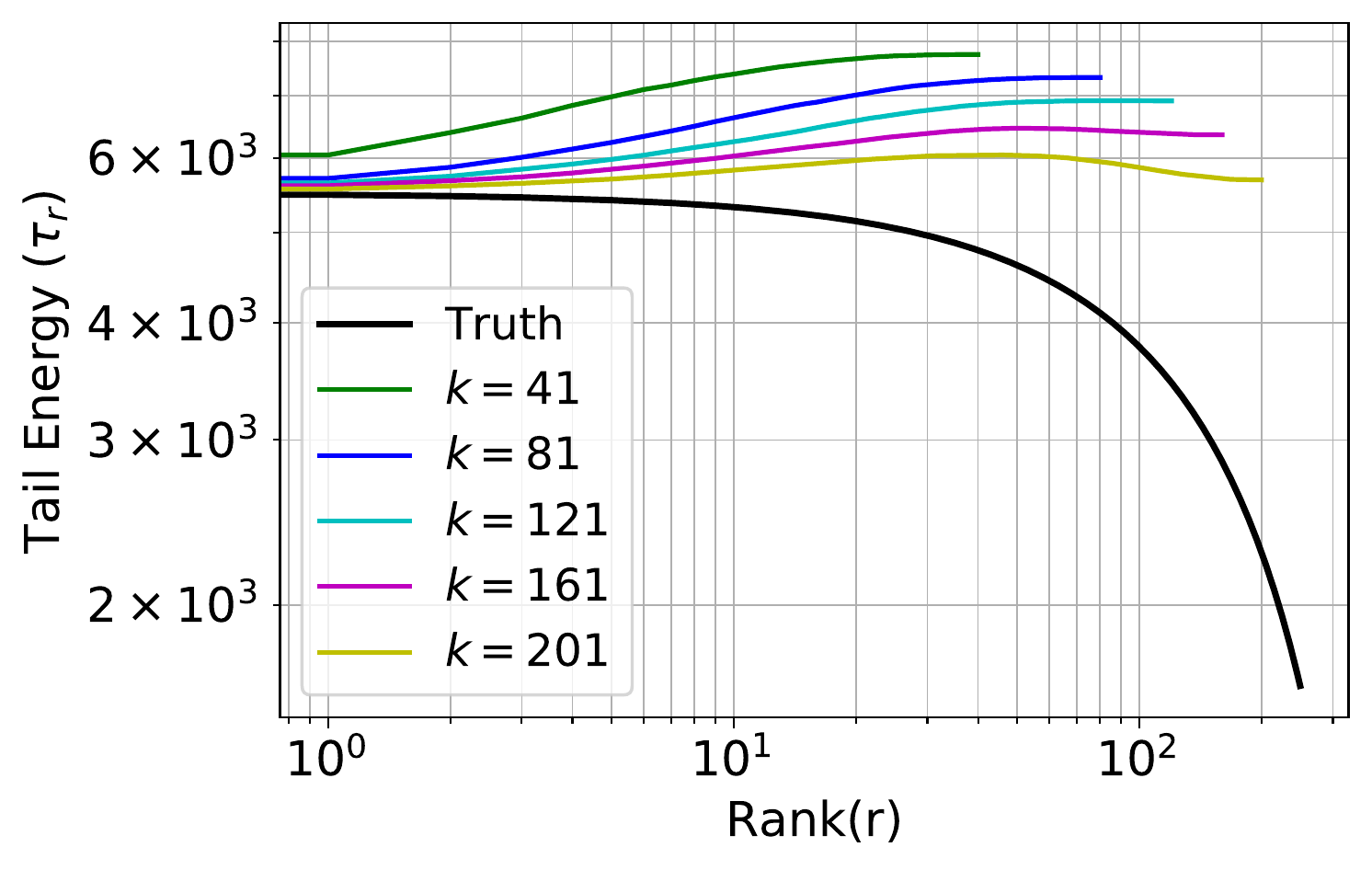}  
      \caption{Tail Energy}
    \end{subfigure}
    \begin{subfigure}{.32\linewidth}
      \centering
      \includegraphics[width=\linewidth]{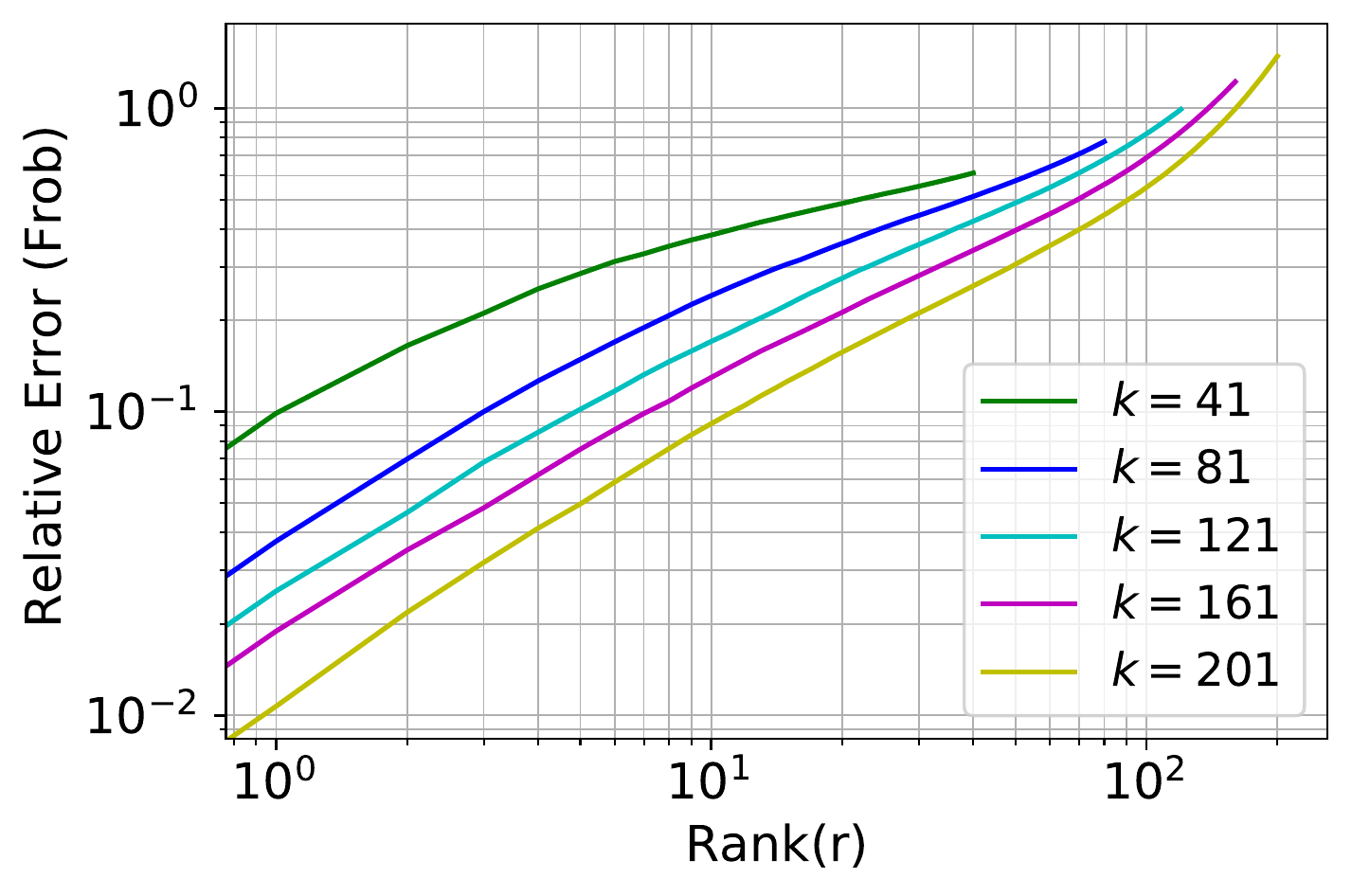}  
      \caption{Relative Tail Energy}
    \end{subfigure}
    \caption{SketchySVD approximation errors for a synthetic dataset with linear decay in singular value spectrum, corresponding to $r=[10,20,30,40,50]$ with corresponding oversampled rank $k=[41,81,121,161,201]$. Since the decay is non-exponential, SketchySVD accrues large approximation errors, hence impractical for real datasets with similar behavior.}
    \label{fig:sketchy_err}
\end{figure}

\textbf{Fig. \ref{fig:sketchy_err}} shows the singular values extracted by SketchySVD for a linearly decaying spectrum with corresponding errors in absolute and relative tail energies. The reader is referred to \textbf{Appendix \ref{app:prel}} for the definitions of tail energy and relative tail energy. Note that the synthetic data is a diagonal matrix chosen specifically so that the exact tail energies can be computed using Frobenius norm as $\left( \sum \limits_{i=r+1}^k x_{ii}^{2}\right)^{\frac{1}{2}}$. For a rank-$r$ approximation, SketchySVD suggests oversampling by a factor of $k=4r+1$ to extract the rank-$r$ factors correctly. Hence for an oversampled rank $k=[41,81,121,161,201]$ the corresponding top rank $r = [10,20,30,40,50]$ extracted singular values and vectors will have the lowest approximation errors. However as shown in the \textbf{Fig. \ref{fig:sketchy_err} (a)} the extracted singular values have an order of magnitude difference \wrt the ground truth. Consequently, \textbf{Fig. \ref{fig:sketchy_err} (b)} shows that the absolute tail-energies of the extracted features deviate quite substantially from the true tail-energy. Furthermore, we also notice that the deviations remains large as long as the oversampled rank-$k$ is such that $k<f$. For a practical dataset $f$ is either unknown or almost full rank $f \approx \min(m,n)$ or both and can only be detected by performing a full SVD on the dataset. This poses a serious restriction on SketchySVD's reliability for a realistic big data application, due to an exponential decay assumption.

\begin{remark}
A low relative error in tail energy does not imply the extracted singular values and vectors will have similar relative errors at scale. In fact, for practical datasets one can easily confirm that the scale of the relative error (as the name suggests) is subjective to the data at hand. Therefore, for a fair comparison we also show error metrics on the extracted singular factors for all our numerical experiments and forgo relative and absolute tail energy based comparisons.
\end{remark}

We also notice that for smaller values of $r$, the accrued error in both the extracted singular values and tail energy error is worse. \textbf{Fig. \ref{fig:sketchy_err} (b)} shows that for different rank approximations SketchySVD tail energy deviates from the truth quite substantially. This is due to the fact that the oversampled rank $k=4r+1 \ll f$, as pointed out before. Note that for oversampling parameter $k<f$, although the error decreases as $k \rightarrow f$ the memory requirement increases correspondingly as $\mathcal{O}(k)$ for extracting a low rank-$r$ approximation. This implies that for slow decaying spectrum optimal values of $k$ are such that $k\geq f$ even when a rank-$1$ approximation is desired. In \textbf{Fig. \ref{fig:sketchy_err} (c)}, one can observe \textit{relative errors} between $10^{-2}$ and $10^{0}$. Although this implies that the actual rank-$r$ tail-energy approximation error is off by $1\%$ at the best, the extracted singular values and vectors are off by one order of magnitude. As a consequence, the extracted singular vectors no longer represent the features of the dataset. The same issue has been raised by \citet{musco2015stronger}, who state (and we rephrase) that merely upper bounding the tail energy equipped with Frobenius or Spectral norm does not bound the approximation errors in the extracted singular values or vectors.

\begin{figure}[th]
    \centering
    \includegraphics[width=0.5\linewidth]{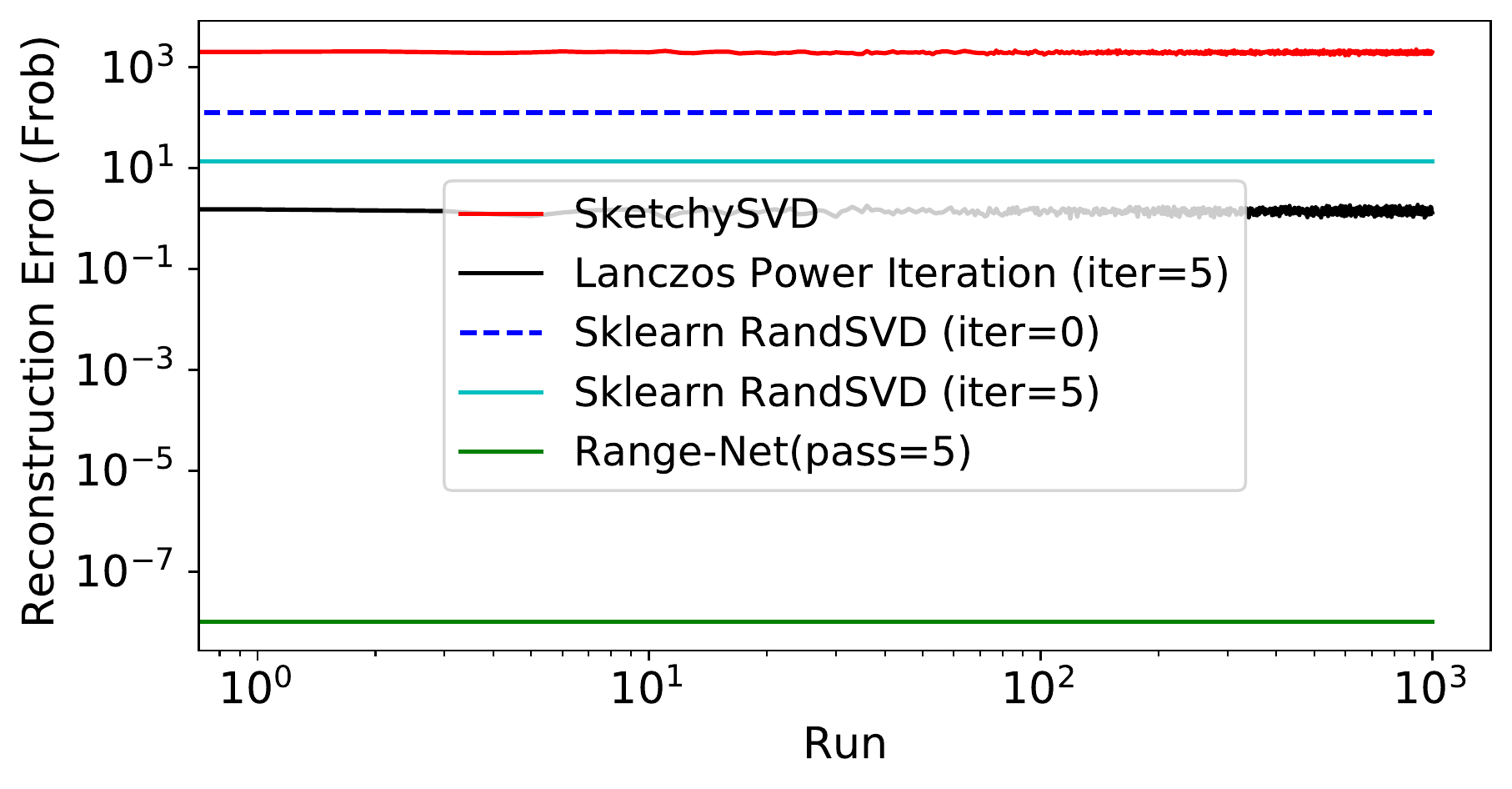}
    \caption{Comparison of reconstruction errors (desired rank 20) for  Range-Net and randomized SVD schemes (with and without power iterations) for the non-exponentially decaying singular value spectrum over 1000 runs.}
    \label{fig:err_trend_synth}
\end{figure}

\textbf{Fig. \ref{fig:err_trend_synth}} shows a comparison of reconstruction errors (see Metrics in \textbf{Section \ref{sec:metric}} for the definition) for SketchySVD \cite{tropp2019streaming} (red line), Block Lanczos with Power Iteration \cite{musco2015stronger} (black line), Sklearn's randomized SVD \cite{skrandsvd} implementation \cite{halko2011finding} with (solid cyan line) and without power iteration (dashed blue line), and Range-Net (green line) over $1000$ runs for this synthetic dataset. \textit{Please note that although power iteration improves the reconstruction error for both Block Lanczos \cite{musco2015stronger} and Sklearn's RandSVD \cite{halko2011finding}, power iteration itself requires a persistent presence of the data matrix $X$ in the main memory.} For a practical big data scenario, power iteration is therefore not a feasible alternative when the data matrix $X$ or it's sketch is itself too big to be loaded into the main memory. 

Also note that a single power iteration requires one pass over $X$ and therefore the number of Range-Net's passes and power iterations are equivalent. As is evident from the \textbf{Fig. \ref{fig:err_trend_synth}} Range-Net requires $5$-passes over the data to achieve an error which is 6 orders of magnitude better than any other Randomized method, with or without power iterations. Further, note that in \textbf{Fig. \ref{fig:err_trend_synth}} the error expectation (upper bound) over multiple runs of Randomized SVD algorithms do not reduce.

\begin{remark}
For a randomized SVD algorithm to converge (without power iterations) to a rank-$r$ approximation $X_{r}$ over multiple runs, we posit that a rank-$r$ sketch matrix $\hat{X}$ for a given rank-$f$ dataset $X$,  for $f \geq r$, be such that $P(span\{\hat{X}\} \cap span\{X_{r}\} = span\{X_{r}\}) \geq 0.5$. However, ensuring this requires substantial amount of prior knowledge or intelligent sampling (a multipass iterative algorithm).
\end{remark}

Range-net with it's explicit minimization of tail energy is capable of intelligent sampling on an arbitrary matrix without requiring any prior information. The key point to note here is that Range-Net relies upon an iterative computation of a near optimal projector instead of arbitrary/user-specified projectors used in Randomized SVD schemes. Even if the tail energy is theoretically upper bounded for some of the Randomized SVD schemes, the target is to find the lower-bound (minimizer) on the tail energy as discussed in \textbf{Section \ref{sec:prob}}. Furthermore, since none of the randomized SVD schemes construct the projector in an iterative manner while minimizing Eq. \ref{eq:eym1}, the relative error in the tail energy remains high. Even if multiple runs of SketchySVD\footnote{Power iterations can be used with SketchySVD to improve the reconstruction error at additional compute and memory cost. However, such a scheme has not been proposed in the literature and therefore it is not our prerogative to show.} or Sklearn's RandSVD are performed, the reconstruction errors in tail energies remain the same at scales shown in \textbf{Fig. \ref{fig:err_trend_synth}}. We would also like to point out that although one must strive for lower errors (relative or otherwise) and tighter theoretical upper and lower bounds, in practice we should also closely monitor if these theoretical bounds deliver us the desired solution. 

\section{Range-Net: A 2-stage SVD solver}
In the following, we present \textbf{Range-Net} that explicitly relies upon solving the minimization problem in Eq. \ref{eq:eym1} to achieve the lower bound on the tail-energy for a desired rank-$r$ approximation of a data matrix $X$ under a streaming setting. The readers are referred to \textbf{Appendix \ref{app:thm}} for some of the preliminaries followed by theorems and lemmas associated with each of the two stages.

\subsection{Network Architecture}\label{sec:arch}

\begin{figure}[h]
    \centering
    \includegraphics[width=0.6\linewidth]{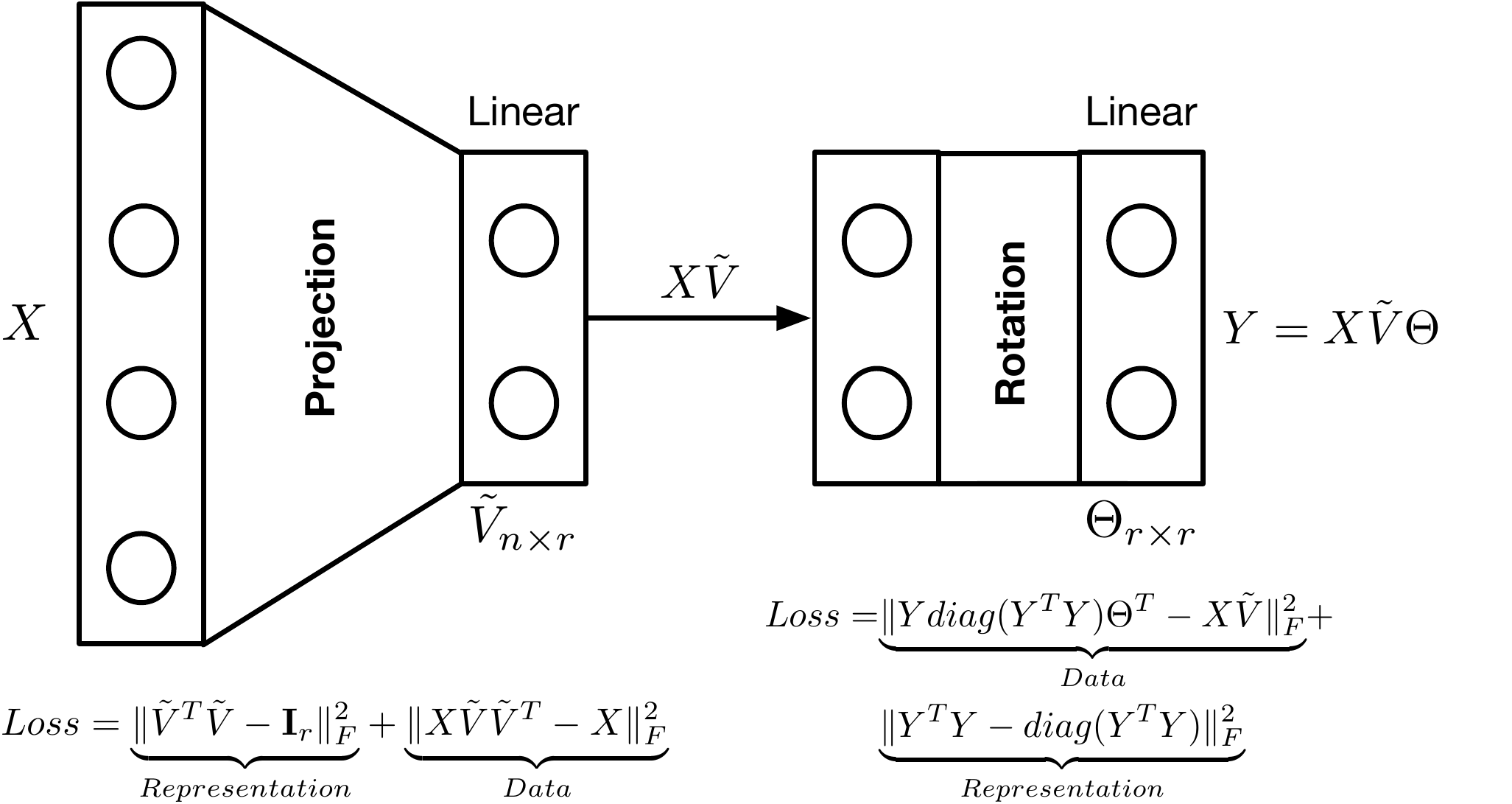}
    \caption{Network Architecture: Projection (Net1) and Rotation (Net2) for a 2-stage SVD}
    \label{fig:arch}
\end{figure}

The proposed network architecture is divided into two stages: (1) Projection, and (2) Rotation, each containing only one dense layer of neurons and linear activation functions with no biases. \textbf{Fig. \ref{fig:arch}} shows an outline of the this two-stage network architecture where all the weights and outputs have a specific meaning enforced using representation and data driven loss terms. Contrary to randomized SVD algorithms the subspace projection (Stage 1) is not specified preemptively (consequently no assumptions) but is computed by solving an iterative minimization problem following EYM theorem corresponding to Eq. \ref{eq:eym1}. The rotation stage (Stage 2) then reuses the EYM theorem again in a modified form to extract the singular vectors and values.

\textbf{Stage 1: Rank-$r$ Sub-space Identification:} The projection stage constructs an orthonormal basis that spans the $r$-dimensional sub-space of a data matrix $X \in \mathbb{R}^{m \times n}$ of an unknown rank $f \leq min(m,n)$. This orthonormal basis ($\tilde{V}$) is extracted as the stage-1 network weights once the network minimization problem converges to a fixed-point. The representation loss $\|\tilde{V}^T\tilde{V}-I_r\|_{F}$ in stage-1 enforces the orthonormal requirement on the projection space (even when $r>f$) while the data-driven loss $\|X-X\tilde{V}\tilde{V}^T\|_{F}$ minimizes the tail energy.  Although the minimization problem is non-convex the tail-energy is guaranteed to converge to the minimum at machine precision. The reader is referred to \textbf{Appendix \ref{app:converge}} for a brief discussion on the minimization problem (loss function) for details regarding the convergence behavior. 
From EYM theorem \cite{eckart1936approximation, mirsky1960symmetric}, we have that for any rank-$r$ matrix $B_{r}$, the tail energy $\|X - B_{r}\|_{F}$ in Frobenius norm $\|\cdot\|_{F}$ is bounded below by,
\begin{align}
    \|X - B_{r}\|_{F} \geq \|X - X_{r}\|_{F},
    \label{eq:eym}
\end{align}
where, $X_{r}$ is the rank-r approximation of $X$. As per \textbf{Theorem \ref{thm:eym_mod}}, the equality holds true only when $span\{B_{r}\} = span\{X_{r}\}$. Further, we define $B_{r}$ as,
\begin{align*}
    B_{r} = X\tilde{V}\tilde{V}^T
\end{align*}
then, $V_{*}V_{*}^{T} = V_{r}V_{r}^{T}$ and $V_{*}^{T}V_{*} = I_{r}$ where $V_{r}$ is a rank-r matrix with column vectors as top-r right singular vectors of $X$. The minimization problem then reads,
\begin{align}
    \underset{\tilde{V}}{\min} \quad \|X - X\tilde{V}\tilde{V}^{T}\|_{F} \quad s.t. \quad \tilde{V}^{T}\tilde{V} = I_{r}
    \label{eq:st1}
\end{align}
with a minimum at the fixed point $V_{*} = span\{v_{1}, v_{2},\ldots,v_{r}\}$ where $v_{i=1,2,\ldots,r}$ are the right singular vectors of $X_{r}$. This minimization problem describes the Stage 1 loss function of our network architecture. Upon convergence, the minimizer $\tilde{V}_{*}$ is such that $V_{*}V_{*}^{T} = V_{r}V_{r}^{T}$ following \textbf{Theorem \ref{thm:eym_mod}} where $V_{r}$ is the matrix with columns as right singular vectors of $X$ corresponding to the largest $r$ singular values of $X$.

\begin{remark}
Note that for $r\leq f$, the orthonormality constraint is trivially satisfied as shown in \textbf{Lemma \ref{lem:1.1}}. However for $r>f$, the orthonormality constraint ensures that the column vectors in $V_{*}$ are orthonormal (see \textbf{Lemma \ref{lem:1.2}}) allowing us to extract orthonormal right singular column vectors of $V_{r}$ from the Stage 2 minimization problem.
\end{remark}

One can easily identify that if a rank-1 approximation of any data matrix is desired then Stage-1 itself is sufficient to extract the singular vectors and values. Further, similar to conventional SVD the singular vectors and values can be extracted one by one in decreasing order by using Stage-1 alone. However, if the data matrix $X$ is large and the desired rank $r\gg 1$ then this naive approach will incur substantial computational and streaming costs. We therefore use Stage-1 as a flexible module to extract rank-$r$ subspace of $X$ or in other words a low rank approximation $\hat{X}=XV_{*}V_{*}$ where $V_{*}$ is the minimizer of the Eq. \ref{eq:st1}. Although not advised, a rough sketch can be generated from \textbf{Stage 1} by enforcing a looser termination criteria on the loss function. This will result in lesser number of passes required over the data matrix resulting in reduced compute time at the cost of solution accuracy.

\textbf{Stage 2: Singular Value and Vector Extraction:} The rotation stage then extracts the singular values by rotating the orthonormal vectors ($V_{*}$) to align with the right singular vectors ($V_{r} = V_{*}\Theta_{r}$). From the fixed point of the Stage-1 minimization problem Eq. \ref{eq:st1} we have $V_{*}V_{*}^{T} = V_{r}V_{r}^{T}$. According to the EYM theorem the tail energy of a rank-r matrix $XV_{*}C_{r}$, where $C_{r}$ is an arbitrary rank-r, real valued, square matrix, with respect to $XV_{*}$ is now bounded below by 0,
\begin{align*}
    \|XV_{*} - XV_{*}C_{r}\|_{F} \geq  0\\
\end{align*}
From \textbf{Theorem \ref{thm:st21}} and \textbf{Lemma \ref{lem:2.1}} we know that, $C_{r} = \Theta_{r}\Theta_{r}^{T}$, where $\Theta{r}$ is a rank-r, unitary matrix in an r-dimensional Euclidean space. Further, from \textbf{Theorem \ref{thm:st22}} we have that $(XV_{*}\Theta_{r})^{T}(XV_{*}\Theta_{r})$ is a diagonal matrix $\Sigma^{2}_{r} = \mathrm{diag}(\sigma^{2}_{1},\sigma^{2}_{2},\cdots,\sigma^{2}_{r})$, where $\sigma_{i}$s are the top-r singular values of $X$ if and only if $V_{*}\Theta_{r} = V_{r}$. Assuming  $Y =  XV_{*}\Theta_{r}$ for convenience of notation, the minimization problem now reads:
\begin{align*}
\underset{\Theta_{r}}{\min} \quad  & \|Y\Theta_{r}^{T} - XV_{*}\|_{F}\\
 & s.t.  \text{ and } \quad Y^{T}Y - \mathrm{diag}(Y^{T}Y) = 0\\
\end{align*} 

\begin{remark}
Note that stage 1 can be verified numerically independently of stage 2 by checking whether the orthonormality condition is met in addition to minimization problem converging to the tail-energy bound. Similarly, stage 2 minimization problem will return a rotation matrix $\Theta_{r}$ $(\Theta_{r}^{T}\Theta_{r} = \Theta_{r}\Theta_{r}^{T} = I_{r}, det(\Theta_{r}) = \pm 1)$ upon convergence that can again be verified numerically.
\end{remark}
As discussed previously, this choice of loss terms equipped with a Frobenius norm ensures a rank-$r$ approximation in accord with the Eckart-Young-Mirsky (EYM) theorem. We are therefore able to preemptively state that the expected values of the stage-1 loss term at the minimum must correspond to the rank $(n-r)$ tail energy. This can be verified by performing a full SVD using conventional solvers and computing a Frobenius norm on a reconstruction of the data using the bottom $(n-r)$ singular values and vectors or alternatively by computing $(\|X - X_{r} \|_{F})$. Further, the second loss term is expected to reach a machine precision zero at the minimum. 

Once, the network minimization problem converges, the singular values are extracted from \textbf{Stage 2} network weights $\Theta_{r}$ as $\Sigma_{r}^{2} = (XV_{*}\Theta_{r})^{T}(XV_{*}\Theta_{r})$. The right singular vectors can now be extracted using \textbf{Stage 2} layer weights given by $V_{r} = V_{*}\Theta_{r}$. Once $V_{r}$ and $\Sigma_{r}$ are known, left singular vectors $U_{r} = XV_{*}\Theta_{r}\Sigma_{r}^{-1}$. Please note that for $r>f$, $f-r$ singular values are zero and therefore $\Sigma_{r}^{-1}$ implies inverting the non-singular values using a threshold of $10^{-8}$.

\subsection{Choice of Activation Functions}

All the activation functions in both stages are \textit{Linear} with no biases. One might argue that this choice is not a neural approach, since all the activations are linear. However, please note that singular vectors are \textit{linearly separable orthogonal features} of a dataset, and therefore any other choice of activation function will result in approximation errors. A simple verification can be done by approximating a straight line with \textit{tanh} activation, only to realize that the approximation error $\rightarrow 0$ as the number of neurons $\rightarrow \infty$. Since singular vectors have entries in the range $[-1,1]$, a choice of \textit{relu} activation is also problematic. These arguments can also be verified by replacing linear activation in Stage 1 by any non-linear activation only to find that the tail energy bound cannot be satisfied. Note that given a small matrix, one can calculate the right singular vectors and substitute them directly as our network weights to confirm this tail energy bound.

\subsection{Data Streaming}

Given a data matrix $X^{m \times n}$, we stream the data along the smaller dimension assuming the user prescribed rank-$r$ is such that $r\leq min(m,n)$. For the sake of simplicity we assume that the data matrix has $m$ samples and $n$ features, where $m>n$, and consequently feature vectors of samples are streamed in batches. We rely upon the built-in Keras \textbf{\code{fit\_generator}} class for data streaming from the secondary memory (HDD). For a big data matrix $X$ that cannot be loaded into the main memory, this allows us to mimic the modality of data residing on an external server. Thus, given a pointer to the data, the function yields a batch of specified size for the network to train on for specific epochs. This ability saves main memory load and allows us to process bigger datasets on smaller main-memory machines than reported in prior works.

Note that for the stage-1 network to converge to a desired tolerance, we require multiple passes (empirically $\leq 5$) over the original data streamed batchwise. \textit{Therefore for Stage 1, the input data is never persistently present in the main memory of the remote machine.} The output data is dumped onto the secondary memory assuming that storing a low rank approximation is still main memory intensive. For Stage 2, this low rank approximation in the secondary memory is streamed as input, and the extracted singular values and vectors are saved in main memory. Another alternative here for the stage-2 minimization problem is to not store $XV_{*}$ on the remote machine and stream $X$ directly from the server and use the stage-1 trained network as a projector to construct $XV_{*}$ in a streaming manner. This completely removes the burden of storing the full size, low rank approximation on the remote machine performing the computations.

\subsection{Network Interpretability} \label{sec:inter}
As described before in Fig. \ref{fig:arch}, our network weights and outputs are strictly defined and incorporated as losses in the network minimization problem. In order to create a distinction, we refer to the problem informed (SVD) restrictions on the network weights as representation driven losses. This is in contrast to kernel regularization loss often considered to impose a weak requirement on the network weights to be small. The representation driven, orthonormality loss term, in Stage 1 enforces that the weights $\tilde{V}$ must be orthonormal or $(V_{*}^TV_{*}=I_r)$ for a desired rank-$r$ which is greater than the rank-$f$ of the data matrix $X$. We numerically verify the interpretability of the layer outputs and weights by considering two networks: (1) with, and (2) without the aforementioned orthonomality loss. For each of these two cases, two synthetic datasets are considered corresponding to $r\leq f$ and $r>f$. Please note that in a practical scenario $f$ is an unknown and can be determined only by performing a full SVD of $X$. Therefore, numerically testing this aspect for our solver is necessary.

\begin{figure}[h]
    \centering
    \includegraphics[width=0.8\linewidth]{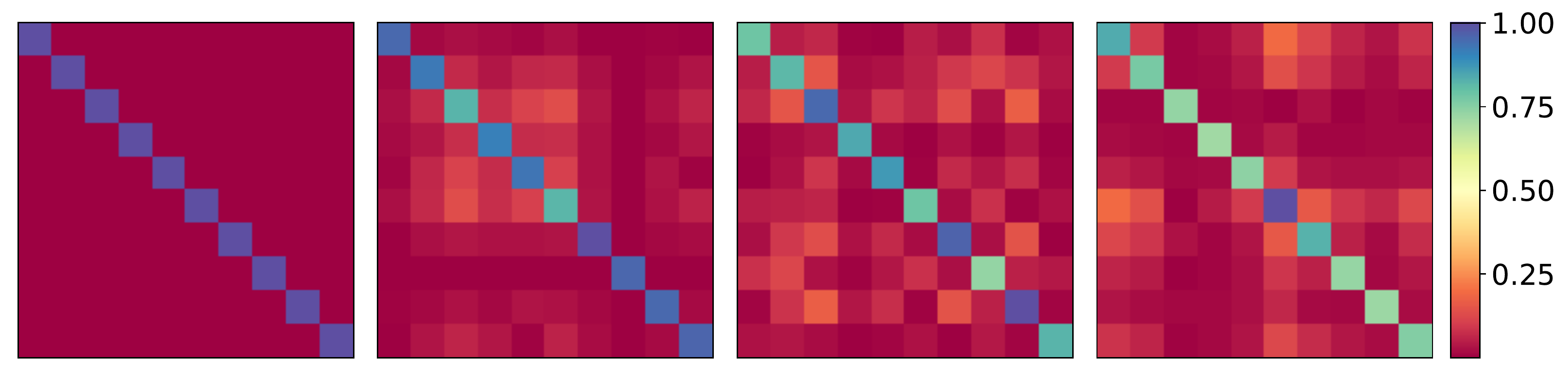}
    \caption{Synthetic Low Rank Scenario: Correlation Map of extracted vectors $V_{*}$ over four runs. Orthonormality is imposed for the first run resulting in a diagonal structure. Absence of this condition results in a scatter, as expected.}
    \label{fig:synth1}
\end{figure}

For the first case, we consider a synthetic data matrix $X_{15 \times 15}$ where the top $5$ singular values are positive ($f=5$) while the rest are zero. The objective is to extract the top $10$ ($r=10$) singular vectors where the desired rank is higher the the rank of the system itself. A total of four training runs are considered: one run for a network with the orthonormality condition imposed and three runs for another network without this additional constraints. \textbf{Fig. \ref{fig:synth1}} shows the correlation map between the recovered vectors $V_{*}$ for each of the four runs. Notice that only when the orthonormality criteria is not imposed, we get scatter away from the diagonal matrix, although all four runs converged to the same tail energy. Since the true rank of $X$ is $5$, the null space of $X$ is of dimension $10$. The absence of this orthonormality imposing, representation loss results in non-orthonormal vectors $V_{*}$ spanning the low-rank range space.

\begin{figure}[h]
    \centering
    \includegraphics[width=0.8\linewidth]{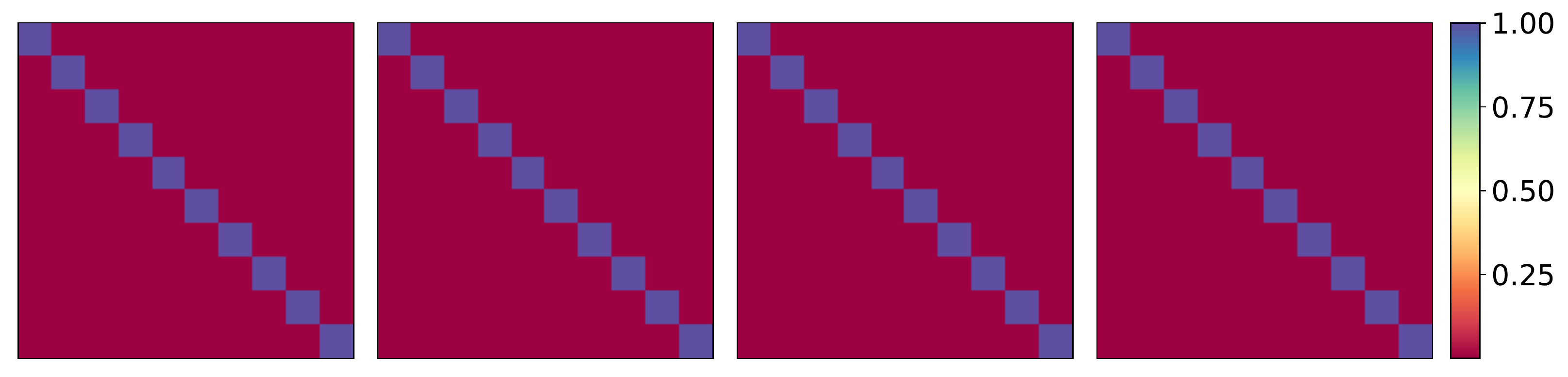}
    \caption{Synthetic Full Rank Scenario: Correlation Map of extracted vectors $\tilde{V}$ over four runs. Orthonormality loss now does not contribute and therefore all runs have a diagonal structure.}
    \label{fig:synth2}
\end{figure}

For the second case, we consider a full-rank, synthetic data matrix $X_{15 \times 15}$ with $f=15$. As before, we extract the top $10$ singular vectors ($r=10$) using four training runs: one with and three without imposing the orthonormality loss. Since the desired rank $10$ system now itself is full rank, this additional loss term does not contribute, as expected. Fig. \textbf{Fig. \ref{fig:synth2}} shows that the extracted vectors $V_{*}$ remain orthonormal, for all the four runs, so as to minimize tail energy ($\|X (\tilde{V}\tilde{V}^T - I)\|_{F}$), as described in \textbf{Section \ref{sec:arch}} above. In fact, for any rank $r$ approximation of a rank $f$ system such that $r\leq f$, an arbitrary non-orthonormal matrix $V_{*}$ will increase the tail energy and hence will not be a fixed point (solution) of our network minimization problem.

\section{Results}
We begin this section by first defining error metrics for comparison and bench-marking purposes. In the subsequent subsections, we present our training setup, results and analysis for various synthetic and real datasets, which vary in scale from small to big data. Please note that our numerical results do not show large variations in these error metrics over multiple runs due to a precise low-weight architecture. For all of the following numerical experiments, the \textbf{Stage 1} of our neural SVD solver requires at most 5 passes (empirical observation) over the data matrix to converge. Further, the error metrics rely upon conventional SVD as the baseline for a fair comparison.

\begin{figure}[h]
    \centering
    \begin{subfigure}{.32\linewidth}
      \centering
      \includegraphics[width=\linewidth]{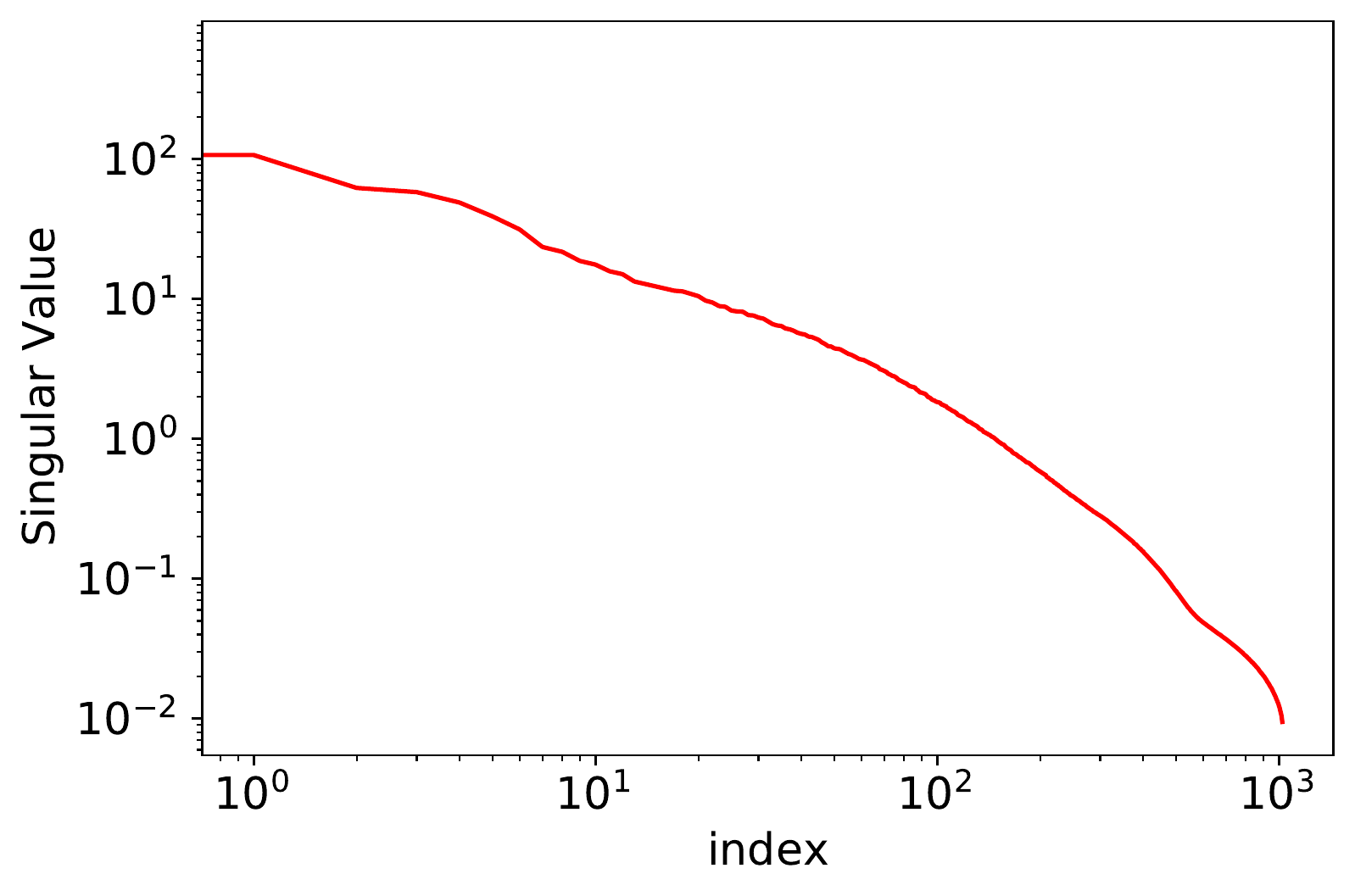}  
      \caption{Parrot}
    \end{subfigure}
    \begin{subfigure}{.32\linewidth}
      \centering
      \includegraphics[width=\linewidth]{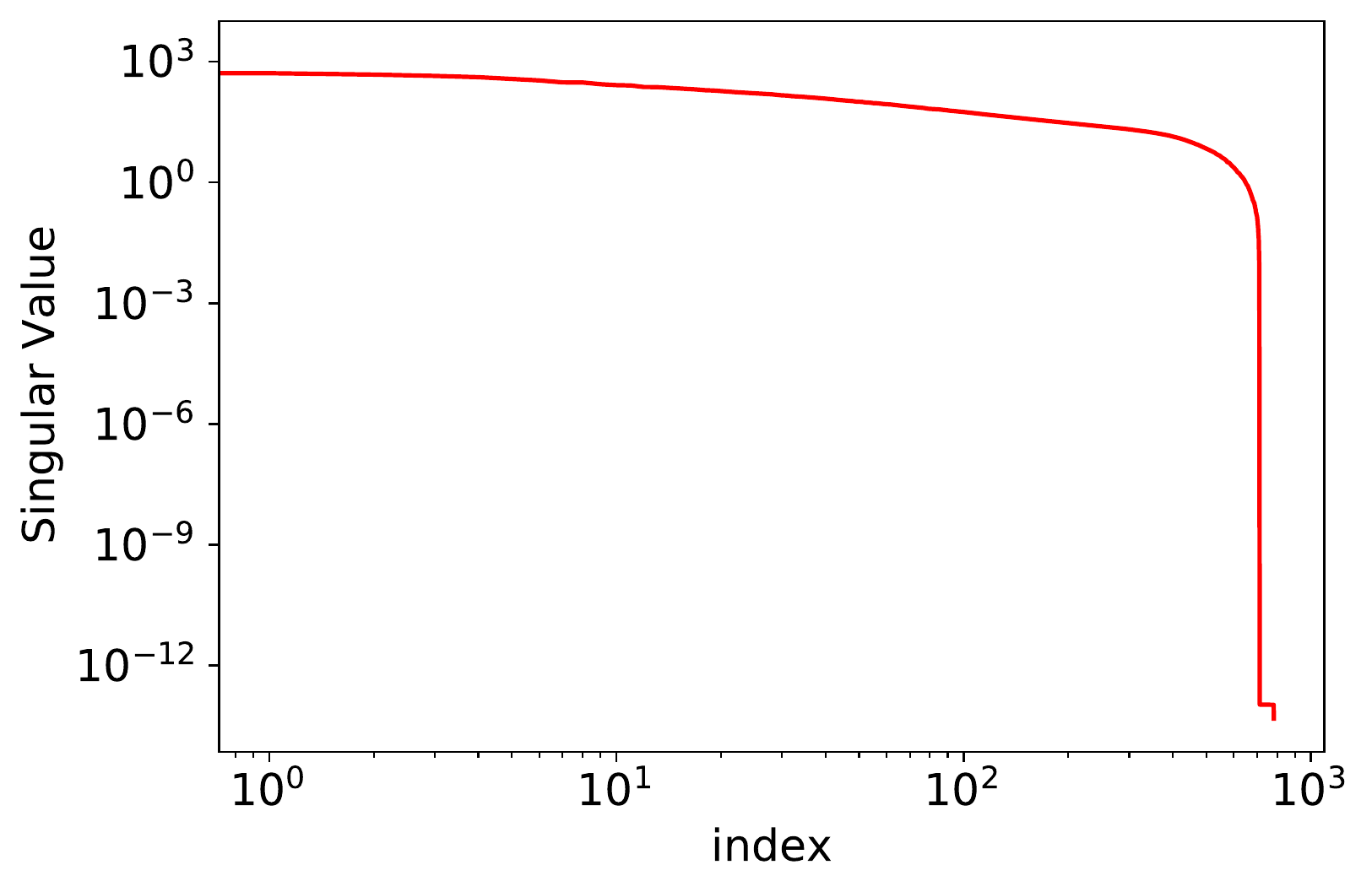}  
      \caption{MNIST}
    \end{subfigure}
    \begin{subfigure}{.32\linewidth}
      \centering
      \includegraphics[width=\linewidth]{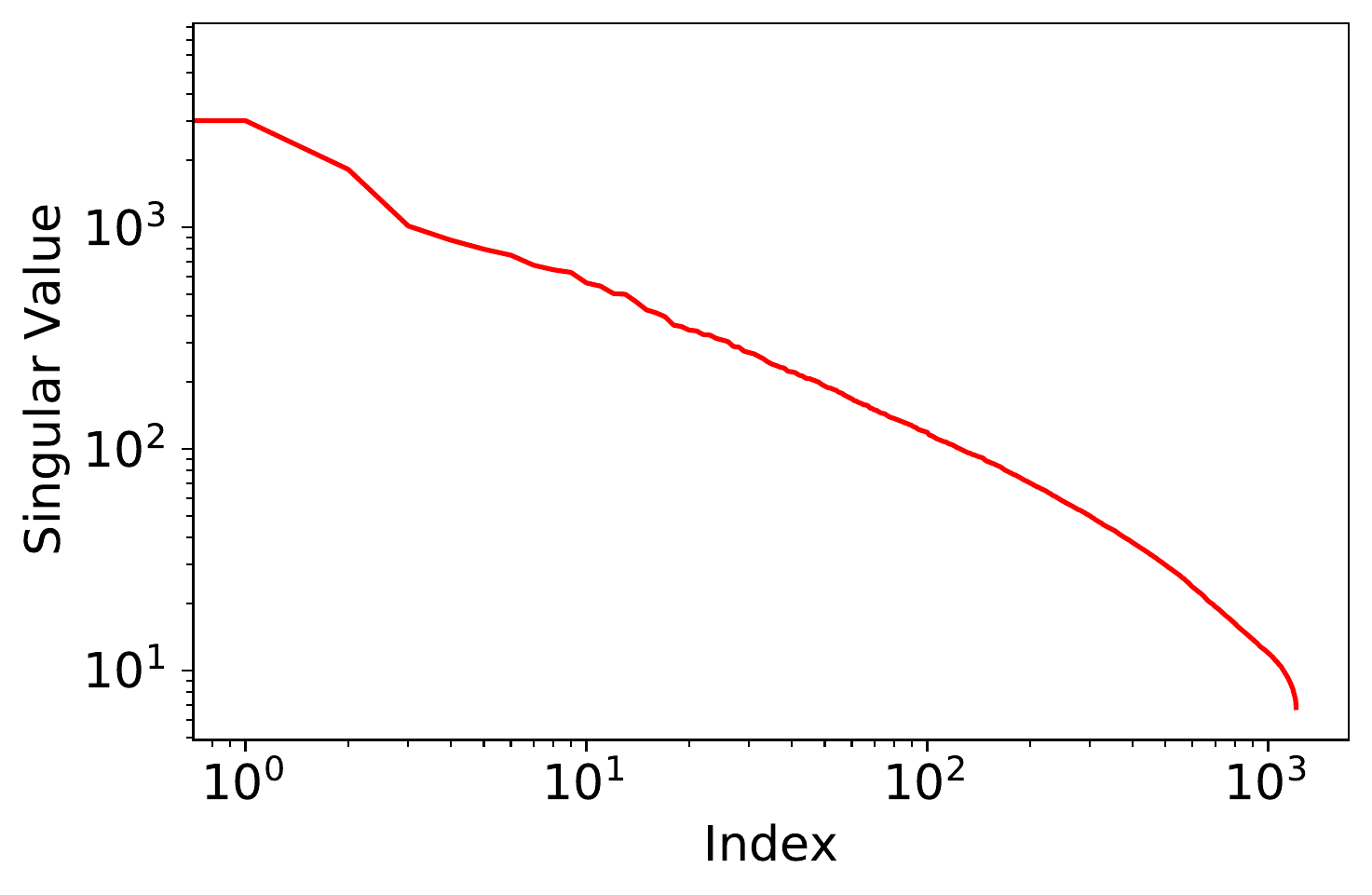}  
      \caption{Sandy}
    \end{subfigure}
    \caption{Singular value spectrum for the three practical datasets considered in this work. One can visually assess that the decay rate of the singular values is non-exponential.}
    \label{fig:spectrum}
\end{figure}

\subsection{Metrics} \label{sec:metric}

As discussed before in \textbf{Section \ref{sec:practical}} since relative errors in tail energies do not imply similar errors at scales in extracted singular factors, we rely upon additional error metrics on the extractor factors for performance comparison and benchmarking. In the following $X$ and $\hat{X}$ are used to denote the true and the reconstructed data matrices.
\begin{itemize}
    \item \textbf{Scree Error:} Absolute values of the difference between singular values from our proposed approach and rank-$r$ SVD.
    \begin{align*}
        scree_{err}(r) = |\sigma_i(X) - \hat{\sigma}_i(\hat{X})| \quad \forall \, i \in [1,r]
    \end{align*}
    % \item \textbf{Trace Error:} Sum of the differences in the top-$r$ singular values using our proposed approach and SVD.
    % \begin{align*}
    %     trace_{err}(r) = \sum_{i=1}^r |\sigma_i(X) - \hat{\sigma}_i(\hat{X})|
    % \end{align*}
    \item \textbf{Reconstruction Error:} Frobenius norm of the element wise error of the true data and its rank-$r$ approximation.
    \begin{align*}
        frob_{err}(r) = \|X - \hat{X}\|_F^2 - \|X - X_r\|_F^2
    \end{align*}  
    \item \textbf{Spectral Error}: 2-norm of the singular value of the true data and its rank-$r$ approximation.
    \begin{align*}
        spectral_{err}(r) = \|X-\hat{X}\|_2 - \|X-X_r\|_2 
    \end{align*}
    \item \textbf{Chi Square Statistic}: Deviation of the network recovered singular vectors from the true singular-vectors.
    \begin{align*}
        \chi^2_err(r) = 1 - \frac{1}{r} \|Corr(v_{[:r]}(X),\hat{v}_{[:r]}(\hat{X}))\|_F
    \end{align*}
\end{itemize}
Here, $\sigma_{i}$s are the true singular values and $X_{r}$ is the desired rank-$r$ approximation of $X$ using conventional SVD as the baseline for benchmarking. Under perfect recovery, all the error metrics are expected to approximately achieve zero at machine precision. All of our numerical experiments were performed on a GPU using single (32-bit) precision floating point operations. Therefore, the tail energies are expected to be correct to upto 8 significant digits approximately in all the subsequent calculations. In the following sections, $\hat{X}$ is replaced by approximations from Randomized SVD algorithms and Range-Net.

\subsection{Setup and Training}

All experiments were done on a setup with Nvidia 2060 RTX Super 8GB GPU, Intel Core i7-9700F 3.0GHz 8-core CPU and 16GB DDR4 memory. We use Keras \cite{chollet2015} library running on a Tensorflow 2.0 backend with Python 3.7 to train the networks presented in this paper. For optimization, we use AdaMax \cite{kingma2014adam} with parameters (\textit{lr}= 0.001) and $2000$ steps per epoch. The batch-sizes vary with dataset sizes and are therefore not reported explicitly.

\subsection{Image Compression: Parrots (SVD)}\label{sec:parrot}

\begin{figure}[th]
    \centering
    \begin{subfigure}{.35\linewidth}
      \centering
      \includegraphics[width=0.8\linewidth]{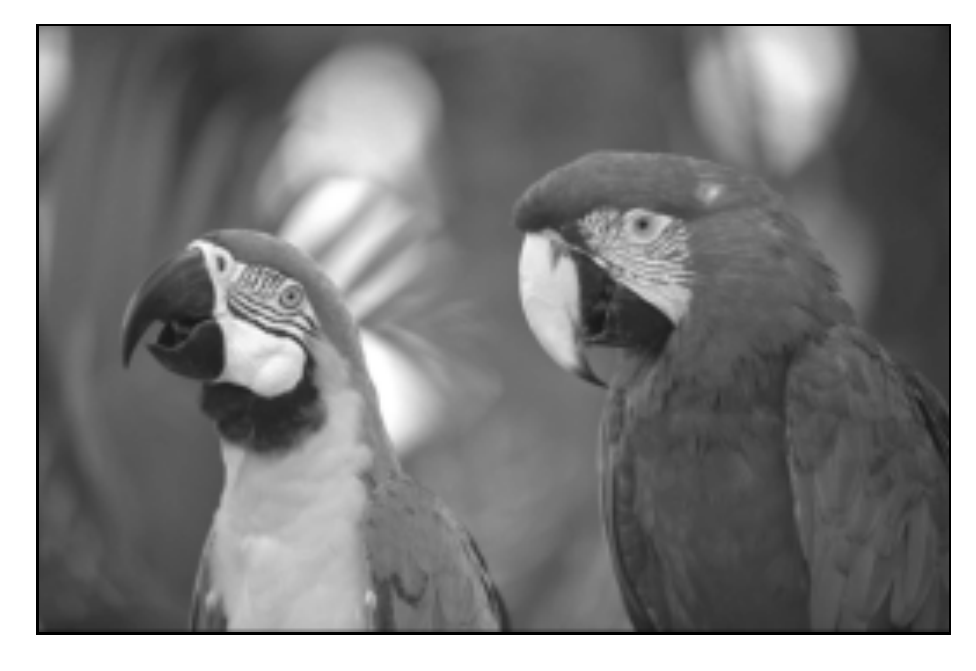}  
      \caption{True Image}
    \end{subfigure}
    \begin{subfigure}{.35\linewidth}
      \centering
      \includegraphics[width=0.8\linewidth]{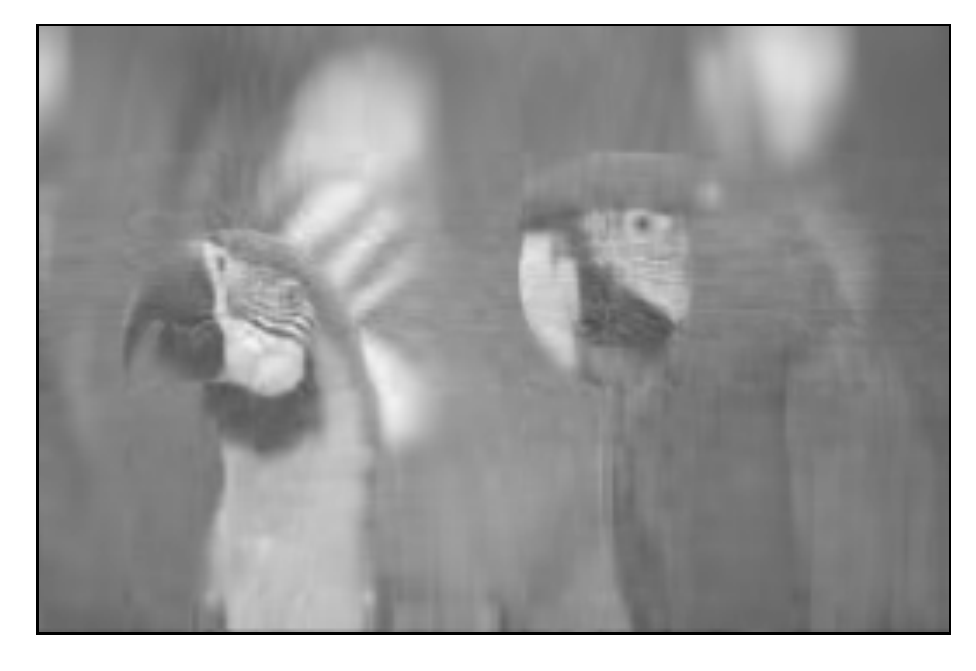}  
      \caption{$X_r$ for $r=20$}
    \end{subfigure}
    \begin{subfigure}{.35\linewidth}
      \centering
      \includegraphics[width=0.8\linewidth]{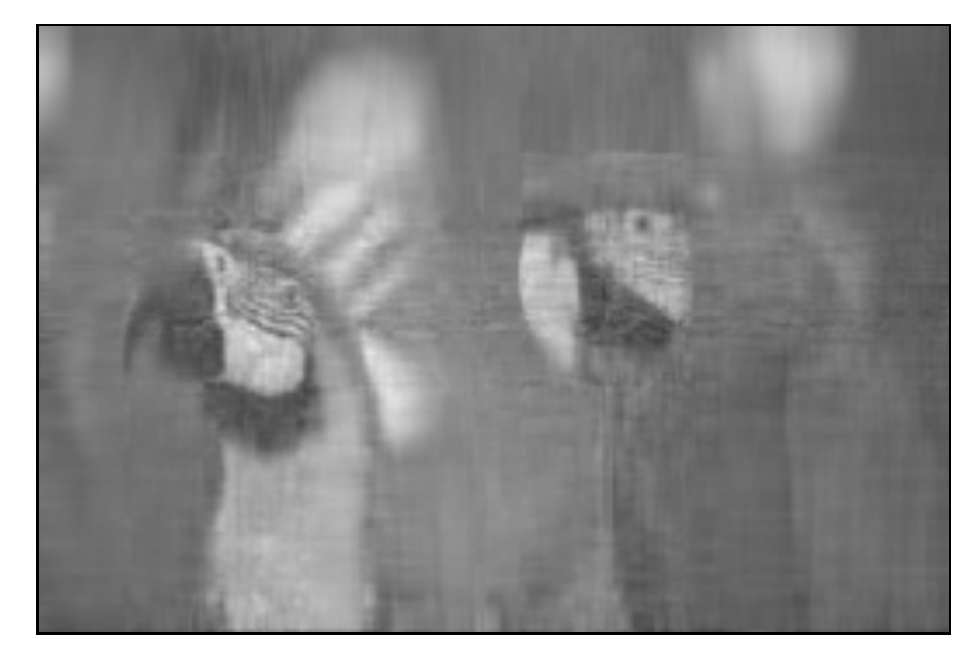}  
      \caption{SketchySVD $r=20$}
    \end{subfigure}
    \begin{subfigure}{.35\linewidth}
      \centering
      \includegraphics[width=0.8\linewidth]{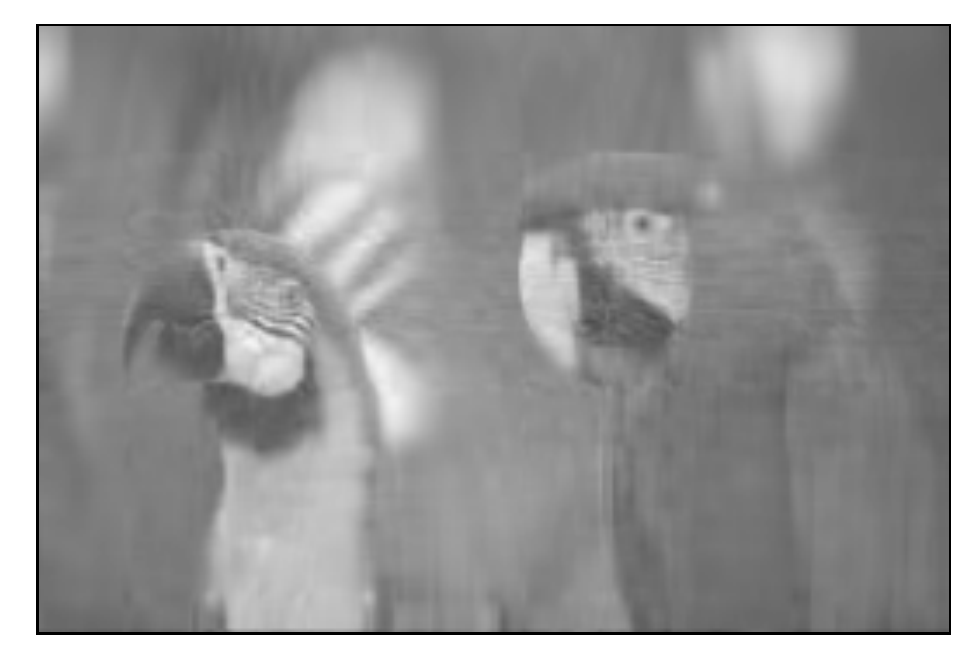}  
      \caption{Range-Net $r=20$}
    \end{subfigure}
    \begin{subfigure}{.35\linewidth}
      \centering
      \includegraphics[width=0.9\linewidth]{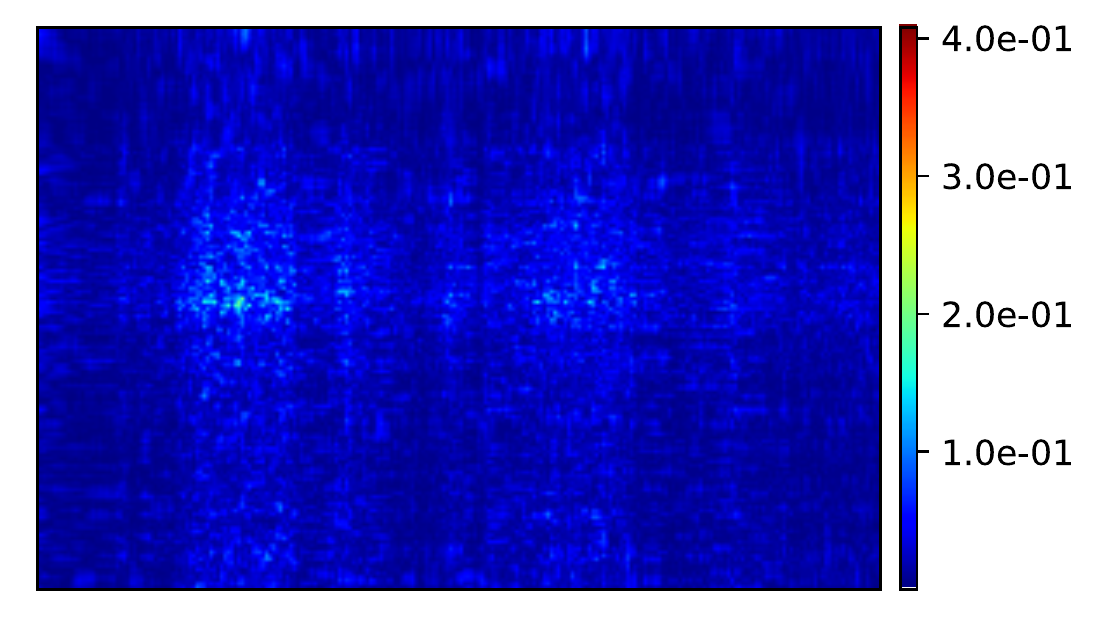}  
      \caption{$|\tilde{X}^{sketchy} - X_r|$}
    \end{subfigure}
    \begin{subfigure}{.35\linewidth}
      \centering
      \includegraphics[width=0.9\linewidth]{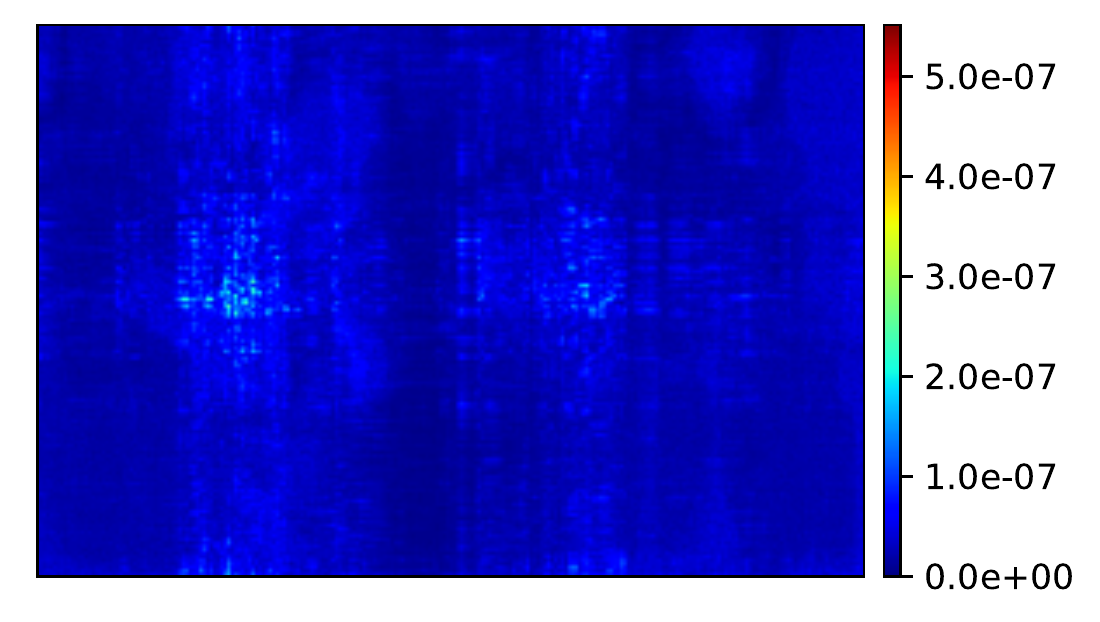}  
      \caption{$|\tilde{X}^{net} - X_r|$}
    \end{subfigure}
    \caption{(a) True image, rank-$20$ reconstruction using (b) conventional SVD, (c) SktechySVD (oversampled rank $k=81$), (d) Range-Net (5-passes). Note that Sketchy SVD reconstruction error $(10^{-1})$ is 6 orders of magnitude apart from Range-Net's reconstruction error $(10^{-7})$.}
    \label{fig:err_parrot}
\end{figure}

As an example for SVD of natural images, we use the well known Parrots image from the image processing domain. The original image is in an RGB format, converted to a gray scale for demonstration purposes followed by normalization between $[0,1]$. This ${1024 \times 1536}$ data matrix is then used to compute a rank $r=20$ approximation for comparison and numerical analysis. \textbf{Fig. \ref{fig:err_parrot}} shows the result of the low rank reconstruction for SVD, SketchySVD and Range-Net. Visually one can verify that Fig. \ref{fig:err_parrot} (b,d) are similar while (c) is different. To make the error in approximation more clear, we plot the absolute difference of SketchySVD and our net from the truncated rank image. Fig. \ref{fig:err_parrot} (e,f) shows the the corresponding plots with heatmaps imposed for clarity. Notice that while the reconstruction error for our network $(\approx 10^{-7})$ is close to the GPU precision, SketchySVD has significantly higher error scale $(\approx 10^{-1})$, validating the artifacts in the approximated image.

\begin{figure}[h]
    \centering
    \begin{subfigure}{.3\linewidth}
      \centering
      \includegraphics[width=\linewidth]{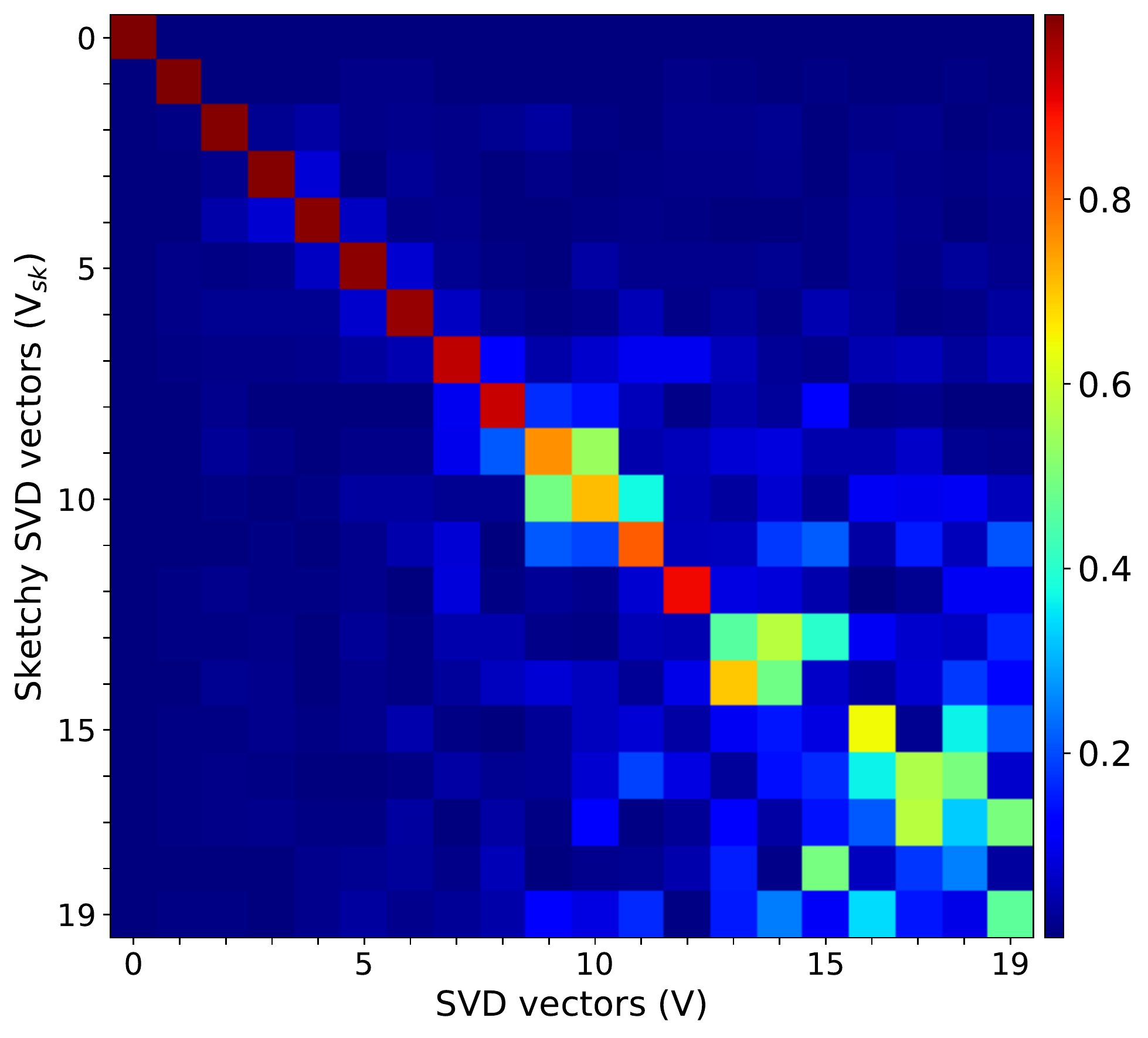}  
      \caption{SketchySVD}
    \end{subfigure}
    \qquad
    \begin{subfigure}{.3\linewidth}
      \centering
      \includegraphics[width=\linewidth]{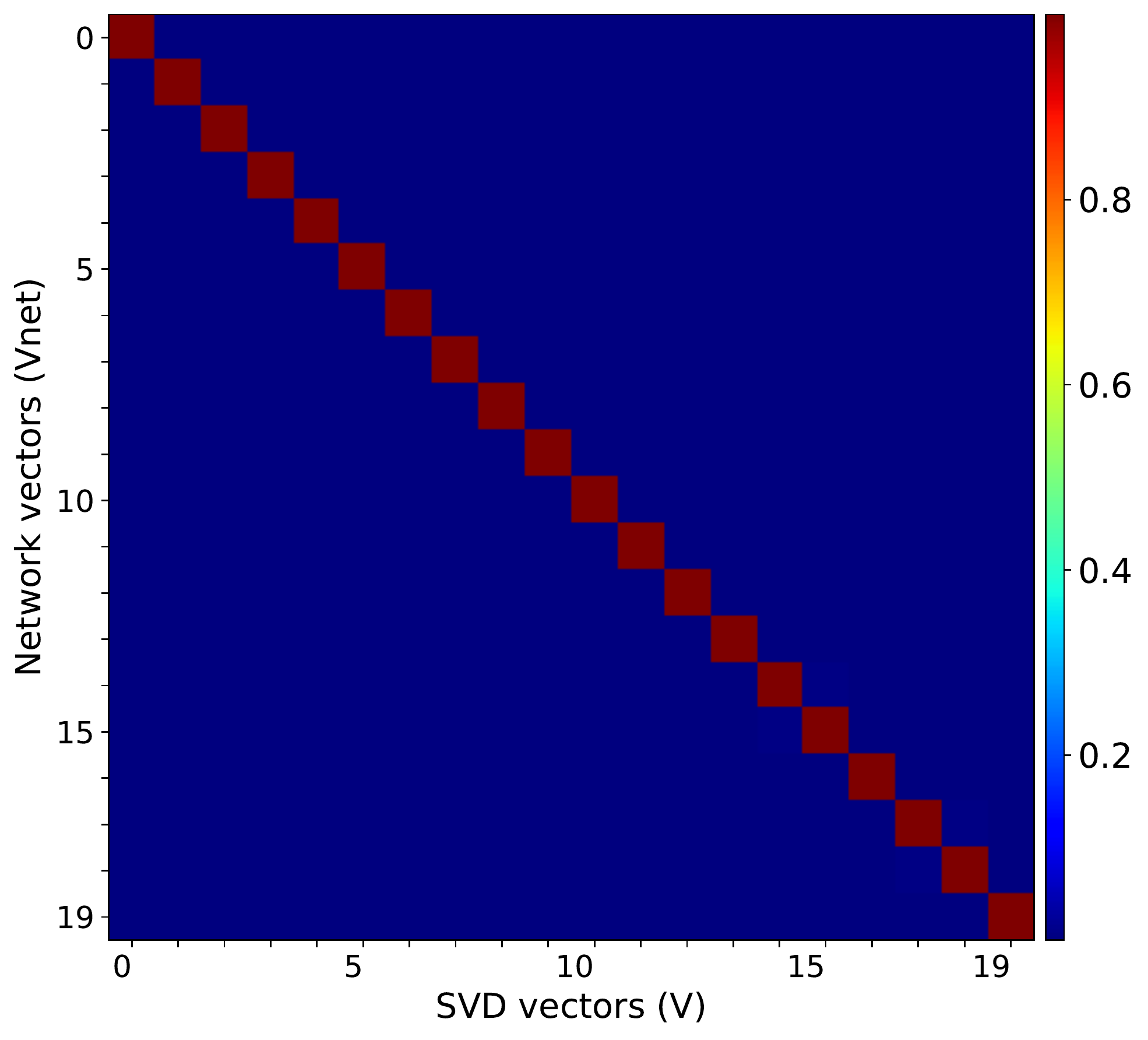} 
      \caption{Range-Net}
    \end{subfigure}
    \caption{Cross-correlation between true (conventional SVD) and extracted right singular vectors from (a) SketchySVD (b) Range-Net for a rank-20 ($r$ = 20) approximation of the full rank ($f$ = 1024) Parrot image.}
    \label{fig:corr_parrot}
\end{figure}

\textbf{Fig. \ref{fig:corr_parrot}} shows the cross-correlation between extracted right singular vectors from SketchySVD (left) and Range-Net (right) against conventional SVD for a rank-$20$ approximation of the Parrot image. Notice that SketchySVD oversampled rank is $k=81$ and still the extracted right singular vectors deviate substantially. This implies that the extracted vectors from SketchySVD do not span the top rank-$20$ subspace of $X$ as opposed to Range-Net where stage-1 explicitly ensures this span without any oversampling. The lower the vector index, the higher the spread, owing to random projections. Our method on the other hand has a near-perfect cross-correlation with the true vectors, indicated by the solid diagonal and zero off-diagonal (near GPU-precision) entries.

Note that the singular value spectrum does not decay exponentially (\textbf{Fig. \ref{fig:spectrum}}) and the data matrix is near-full rank ($f \approx 1024$). For SketchySVD to generate low-error solutions either the oversampling parameter $k$ has to be chosen such that $k\geq f$ or power iterations are required. The former is more expensive than performing a full, conventional SVD while the latter requires expensive power iterations that are not feasible for big data applications. 

\begin{figure}[h]
    \centering
    \begin{subfigure}{.4\linewidth}
      \centering
      \includegraphics[width=\linewidth]{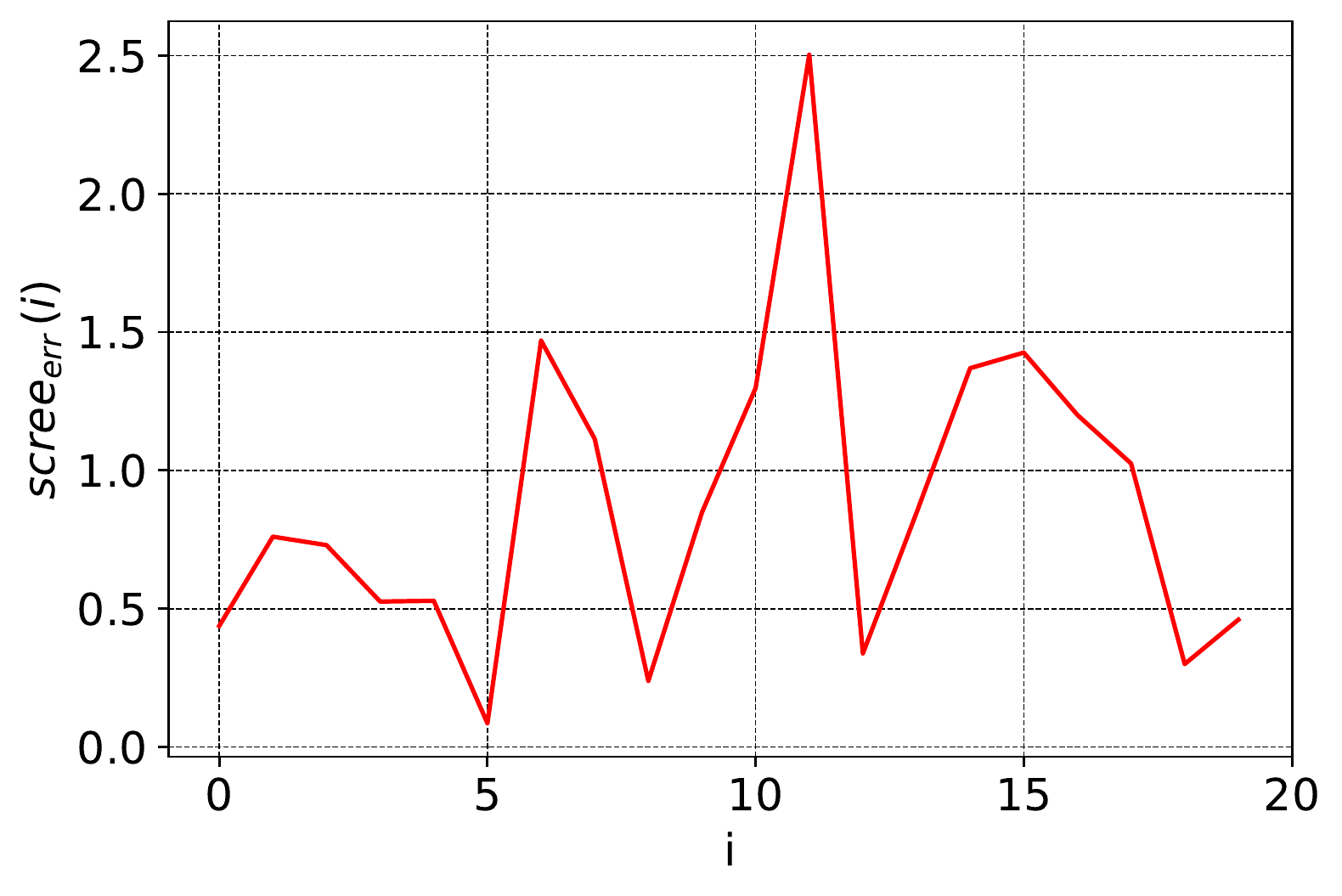}  
      \caption{SketchySVD}
    \end{subfigure}
    \qquad
    \begin{subfigure}{.4\linewidth}
      \centering
      \includegraphics[width=\linewidth]{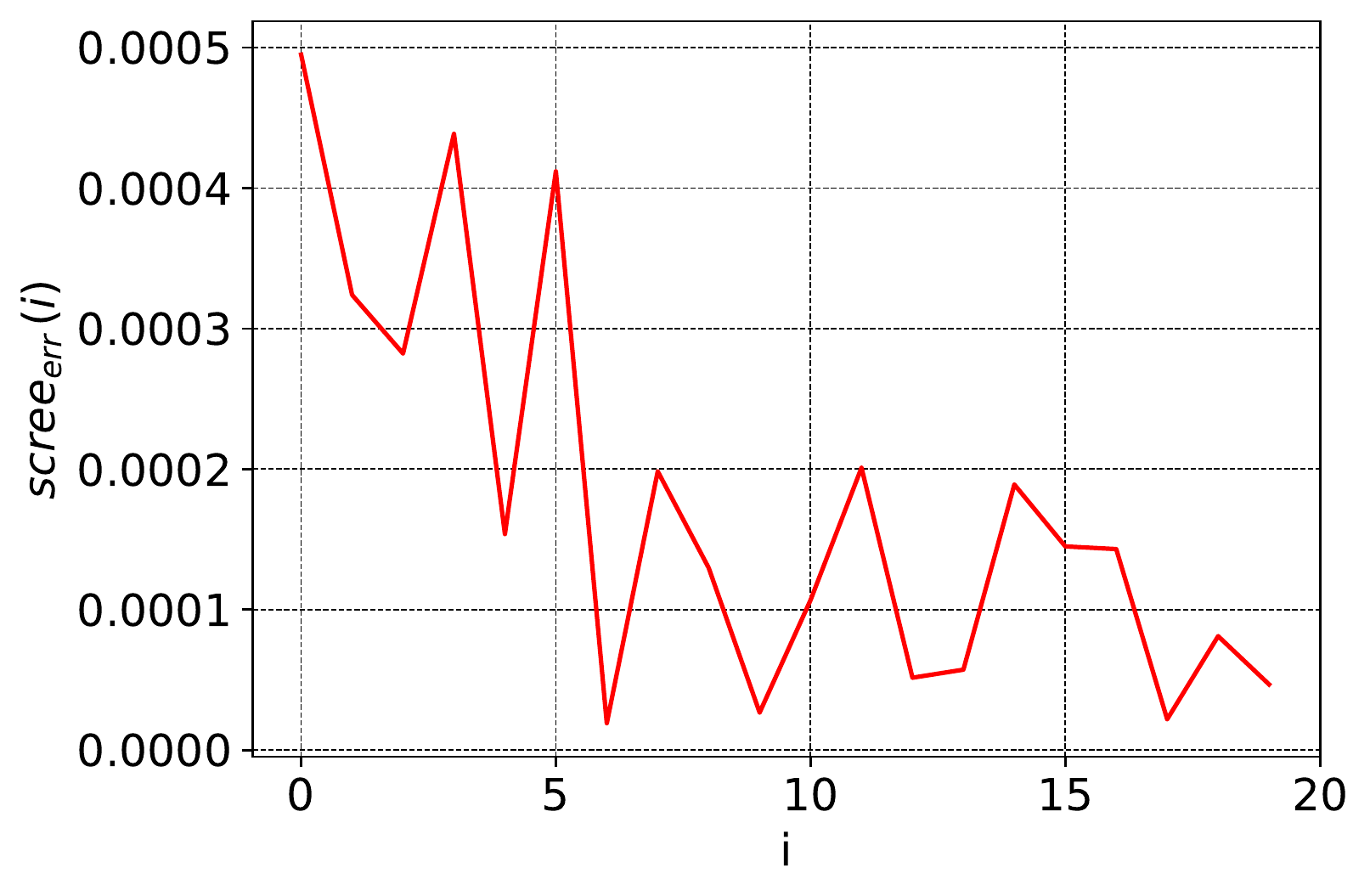}  
      \caption{Range-Net}
    \end{subfigure}
    \caption{Scree error in the extracted singular values from SketchySVD (left) ($\approx 1$) and Range-Net (right) ($\approx 10^{-4}$). Notice the scale difference in the errors. }
    \label{fig:scree_parrot}
\end{figure}

\textbf{Fig. \ref{fig:scree_parrot}} shows the scree-errors for the two methods, with the absolute difference between the predicted and the true singular values indexed by the decreasing order of singular values. For SketchySVD the error fluctuates across the top $r=20$ values, but also the scale of fluctuations is around $1$. Comparably, our method incurs significantly lower errors in singular values at scale of $(10^{-4})$.

We tabulate all metric scores for SketchySVD and the proposed method on the same image for various ranks $(10,20,50,100)$ in \textbf{Table. \ref{tab:errr-parrot}}. One can easily notice that as the expected rank increases, SketchySVD's performance keeps on deteriorating for the $\chi^2_{err}$ metric, as it is only guaranteed high precision solution for only few top values, and rest spread their energy as evident from Fig. \ref{fig:corr_parrot}. Our method gives consistent low errors on all the metrics showing it accuracy in approximating the true SVD. 
\begin{table}[h]
  \centering
  \caption{Metric Performance of SketchySVD  \vs Range-Net}
  \label{tab:errr-parrot}
  \begin{tabular}{l|c|ccc}
    \toprule
    Method & rank & $err_{fr}$ & $err_{sp}$ & $\chi^2_{err}$ \\
    \midrule
    SketchySVD & \multirow{2}{*}{r=10} & 27.904 & 18.071 & 0.492\\
    Range-Net &  & 0.0 & 0.0 & 0.018  \\ \midrule
    SketchySVD & \multirow{2}{*}{r=20} & 11.974 & 2.201 & 0.662\\
    Range-Net & & 0 & 0 & 0.023  \\ \midrule
    SketchySVD & \multirow{2}{*}{r=50} & 2.772 & 0.181 & 0.762\\
    Range-Net & & 1.91e-7 & 0 & 0.027 \\ \midrule
    SketchySVD & \multirow{2}{*}{r=100} & 0.614 & 2.14e-2 & 0.923\\
    Range-Net &  & 2.32e-7 & 1.08e-7 & 0.033 \\
  \bottomrule
\end{tabular}
\end{table}

\begin{remark}
Note that for image compression where the application is solely compressed storage and latter visualization for human interaction, \textbf{Fig. \ref{fig:err_parrot} (b), (c), (d)} look similar. Therefore one can justify using any randomized SVD scheme for this application due to limited human vision acuity. However, for scientific applications visual distinction is the least form of requirement and accurate feature resolution is of utmost importance. 
\end{remark}

\begin{figure}[th]
    \centering
    \includegraphics[width=0.5\linewidth]{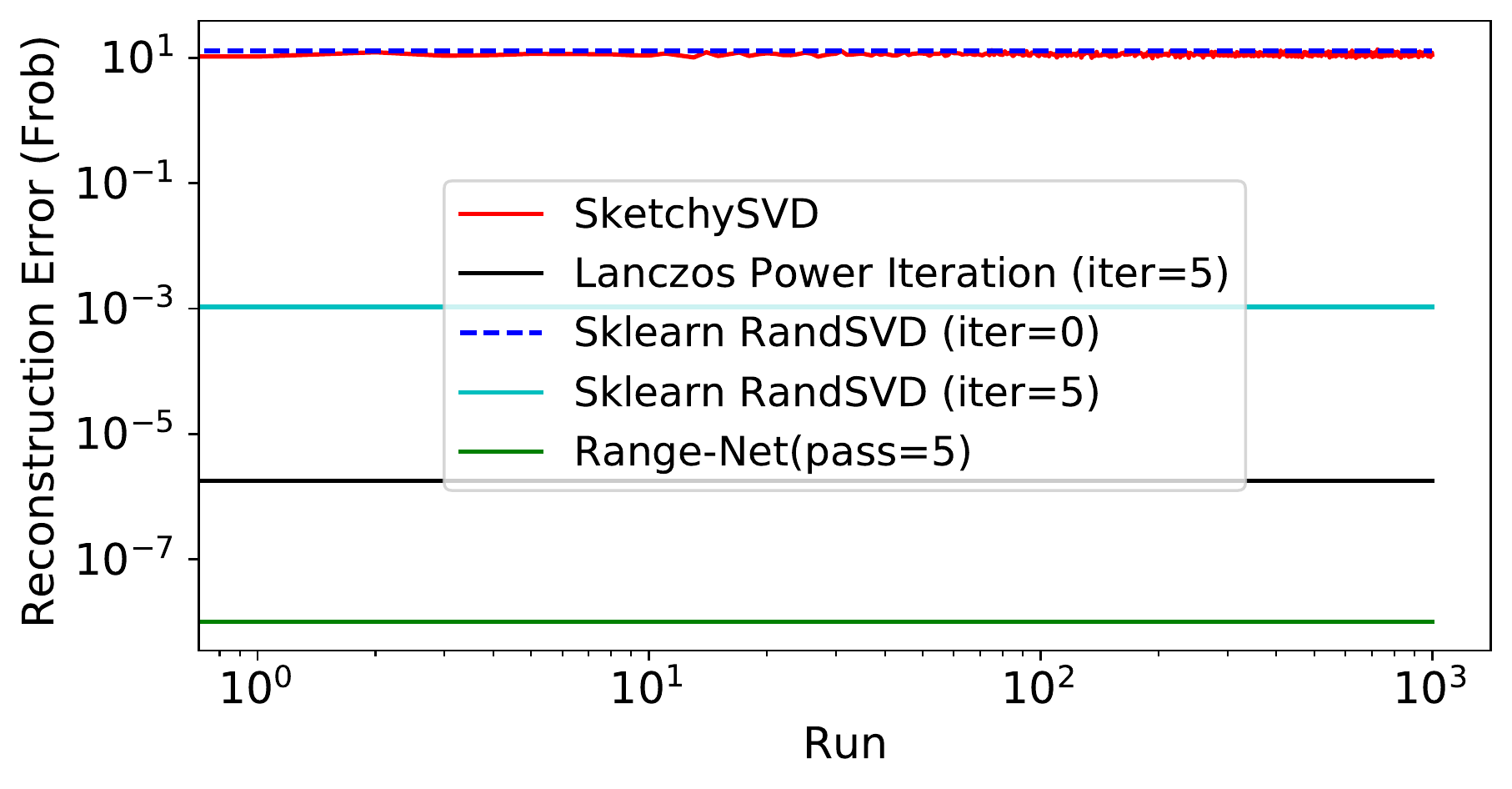}
    \caption{Comparison of reconstruction errors (rank $r=20$) for  Range-Net and randomized SVD schemes (with and without power iterations) for Parrot image over 1000 runs.}
    \label{fig:err_trend_parrot}
\end{figure}

\textbf{Fig. \ref{fig:err_trend_parrot}} shows the reconstruction errors in Frobenius norm for SketchySVD \cite{tropp2019streaming} (red line), Block Lanczos with Power Iteration \cite{musco2015stronger} (black line), Sklearn's randomized SVD \cite{skrandsvd} implementation \cite{halko2011finding} with (solid cyan line) and without power iteration (dashed blue line), and Range-Net (green line) over $1000$ runs on the Parrot image data. This shows that in order to gain lower reconstruction errors a power iteration is necessary that quickly becomes intangible in a big-data setting. Further, note that in \textbf{Fig. \ref{fig:err_trend_synth}} the expected error (upper bound) over multiple runs of Randomized SVD algorithms does not contract (reduce). 

\subsection{Navier-Stokes Simulated Data (SVD)}

For our next numerical experiment, we rely upon synthetic data generated using a Navier Stokes flow simulator for an incompressible fluid. The flow data is available on tensor-product grid on two-spatial and one temporal dimensions of size $(w,h,t) = (100 \times 50 \times 200)$. For each point on the grid, velocity vector values are available in both $x$ and $y$ spatial dimensions for $200$ time instances. The fully-developed, flow pattern exhibits a periodicity in the time dimension at approximately every $\sim 60$ time step that can be identified using SVD as characteristic modes. The data is therefore reshaped into a spatial vector for each time instance resulting in a spatio-temporal matrix $X \in \mathbb{R}^{5000 \times 200}$. For comparison purposes, we use only the x-direction stream velocity.

\begin{figure}[ht]
    \centering
    \includegraphics[width=0.8\linewidth]{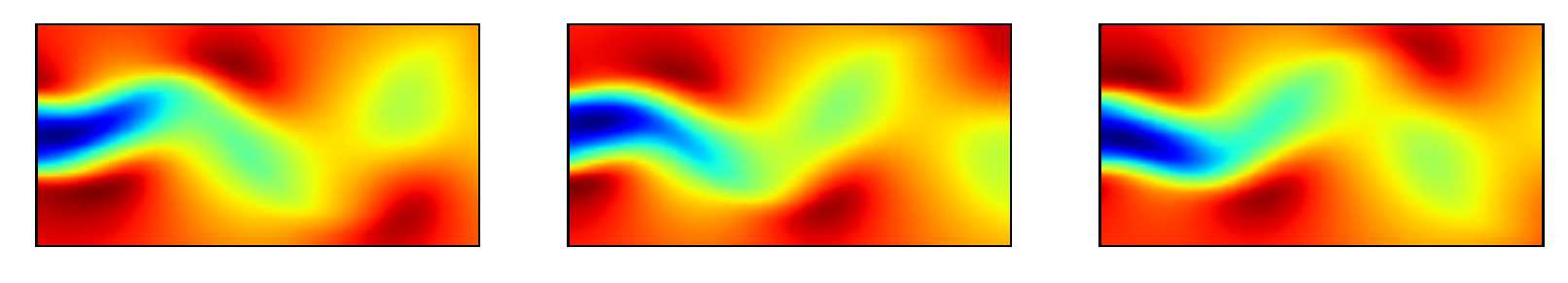}
    \caption{From left to right: x-direction stream velocity at times t = 0, 100, and 200}
    \label{fig:ns-flow}
\end{figure}

\begin{figure}[ht]
    \centering
    \begin{subfigure}{.22\linewidth}
      \centering
      \includegraphics[width=\linewidth]{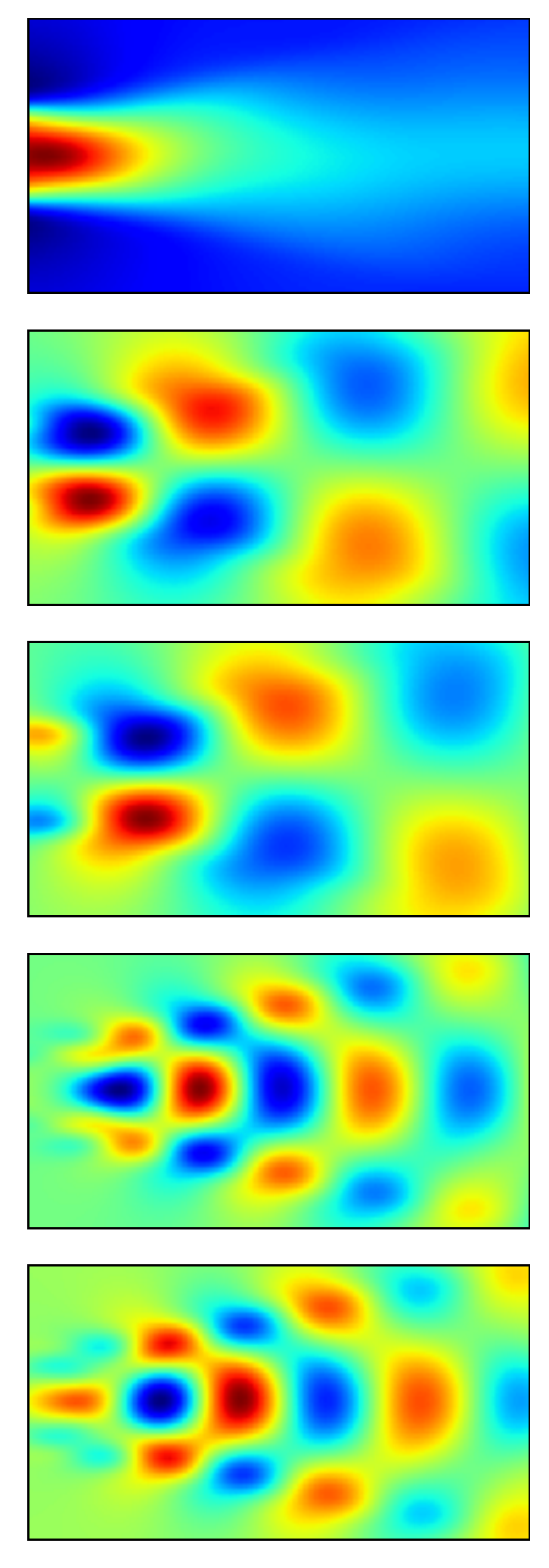}
      \caption{SketchySVD}
    \end{subfigure}
    \qquad
    \begin{subfigure}{.22\linewidth}
      \centering
      \includegraphics[width=\linewidth]{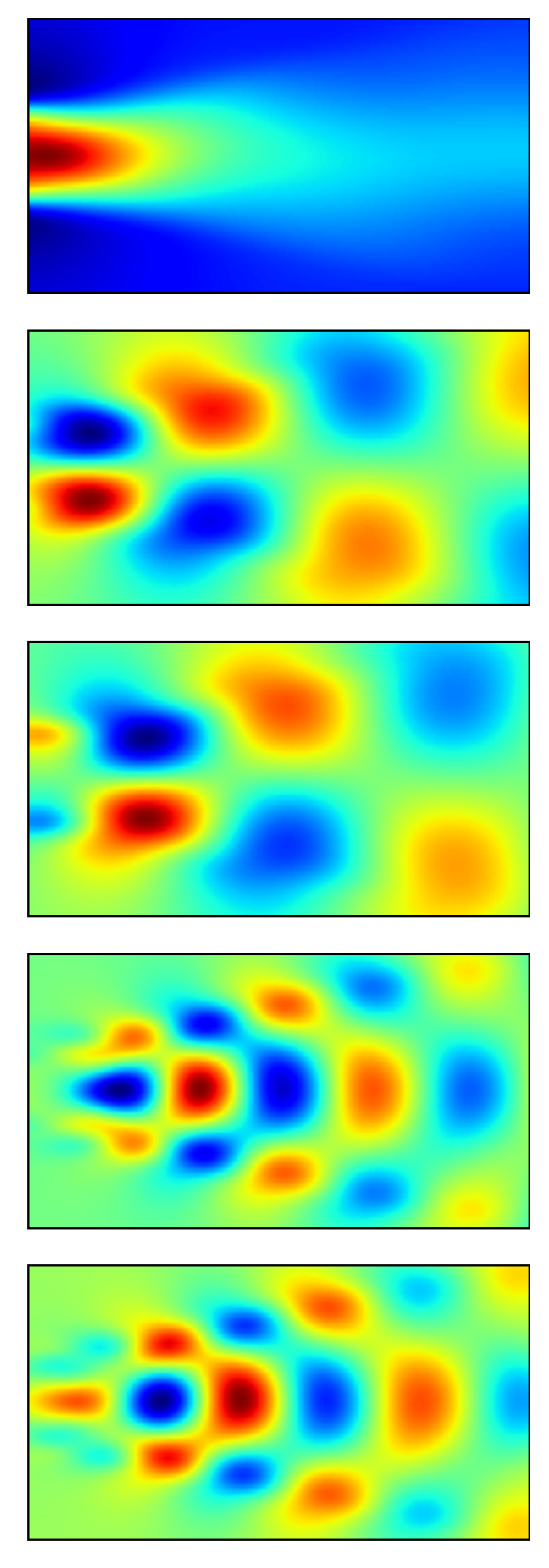}
      \caption{Truth}
    \end{subfigure}
    \qquad
    \begin{subfigure}{.22\linewidth}
      \centering
      \includegraphics[width=\linewidth]{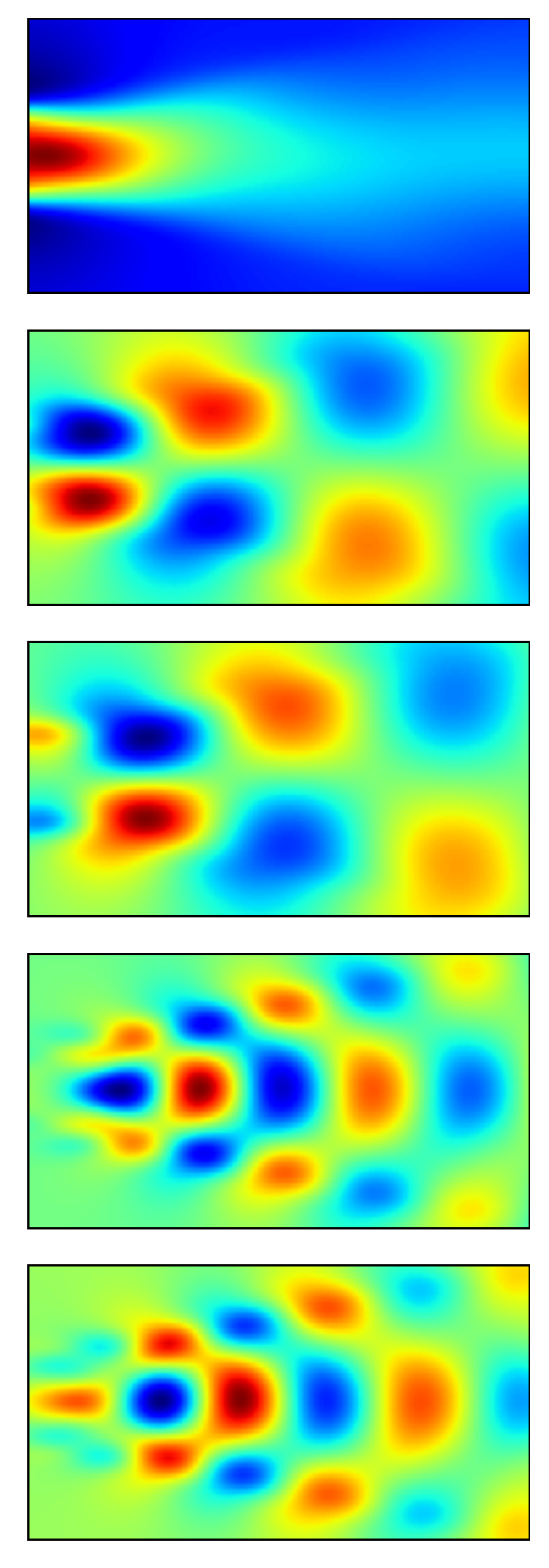}
      \caption{Range-Net}
    \end{subfigure}
    \caption{Reshaped $U$ indicative of dynamic modes, corresponding to top six right singular vectors for $r=5$. %The dynamic modes of SketchySVD has visible errors at indices $2,4,5$. Our proposed method does not have such artifacts.}
    }
    \label{fig:ns-vec}
\end{figure}

\begin{figure}[ht]
    \centering
    \begin{subfigure}{.54\linewidth}
      \centering
      \includegraphics[width=0.8\linewidth]{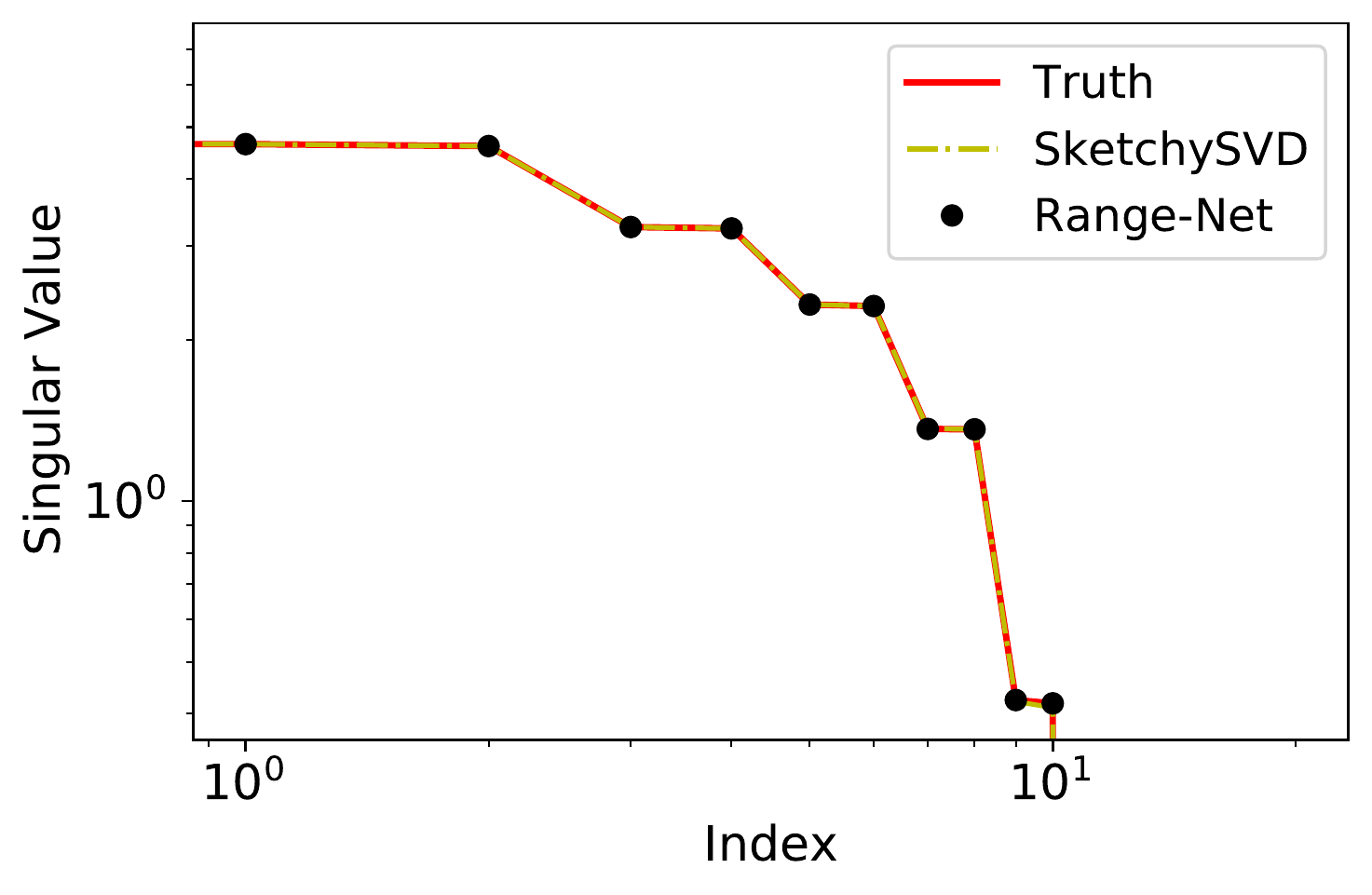}
      \caption{Singular Values}
    \end{subfigure}
    \begin{subfigure}{.39\linewidth}
      \centering
      \includegraphics[width=0.8\linewidth]{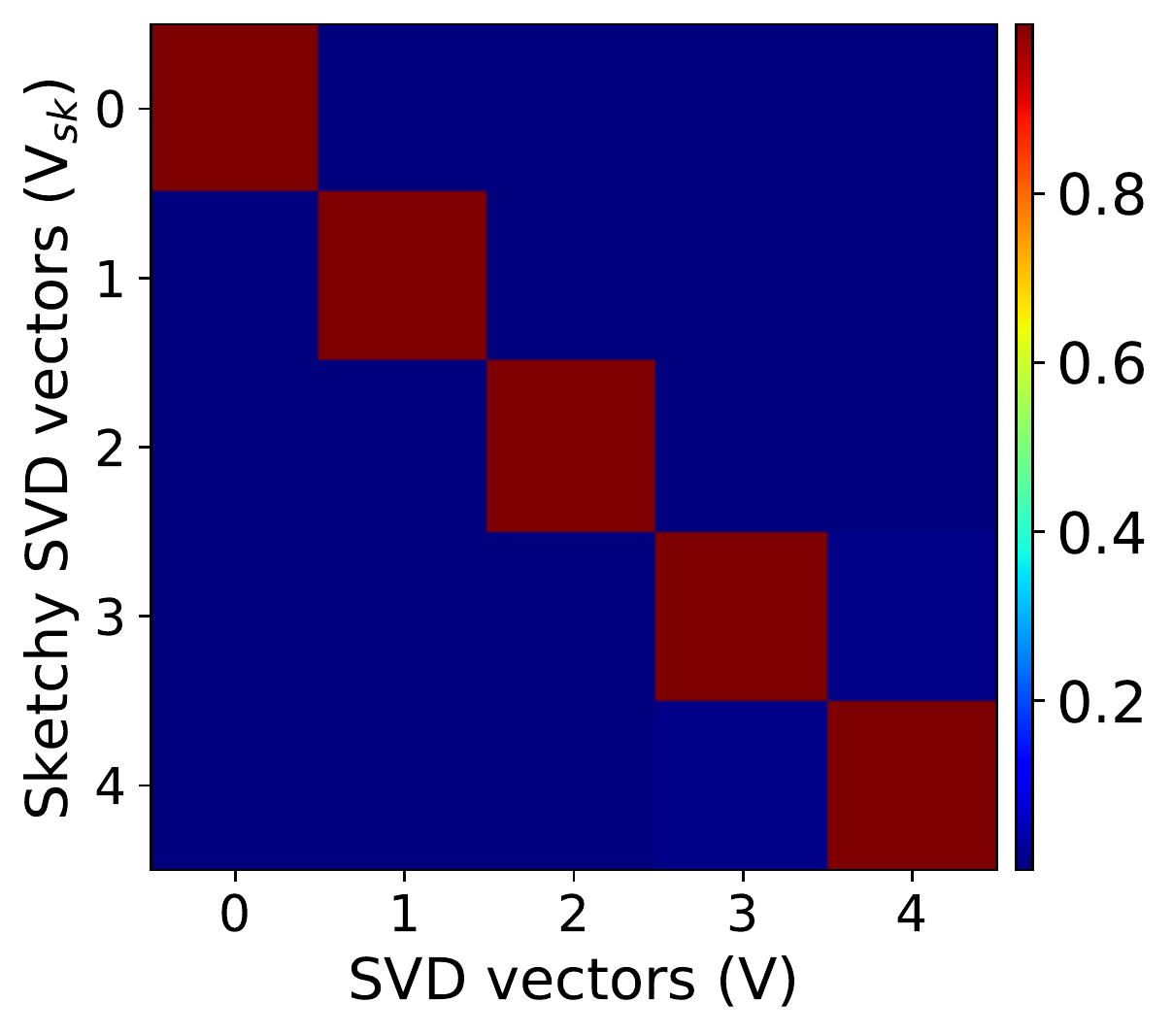}
      \caption{Correlation}
    \end{subfigure}
    \caption{Sketchy SVD and Range-Net (a) extracted singular values and (b) cross correlation between estimated and true (conventional SVD) right singular vectors for $r=5$ on the Navier-Stokes data.}
    \label{fig:ns-spec}
\end{figure}

\textbf{Fig. \ref{fig:ns-flow}} shows the evolution of the stream velocity over three time instances and the inherent time-periodic nature of the data. Notice the central-left region of primary flow across all three images. \textbf{Fig. \ref{fig:ns-vec}} shows the reshaped $U$ vectors for SketchySVD, conventional SVD and Range-Net. Notice that the images corresponding to the first left singular mode (also called dynamic modes) captures a notion of the primary flow in the left-center part. The second one captures spatial variations of the flow as time progresses. For all the three methods, all the modes have similar solution visually. Note that Range-Net computes a rank $r=5$ approximation without oversampling, wherein SketchySVD relies upon memory intensive sketchy projections of ranks $k=21,s=43$ to arrive at the solution. Overall, for this low-rank dataset both SketchySVD and Range-Net perform reasonably well qualitatively looking at the features in \textbf{Fig. \ref{fig:ns-vec}} and the singular value spectrum in \textbf{Fig. \ref{fig:ns-spec}}, because this synthetic dataset is extemely low-rank ($f=10$). 

\begin{remark}
For a given data matrix $X \in \mathbb{R}^{m\times n}$ of rank-$f$ ($f \leq \min(m,n)$) SketchySVD generates low error approximations if the oversampled rank $k$ is such that $k\geq f$. For full rank tall skinny matrices, this implies that $k\geq \min(m,n)$. For cases where $k< f$ (full rank or otherwise), SketchySVD accrues large approximation errors resulting in incorrect SVD factors.
\end{remark}

From a use-case point of view, randomized SVD generates low-error factors for full-rank, tall skinny matrices ($X \in \mathbb{R}^{m\times n}$) only when the oversampled rank $k\geq n$. This poses a serious limitation for all applications where this requirement is not met and consequently randomized SVD algorithms accrue large approximation errors in SVD factors as shown in the Sandy Big Data case study below and \textbf{Appendix \ref{app:sandy_low}}. Note that given a big data matrix $X$ determining the rank $f$ of $X$ is unknown and therefore selecting an oversampled rank $k$ such that $k \geq f$ is impractical in such cases. 

\begin{remark}
Note that once the SVD factors are extracted SketchySVD cannot be independently verified without performing a full SVD. In contrast, Range-Net is independently verifiable since the stage-1 of Range-Net cannot return orthonormal vectors if the vectors do not span a rank-$r$ subspace of a given data matrix $X$.
\end{remark}

\subsection{Dimension Reduction: MNIST (Eigen / PCA)} \label{sec:mnist}

Principal Component Analysis (PCA) is a special variant of Eigen decomposition, where the samples are mean corrected before constructing a feature covariance matrix followed by Eigen decomposition. Note that Range-Net does not require construction of the feature covariance matrix and can directly extract the eigenvectors and values without any modification. This is due to the fact that the Eigenvalues are the square of the singular values for any non-square data matrix $X \in \mathbb{R}^{m \times n}$ where right singular vectors are exactly the same as the eigenvectors.  

The MNIST dataset consists of handwritten digits images of size $28 \times 28$ with 60k images in the training dataset. To mold the dataset, we reshape each image into a $784$-dim vector to obtain the data matrix $X \in \mathbb{R}^{60000 \times 784}$, as a \textit{tall skinny matrix}. In a streaming setting, the mean feature vector computation requires one pass over the data matrix. This can be subsequently used during the network training (Stage-1) to mean correct streamed input vectors.

\begin{table}[ht]
\centering
  \caption{Error metrics and peak memory load for MNIST}
  \label{tab:errr-mnist}
  \begin{tabular}{c|ccc|c}
    \toprule
    rank & $err_{fr}$ & $err_{sp}$ & $\chi^2_{err}$  & SketchySVD | Range-Net \\
    \midrule
    20 & 0 & 0 & 0.012 & 0.47 GB | 0.52 MB \\ 
    50 & 0 & 0 & 0.025 & 1.18 GB | 1.31 MB \\ 
    100 & 1.12e-7 & 1.01e-7 & 0.052 & 2.38 GB | 2.78 MB\\ 
    200 & 2.36e-7 & 2.52e-7 & 0.071 & 4.83 GB | 6.19 MB \\
  \bottomrule
\end{tabular}
\end{table}

For this dataset, it is well known that $r=200$ captures around $97\%$ variance in the dataset. For SketchySVD, this results in projection matrices of ranks $k=4r+1=801$ and $s=2k+1=1602$. Since MNIST (tall and skinny matrix) data only has $n=784$ features, SketchySVD (Alg. \ref{alg:sketchy}) is almost equal if not more memory intensive than conventional SVD for such tall and skinny matrices. \textbf{Table. \ref{tab:errr-mnist}} shows the error metrics under different rank setting. Range-Net on the other hand with an exact memory requirement $(r(n+r))$ can handle much larger full rank tall and skinny matrices without incurring extraneous memory load. As discussed before, since this data matrix is tall and skinny ($60k \times 784$) we already know that for SketchySVD any rank-$r$ \st ($r\geq 196$) will result in the oversampling parameters $k \geq 784$ and $s \geq 1569$. SketchySVD will now extract lower-error SVD factors since the oversampled rank redundantly exceeds the feature dimension. 

\subsection{Scientific Computing: Sandy Big Data (SVD)}

\begin{figure}[h]
    \centering
    \includegraphics[width=0.6\linewidth]{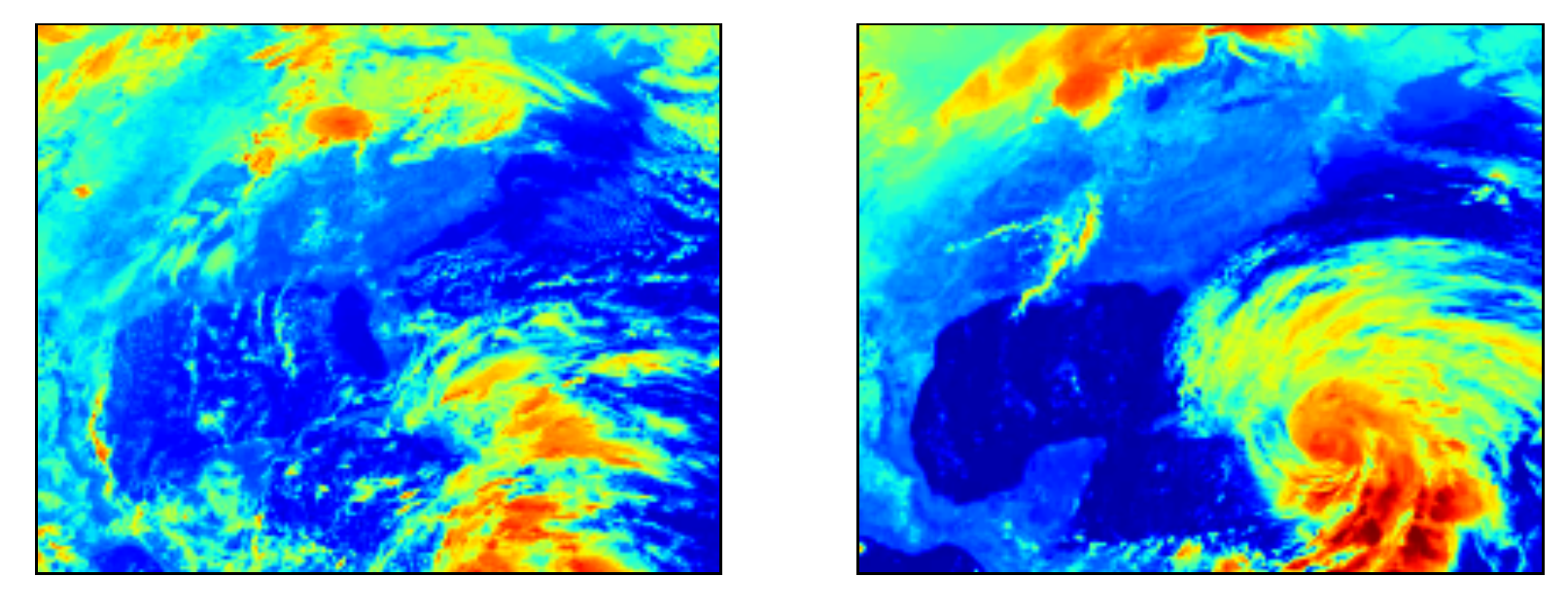}
    \caption{Satellite image captures of hurricane Sandy over the Atlantic ocean at $t = 0$ (left) and $t = 200$ minutes approximately (right).}
    \label{fig:sandy_snaps}
\end{figure}

Satellite data gathered by NASA for Hurricane Sandy over the Atlantic ocean represents the big data counter-part for scientific computations. The data-set is openly available\footnote{https://www.nasa.gov/mission\_pages/hurricanes/archives/2012/h2012\_Sandy.html} and comprises of RGB snapshots captured at approximately one-minute interval. The full data-set consists of $896 \times 719$ pixel images for $1208$ time-instances is of size \textbf{$24$ GB}. We chose this particular big data so that a conventional SVD can be performed on our machine ($16$ GB RAM) for benchmarking. Please note that this restriction is imposed by conventional SVD method due to its high main memory requirements. In contrast, our neural SVD solver can handle data sets that are orders of magnitude larger in size with the same hardware specification. \textbf{Fig. \ref{fig:sandy_snaps}} shows evolution of Hurricane Sandy for two time instances.

% We would also like to point out that although the data matrices $A_{50k \times 50k}$ for network graph problem appears big, the non-zero entries are highly sparse.
% In comparison, the data matrix $A_{160k \times 500}$ for Hurricane Sandy is still big since the matrix itself is dense. Therefore, a benchmark case cannot be setup using \textbf{irlba} sparse SVD solver in this case.

\begin{figure}[h]
    \centering
    \begin{subfigure}{.25\linewidth}
      \centering
      \includegraphics[width=\linewidth]{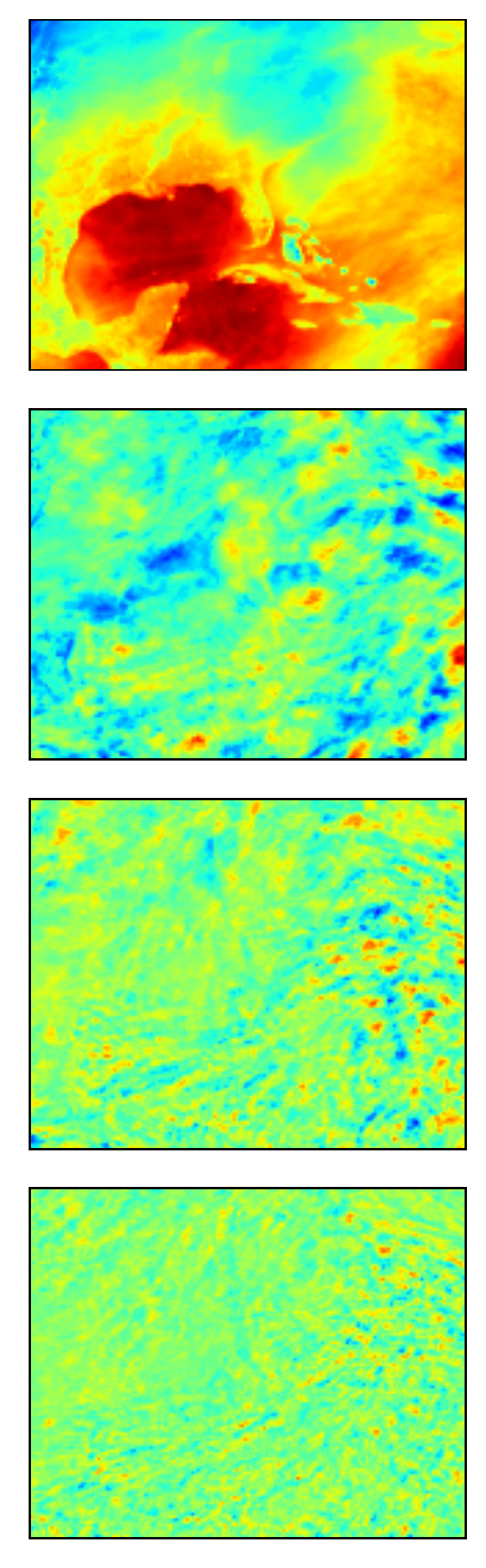}
      \caption{SketchySVD}
    \end{subfigure}
    \begin{subfigure}{.25\linewidth}
      \centering
      \includegraphics[width=\linewidth]{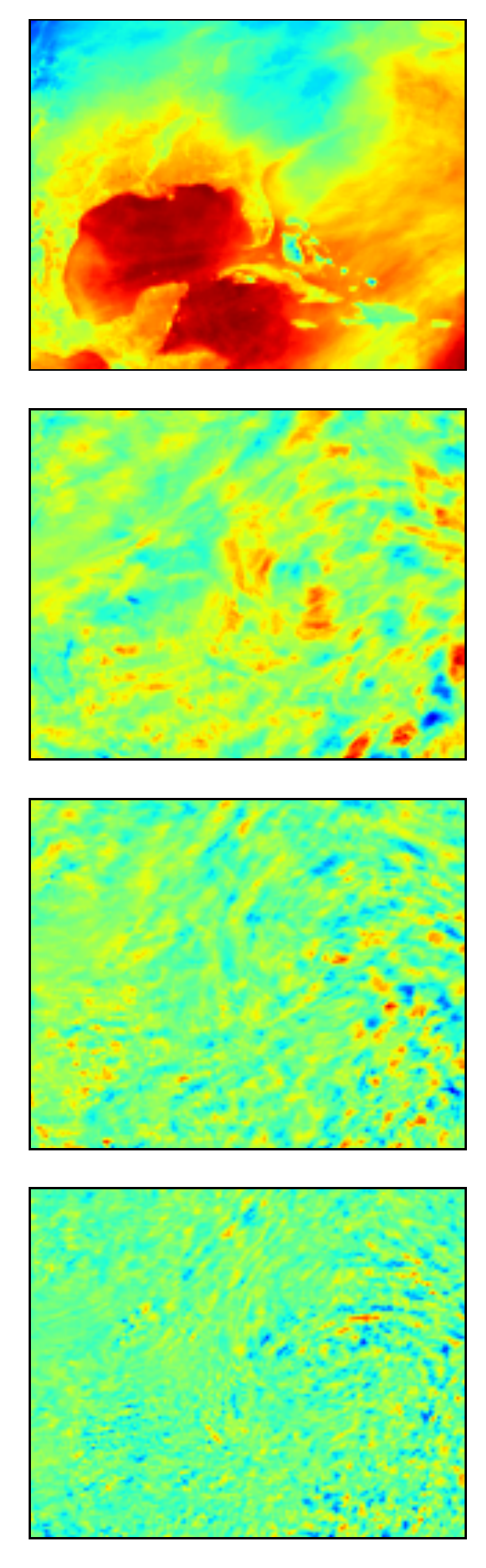}
      \caption{Truth}
    \end{subfigure}
    \begin{subfigure}{.25\linewidth}
      \centering
      \includegraphics[width=\linewidth]{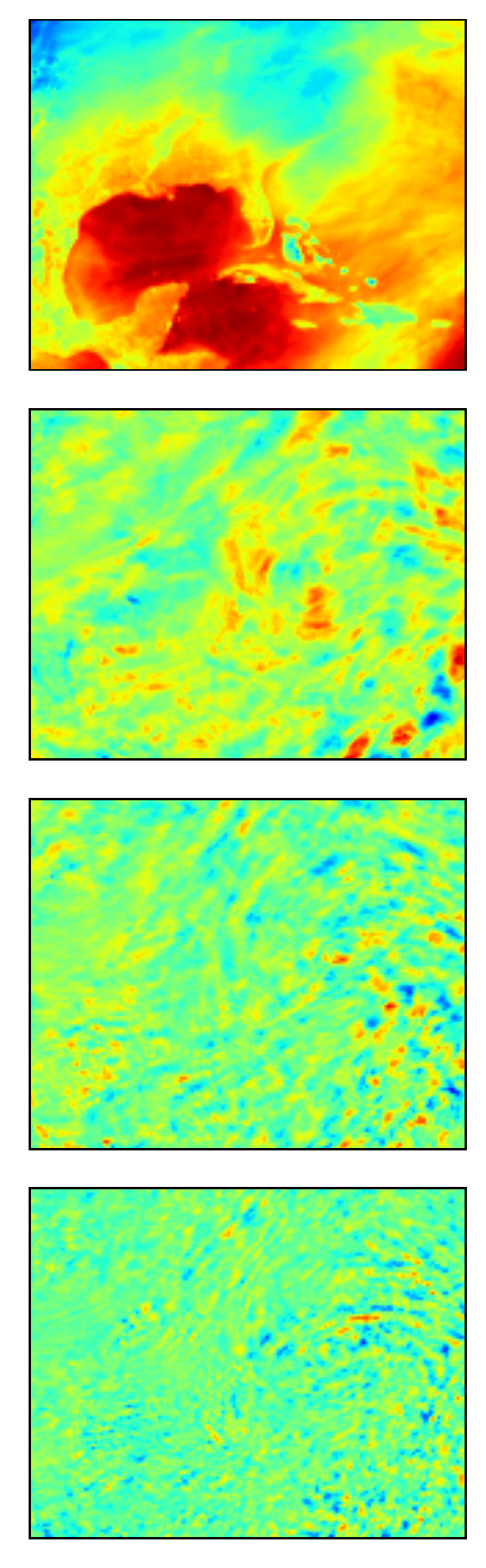}
      \caption{Range-Net}
    \end{subfigure}
    \caption{Reshaped $U_i$ indicative of dynamic modes, corresponding to $i=1,20,50,100$ for $r=100$ (oversampled rank $k= 401$ for SketchySVD. The dynamic mode approximation error stand out visually for SketchySVD for indices $20,50,100$. Our method does not have such artifacts.}
    \label{fig:sandy_sing}
\end{figure}

Range-Net can not only handle datasets of larger sizes similar to SketchySVD, but also ensures lower errors in approximating the SVD factors. For SketchySVD, the user cannot directly verify the approximation errors in the extracted factors without relying on a conventional SVD solver (infeasible for big data). On the other hand, Stage-1 of Range-Net generates near-orthonormal vectors (at GPU precision) only when the vectors $V_*$ spans the desired rank-$r$ subspace of a given data matrix (see \textbf{Lemma \ref{lem:1.2}}). The user can verify this numerically by computing the deviation of $V_*^TV_*$ from the rank-$r$ identity ($I_r$).

Similar to the Navier-Stokes simulation data, \textbf{Fig. \ref{fig:sandy_sing}} shows three dynamic modes corresponding to rank $1,20,50,100$ singular values. As shown, our results are in good agreement with conventional SVD whereas, Sketchy SVD shows substantial deviations after the first 50 dynamic modes. Here, we point out that accuracy is a matter of special concern in scientific computations. Any compression that results in substantial loss of information or obscuring an otherwise identifiable feature in recorded observations directly culls our capacity to make scientific improvements. Consequently, any exploratory data analysis, however big or small, must accurately identify the underlying features. Our algorithm achieves the the lower bound on tail-energy given by EYM theorem  ensures an accurate resolution in this big data setting. Please note that increasing the sensor resolution implies that we are interested in exploring and understanding the high-frequency features (lower singular values) of the data.

\begin{table}[ht]
\centering
  \caption{Error metrics and peak memory load for Sandy}
  \label{tab:errr-sandy}
  \begin{tabular}{c|ccc|c}
    \toprule
    rank & $err_{fr}$ & $err_{sp}$ & $\chi^2_{err}$  & SketchySVD | Range-Net \\
    \midrule
    10 & 0 & 0 & 0.011 & 2.56 GB | 0.39 MB \\
    50 & 0 & 0 & 0.018 & 12.48 GB | 2.01 MB \\
    100 & 1.12e-7 & 1.24e-7 & 0.021 & 24.91 GB | 4.19 MB\\
  \bottomrule
\end{tabular}
\end{table}

\begin{figure}[h]
    \centering
    \begin{subfigure}{.3\linewidth}
      \centering
      \includegraphics[width=\linewidth]{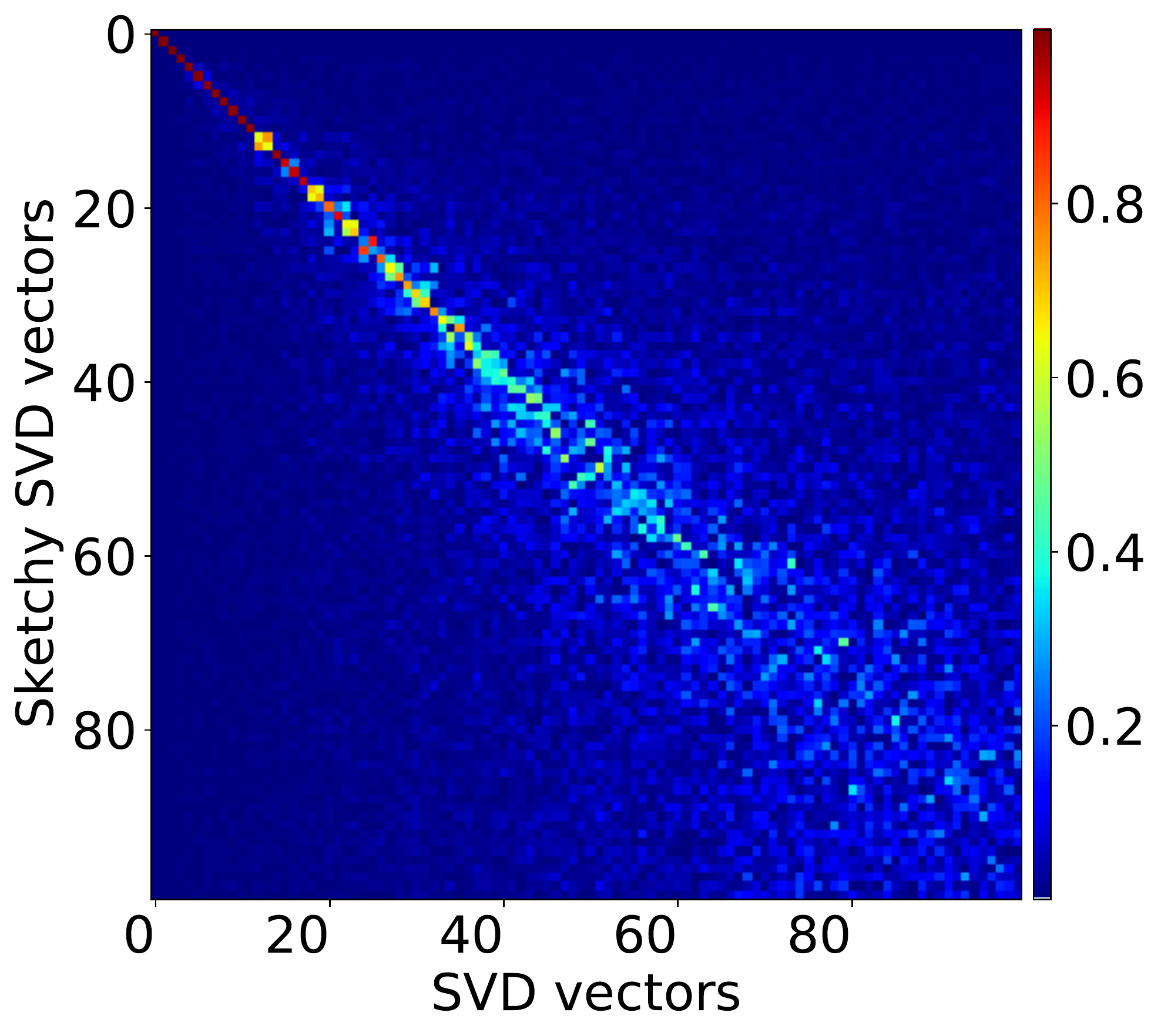}
      \caption{SketchySVD}
    \end{subfigure}
    \qquad
    \begin{subfigure}{.3\linewidth}
      \centering
      \includegraphics[width=\linewidth]{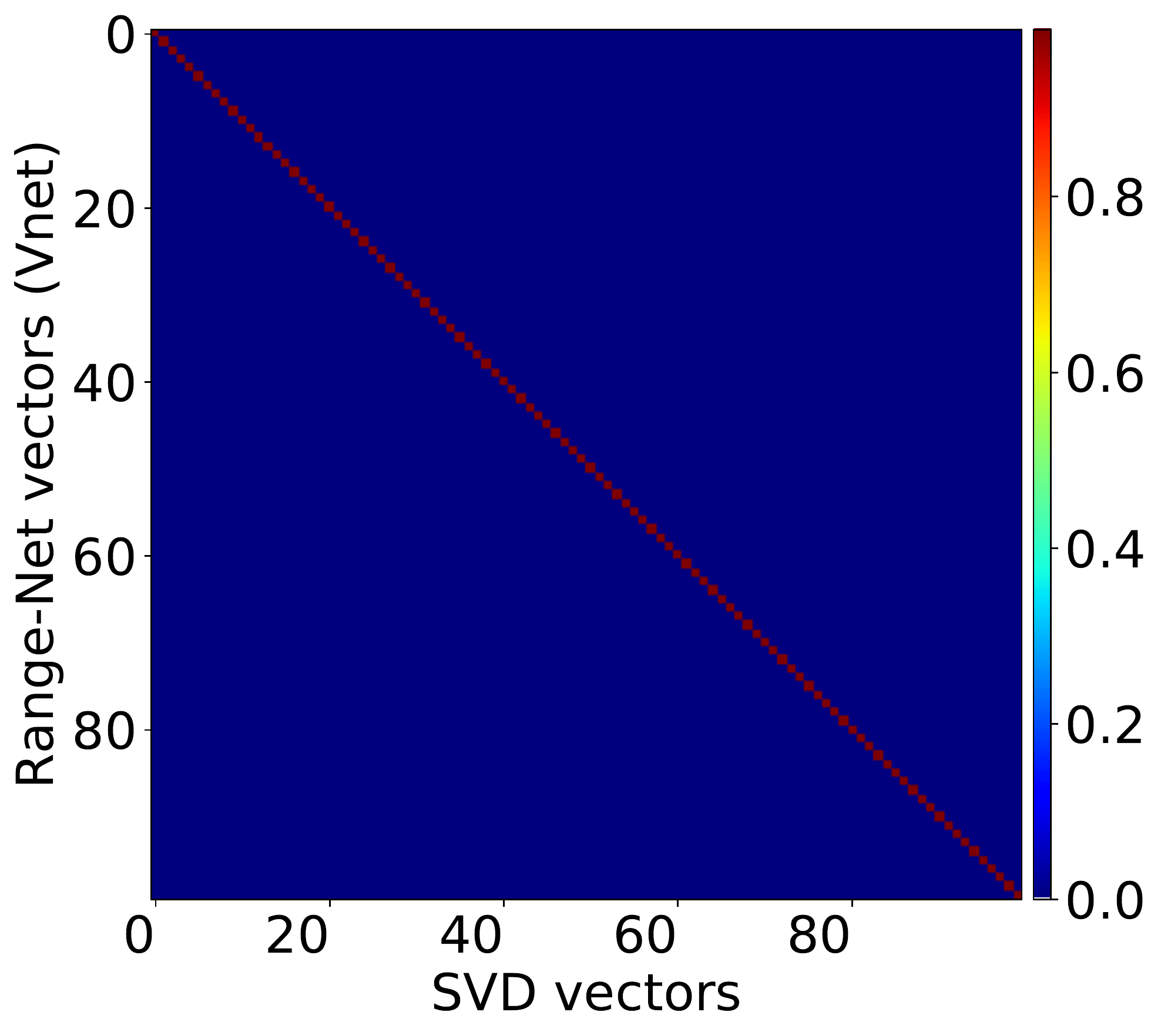}
      \caption{Range-Net}
    \end{subfigure}
    \caption{Cross-correlation between extracted and true (conventional SVD) right singular vectors for (a) SketchySVD and (b) Range-Net for a rank $r = 100$ approximation. SketchySVD deviates substantially after index 10 (although sketching at sizes $k=401$ and $s=803$) while Range-Net is in good agreement for all the 100 indices.}
    \label{fig:sandy_corr}
\end{figure}

\begin{figure}[h]
    \centering
    \begin{subfigure}{.4\linewidth}
      \centering
      \includegraphics[width=\linewidth]{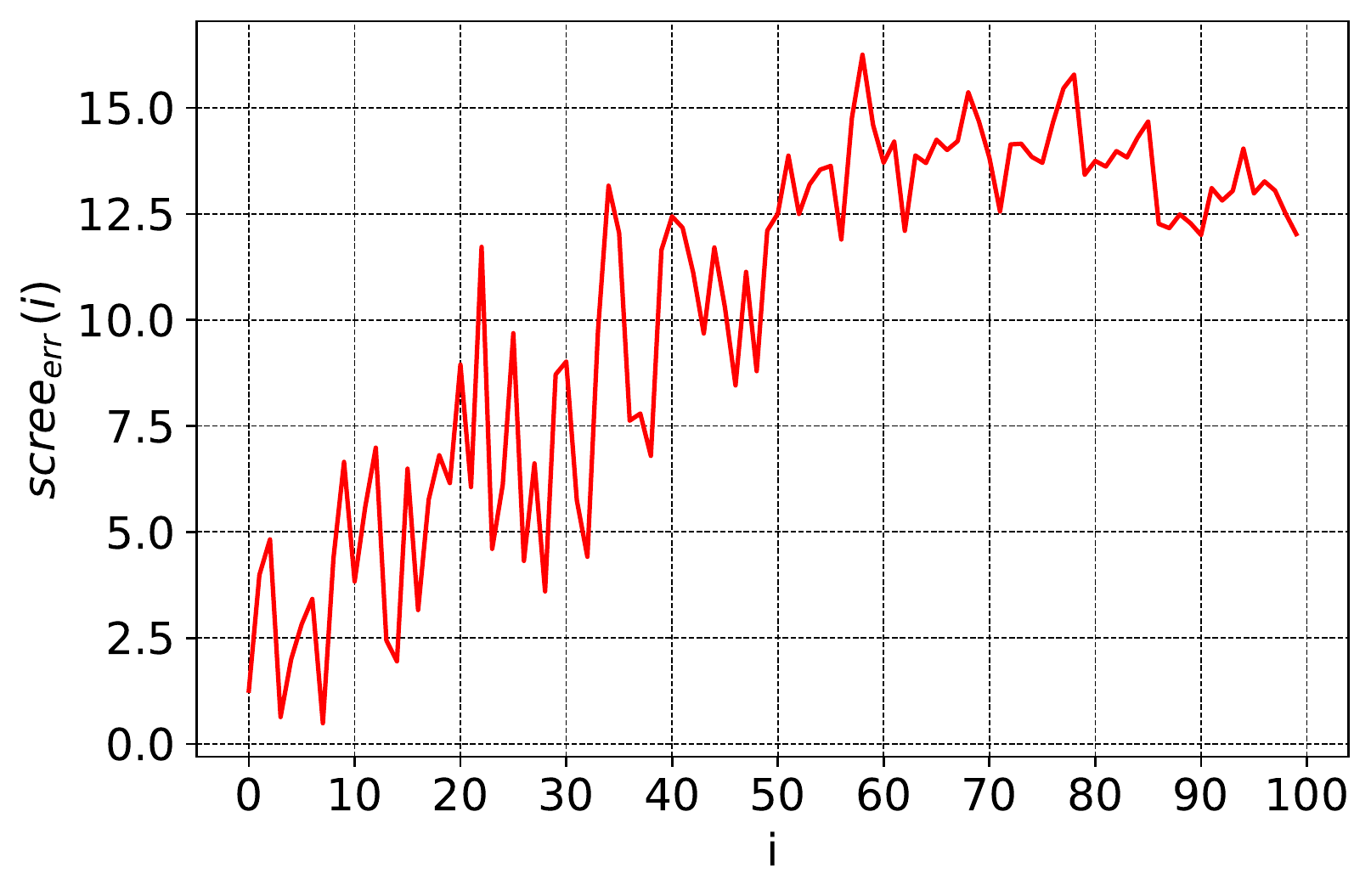}
      \caption{SketchySVD}
    \end{subfigure}
    \qquad
    \begin{subfigure}{.4\linewidth}
      \centering
      \includegraphics[width=\linewidth]{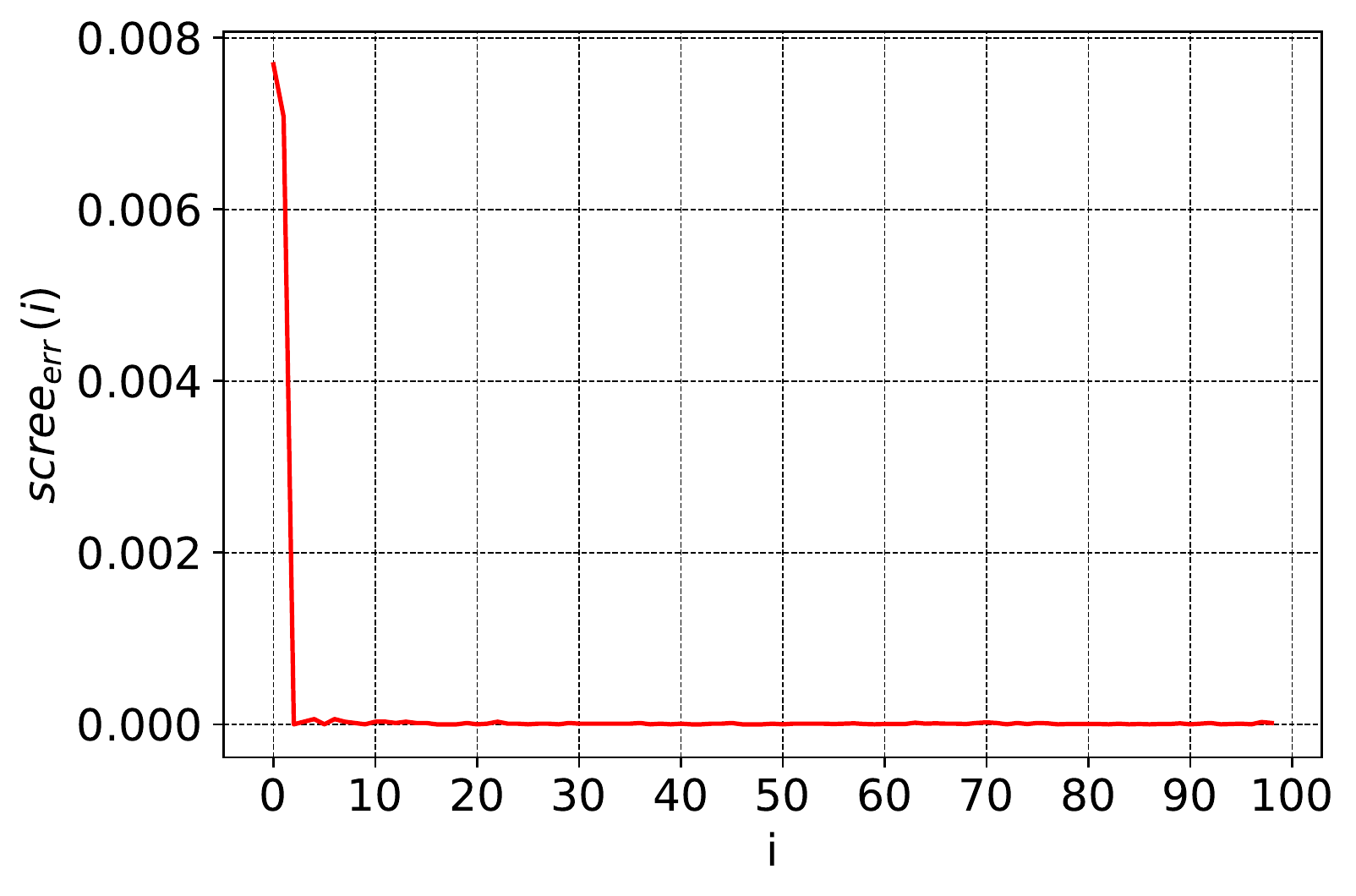}
      \caption{Range-Net}
    \end{subfigure}
    \caption{Scree-error in singular values for (a) SketchySVD and (b) Range-Net where a conventional SVD is used as the baseline in scree-error metric.}
    \label{fig:sandy_scree}
\end{figure}

\textbf{Figs. \ref{fig:sandy_corr}} and \textbf{\ref{fig:sandy_scree}} show the cross correlation of the right singular vectors and scree-error in the corresponding singular values extracted by SketchySVD and Range-Net. Note that for a rank-$100$ approximation, SketchySVD extracted right singular vectors start deviating after rank-$10$ as shown in \textbf{Figs. \ref{fig:sandy_corr}} while the singular values deviate quite substantially from rank-$1$. Range-Net on the other hand is in excellent agreement with the right singular vectors and values for all desired $100$ indices. \textbf{Tab. \ref{tab:errr-sandy}} shows the error metrics for Range-Net with a comparison of peak main-memory load between SketchySVD and Range-Net for ranks $r=[10,50,100]$. 

As mentioned before, Randomized SVD extracted factors deviate quite substantially when the user specified rank $r$ is such that the oversampled rank $k$ is much lower than the unknown rank $f$ of a given data matrix. The reader is referred to \textbf{Appendix \ref{app:sandy_low}} where a rank-$10$ approximation is extracted. Here, SketchySVD deviates quite substantially after rank-$3$ while Range-Net still remains in excellent agreement with the baseline singular vectors and values. From an application point of view, a lower rank is desired in graph-datasets where the second singular value and associated vector is used for community detection. The reader is referred to \textbf{Appendix \ref{app:graph}} where Range-Net is used on sparse adjacency matrices from a number of publicly available graph-datasets . 

\subsection{Storage Complexity Analysis}

To get an estimate of the memory efficiency of Range-Net, let us consider the peak main memory (RAM) requirement for the compuation of SVD factors. Range-Net has two layers in succession, one corresponding to the low-rank projector $\tilde{V}^{r \times n}$ and the rotation matrix $\Theta^{r \times r}$. For Sketchy SVD, the peak memory load occurs during the construction of a core matrix $C^{s \times s}$ (see \textbf{Alg. \ref{alg:sketchy}}). This requires that the two projection matrices $\Phi^{m \times s}, \Psi^{n \times s}$, one projected data matrix $Z^{s \times s}$, and two rank-$k$ decomposition $Q^{m \times k}, P^{n \times k}$ and the core matrix $C^{s \times s}$,  be present in the memory simultaneously. Thus, the overall memory efficiency factor $(s_{eff})$ between Sketchy SVD and Range-Net for a rank-$r$ approximation is:

\begin{align*}
    s_{eff} &= \frac{\textrm{SketchySVD (Peak Load)}}{\textrm{Range-Net}} 
    =\frac{ns+ms+s^2+mk+s^2+nk}{rn+r^2}\\
    &= \frac{(m+n)(k+s) + 2s^2}{r(n+r)} 
    % &\approx \frac{3(m+n)(4r + 8r) + (8r)^2 + (4r)^2}{r(n+r)} \\
    \approx \frac{12(m+n) + 128r}{(n+r)} \quad (k=4r+1,s=2k+1) \\
    &\approx 7.67e2 \quad \textrm{for MNIST}(m=60k,n=784,r=200)
\end{align*}

\begin{figure}[h]
    \centering
    \includegraphics[width=0.4\linewidth]{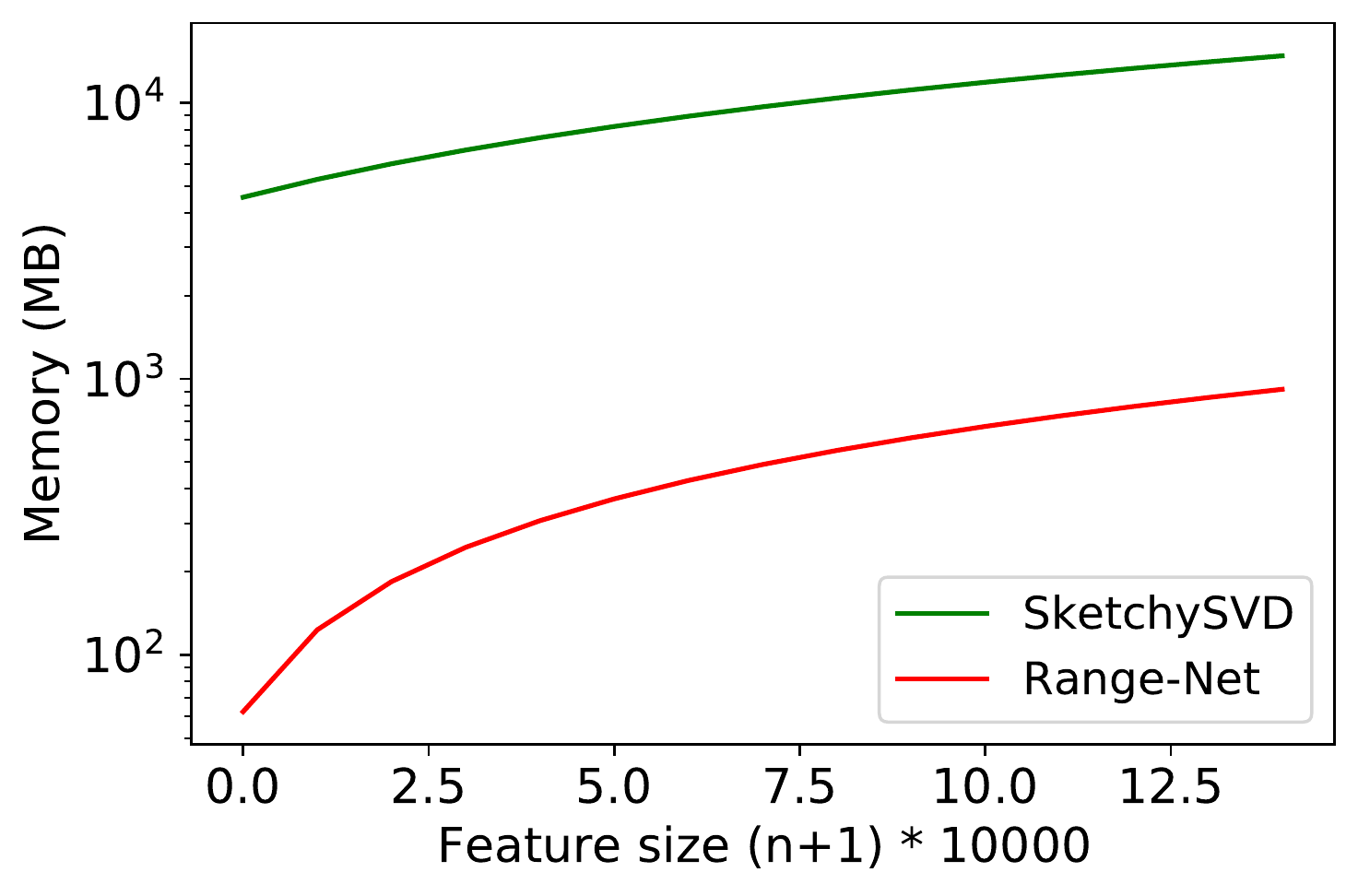}
    \caption{Peak Memory Requirement comparison of SketchySVD \vs Range-Net, on a synthetic matrix of size $m=50000$ and $r=200$, while $n$ varies between $[10k-150k]$.}
    \label{fig:mem_req}
\end{figure}

To validate the above ratio, we constructed a synthetic dataset of $m=50000$ rows and the number of columns were varied starting at $n=10000$ with increments of $10000$. The expected rank from the algorithms was held at $r=200$. \textbf{Fig. \ref{fig:mem_req}} shows the memory allocation (in Megabytes (MB)), where the storage efficiency of our method is evident. When $m=50000, n=150000$ and $r=200$, SketchySVD has a peak memory consumption of $14$GB due to oversampling parameters of $k=801,s=1603$, while Range-Net only requires $916$MB.

\section{Conclusion}

We present Range-Net as a low-weight, high-precision, fully interpretable neural SVD solver for big data applications that is independently verifiable without performing a full SVD. We show that our solution approach achieves lower errors metrics for the extracted singular vectors and values compared to Randomized SVD methods. A discussion is also provided on the limiting assumptions and practical consequences of using Randomized SVD schemes for big data applications. Our numerical experiments on real and synthetic datasets confirm that Range-Net achieves the theoretical lower bound on the tail energy given by EYM theorem. We also verify that our network minimization problems converges to this tail energy bound in Frobenius norm at machine precision. A number of Big Data problems are considered, where SVD or Eigen decompositions are required, that demonstrate the applicability of Range-Net to large scale practical datasets. A fair comparison is also provided against a state of the art randomized, streaming SVD algorithm with conventional SVD solution as the baseline for computational benchmarking and verification. Finally, Range-Net is independently verifiable without requiring a full SVD as opposed to Randomized SVD algorithms where upper bounds on tail-energy errors are subjective to the dataset.

%\bibliographystyle{ACM-Reference-Format}
%\bibliography{References}


%
\appendix

\section{Theoretical Guarantees} \label{app:thm}

\subsection{Preliminaries}\label{app:prel}

The Frobenius norm of a matrix $A$ is given by,
\begin{align*}
    \|A\|_{F} = \left(\sum_{i}\sum_{j}a_{ij}^2\right)^{\frac{1}{2}} = \left(\Tr(A^{T}A)\right)^{\frac{1}{2}} = \left(\Tr(AA^{T})\right)^{\frac{1}{2}} 
\end{align*}
Further, the Frobenius norm can be used to bound \cite{trefethen1997numerical} the norm of a matrix product as,
\begin{align*}
    \|AB\|_{F}\leq \|A\|_{F}\|B\|_{F}
\end{align*}
For a Frobenius norm we have that,
\begin{align*}
    \|A+B\|_{F} \leq \|A\|_{F} + \|B\|_{F}\\
    \|A-B\|_{F} \geq \big|\|A\|_{F} - \|B\|_{F}\big|
\end{align*}
Also the Frobenius norm of a rank-$f$ matrix $A$,
\begin{align*}
    \|A\|_{F} = \|\Sigma_{f}\|_{F}
\end{align*}
where $\Sigma_{f} = \mathrm{diag}(\sigma_{1},\sigma_{2},\cdots,\sigma_{f})$ and $\sigma_{i}$s are the $f$ non-zero, singular values of $A$.

Let $A$, $B$, $C$ be matrices such that the following matrix products are feasible. The cyclic property of the linear trace operator is,
\begin{align*}
    \Tr(ABC) = \Tr(BCA) = \Tr(CAB)
\end{align*}

\begin{definition}
The tail energy of an arbitrary matrix $B \in \mathbb{R}^{m\times n}$ with respect to a given matrix $X\in \mathbb{R}^{m\times n}$ equipped with a Frobenius norm is defined as,
\begin{align*}
    \tau = \|X-B\|_{F}
\end{align*}
\end{definition}

\begin{definition}
Let $r,f\in\mathbb{Z}^{+}$ be positive integers such that $0 < r \leq f$, then a rank-r truncation $X_{r}$ of a rank-f matrix $X$ is defined as,
\begin{align*}
    X_{r} = U_r\Sigma_r V^T_r = XV_{r}V_{r}^{T}
\end{align*}
where, $\Sigma_{R} = \mathrm{diag}(\sigma_{1},\sigma_{2},\cdots,\sigma_{r})$ and $\sigma_{i}$s are the top $r$ singular values of $X$ and $V_{r} = [v_{1}, v_{2},\cdots, v_{r}]$, and $U_{r} = [u_{1}, q_{2},\cdots, q_{r}]$ are matrices such that $v_{i}$s and $u_{i}$s are the corresponding right and left singular vectors, respectively.
\end{definition}

The relative tail energy of a rank-$r$ matrix $B_{r}$ with respect to a rank-$f$ matrix $X$ ($r\leq f$) is then defined as,
\begin{align*}
    \tau_{rel,r} = \frac{\|X-B_{r}\|_F}{\|X-X_r\|_F} - 1
\end{align*}

\begin{theorem}
Eckart-Young-Mirsky Theorem \cite{eckart1936approximation,mirsky1960symmetric}: Let $X\in \mathbb{R}^{m\times n}$ be a real, rank-f, matrix with $m \geq n$ with the singular value decomposition as $X = U\Sigma V^{T}$, where the orthonormal matrices $U,V$ contain the left and right singular vectors of $X$ and $\Sigma$ is a diagonal matrix of singular values. Then for an arbitrary rank-$r$, $r\leq f$ matrix $B_{r} \in \mathbb{R}^{m\times n}$,
\begin{align*}
    \|X-B_{r}\|_{F} \geq \|X-X_{r}\|_{F}
\end{align*}
where $X_{r} = U_{r}\Sigma_{r}V_{r}$ with $\Sigma_{r}$ is the diagonal matrix of the largest $r$ singular values and $U_{r},V_{r}$ are the corresponding left and right singular vector matrices.
\label{thm:orig_eym}
\end{theorem}

\subsection{Stage 1}\label{app:thm:stage1}

\begin{theorem}
For any $r,f\in \mathbb{Z}^{+}$, $0< r \leq f$, if the tail energy of a rank-f matrix $X \in \mathbb{R}^{m \times n}$, $f\leq min(m,n)$, with respect to an arbitrary rank-r matrix $B_{r}=X\tilde{V}_{r}\tilde{V}_{r}^{T}$ is bounded below by the tail energy of $X$ with respect to it's rank-r approximation $X_{r} = XV_{r}V_{r}^{T}$ as,
\begin{center}
    $\|X-B_{r}\|_{F} \geq \|X-X_{r}\|_{F}$
\end{center}
where, $V_{r} = span\{v_{1},v_{2},\cdots,v_{r}\}$ and $v_{i}$s are the right singular vectors corresponding to the largest $r$ singular values then the minimizer of $\underset{\tilde{V} \in \mathbb{R}^{(n\times r)}}{\arg \min}\|X-X\tilde{V}\tilde{V}^{T}\|_{F}$ is $V_{*}$ such that $V_{*}V_{*}^{T} = V_{r}V_{r}^{T}$. 
\label{thm:eym_mod}
\end{theorem}

From \textbf{Theorem \ref{thm:orig_eym}} we have,

\begin{equation}
\|X-B_{r}\|_{F} - \|X-X_{r}\|_{F} \geq 0
\label{eqn:eym_mod}
\end{equation}

Let $V_{r} = \{v_{1}, v_{2},\ldots,v_{r}\}$ be the top-r, right-singular vectors of $X$ corresponding to the largest singular values,

\begin{equation}
X_{r} = XV_{r}V_{r}^{T} =  U\Sigma V^{T}V_{r}V_{r}^{T} = U\Sigma V_{r}^{T}   
\end{equation}

Also let $B_{r} = X\tilde{V}_{r}\tilde{V}_{r}^{T}$ where $\tilde{V}_{r}$ is an arbitrary rank-r matrix. From triangle inequality we have that,

\begin{equation}
\|X(V_{r}V_{r}^{T}-\tilde{V}_{r}\tilde{V}_{r}^{T})\|_{F} \geq \|X(I_{n}-\tilde{V}_{r}\tilde{V}_{r}^{T})\|_{F} -  \|X(I_{n}-V_{r}V_{r}^{T})\|_{F} 
\label{eqn:tri}
\end{equation}

Combining Eq. \ref{eqn:eym_mod} and Eq. \ref{eqn:tri} we get,

\begin{equation}
    \|X(V_{r}V_{r}^{T}-\tilde{V}_{r}\tilde{V}_{r}^{T})\|_{F} \geq 0
\end{equation}

Additionally,
\begin{equation}
     \|X\|_{F}\|V_{r}V_{r}^{T}-\tilde{V}_{r}\tilde{V}_{r}^{T}\|_{F} \geq \|X(V_{r}V_{r}^{T}-\tilde{V}_{r}\tilde{V}_{r}^{T})\|_{F}
\end{equation}

Using the above two inequalities we arrive at,
\begin{equation}
    \|X\|_{F}\|V_{r}V_{r}^{T}-\tilde{V}_{r}\tilde{V}_{r}^{T}\|_{F} \geq 0
\end{equation}
Since $\|X\|_{F}>0$,  equality is achieved when $\tilde{V}_{r}\tilde{V}_{r}^{T} = V_{*}V_{*}^{T}=  V_{r}V_{r}^{T}$. In other words, $span\{V_{r}\} = span\{V_{*}\}$ since $(V_{*}\Theta_{r})(V_{*}\Theta_{r})^{T} = V_{*}V_{*}^{T}$ for any rank-r, real valued, unitary matrix $\Theta_{r}$ spanning the top rank-r subspace of $X$. 

\begin{remark}
Theorem \ref{thm:eym_mod} also implies that any matrix $\tilde{V}_{r}\tilde{V}_{r}^{T}$ that does not span the same rank-r subspace of $X$ as $V_{r}V_{r}^{T}$ will result in a higher tail-energy than given by the EYM tail-energy bound equipped with a Frobenius norm.
\end{remark}

\begin{lemma}
	If $V_{r}^{T}V_{r} = I_{r}$ and $V_{r}V_{r}^{T} = V_{*}V_{*}^{T}$ then $V_{*}^{T}V_{*} = I_{r}$.
	\label{lem:1.1}
\end{lemma}

Let us assume $V_{r}^{T}V_{r} = I_{r}$ then,
\begin{align*}
	\|V_{r}^{T}V_{r} - I_{r}\|_{F}^{2} = 0\\
	\Tr\left(V_{r}^{T}V_{r}V_{r}^{T}V_{r}\right) + \Tr\left(I_{r}\right) - 2\Tr\left(V_{r}^{T}V_{r}\right) = 0
\end{align*}
Using the cyclic property of the trace operator we have,
\begin{align*}
	\Tr\left(V_{r}V_{r}^{T}V_{r}V_{r}^{T}\right) + \Tr\left(I_{r}\right) - 2\Tr\left(V_{r}V_{r}^{T}\right) = 0
\end{align*}
Using \textbf{Theorem \ref{thm:eym_mod}}, $ V_{*}V_{*}^{T} = V_{r} V_{r}^{T}$,
\begin{align*}
	\Tr\left(V_{*} V_{*}^{T}V_{*}V_{*}^{T}\right) + \Tr\left(I_{r}\right) - 2\Tr\left(V_{*}V_{*}^{T}\right) = 0
\end{align*}
Again using the cyclic property of the trace operator we now get,
\begin{align*}
	\Tr\left(V_{*}^{T}V_{*} V_{*}^{T}V_{*}\right) + \Tr\left(I_{r}\right) - 2\Tr\left(V_{*}^{T}V_{*}\right) = 0
\end{align*}
Hence, $V_{*}^{T}V_{*} = I_{r}$. This shows that following \textbf{Theorem \ref{thm:eym_mod}}, the matrix $V_{*}$ comprises of orthonormal column vectors spanning the same top rank-r subspace of $X$ as the orthonormal column vectors of $V_{r}$. 

\begin{lemma}
	If $X \in \mathbb{R}^{m\times n}$ is a rank $f$ matrix, then for any rank $r > f$, where  $\{r,f\} \leq min(m,n)$, if $V_{*}^{T}V_{*} = I_{r}$ and $V_{*}V_{*}^{T} = V_{r}V_{r}^{T}$ then $V_{r}^{T}V_{r} = I_{r}$.
	\label{lem:1.2}
\end{lemma}

\begin{align*}
	\|V_{*}^{T}V_{*} - I_{r}\|_{F}^{2} = 0\\
	\Tr\left(V_{*}^{T}V_{*}V_{*}^{T}V_{*}\right) + \Tr\left(I_{r}\right) - 2\Tr\left(V_{*}^{T}V_{*}\right) = 0
\end{align*}
Using the cyclic property of the trace operator we have,
\begin{align*}
	\Tr\left(V_{*}V_{*}^{T}V_{*}V_{*}^{T}\right) + \Tr\left(I_{r}\right) - 2\Tr\left(V_{*}V_{*}^{T}\right) = 0
\end{align*}
Using \textbf{Theorem \ref{thm:eym_mod}} $ V_{*}V_{*}^{T} = V_{r} V_{r}^{T}$,
\begin{align*}
	\Tr\left(V_{r} V_{r}^{T}V_{r}V_{r}^{T}\right) + \Tr\left(I_{r}\right) - 2\Tr\left(V_{r}V_{r}^{T}\right) = 0
\end{align*}
Again using the cyclic property of the trace operator,
\begin{align*}
	\Tr\left(V_{r}^{T}V_{r} V_{r}^{T}V_{r}\right) + \Tr\left(I_{r}\right) - 2\Tr\left(V_{r}^{T}V_{r}\right) = 0\\
	\Tr\left(V_{r}^{T}V_{r} V_{r}^{T}V_{r} + I_{r} - 2V_{r}^{T}V_{r}\right) = 0\\
	\|V_{r}^{T}V_{r} - I_{r}\|_{F} = 0
\end{align*}
Hence $V_{r}^{T}V_{r} = I_{r}$. 

\begin{remark}
Lemma \ref{lem:1.2} shows that for a rank-$r$ approximation of a rank-$f$ matrix $X$ such that $r>f$, the extracted right singular vectors are orthonormal when $V_{*}^{T}V_{*} = I_{r}$. This justifies the constraint $\tilde{V}^{T}\tilde{V} = I_{r}$ for the stage-1 minimization problem in Eq. \ref{eq:st1} and is numerically verified in \textbf{Section \ref{sec:inter}}.
\end{remark}

\begin{remark}
Note that Lemma \ref{lem:1.1} and \ref{lem:1.2} does not imply that $V_{*} = V_{r}$ instead $V_{*}\Theta_{r} = V_{r}$, where $\Theta_{r}$ is any real valued unitary matrix for the equality to hold true.
\end{remark}

\subsection{Stage 2}\label{app:thm:stage2}

\begin{theorem}
Given a rank-r matrix $XV_{*} \in \mathbb{R}^{(m\times r)}$ and an arbitrary, rank-r matrix $C  \in \mathbb{R}^{(m\times r)}$, following \textbf{Theorem} \ref{thm:orig_eym}, the tail energy of $XV_{*}$ with respect to $XV_{*}C$ is bounded as,
\begin{align*}
    \|XV_{*} - XV_{*}C\|_{F} \geq 0
\end{align*}
where the equality holds true if and only if  $C = I_{r}$.  
\label{thm:st21}
\end{theorem}

\begin{align*}
     \|XV_{*}(I_{r} - C)\|_{F} \geq 0 \\
     \|XV_{*}\|_{F}\|I_{r} - C\|_{F} \geq \|XV_{*}(I_{r} - C)\|_{F}\\
     \|XV_{*}\|_{F}\|I_{r} - C\|_{F} \geq 0
\end{align*}

Since $XV_{*}>0$, this implies equality is achieved if and only if $C = I_{r}$. 

\begin{lemma}
If $C = \Theta_{r}\Theta_{r}^{T}$, where $\Theta_{r} \in \mathbb{R}^{r\times r}$ is a rank-$r$ matrix such that $C = I_{r}$, then $\Theta_{r}$ is a real-valued unitary matrix in an $r$-dimensional Euclidean space.
\label{lem:2.1}
\end{lemma}

\begin{align*}
    \Theta_{r}\Theta_{r}^{T} = I_{r}\\
    \|\Theta_{r}\Theta_{r}^{T} -I_{r}\|^{2}_{F} = 0\\
    \Tr\left(\Theta_{r}\Theta_{r}^{T}\Theta_{r}\Theta_{r}^{T} + I_{r} - 2\Theta_{r}\Theta_{r}^{T}\right) = 0\\
    \Tr\left(\Theta_{r}\Theta_{r}^{T}\Theta_{r}\Theta_{r}^{T}\right) + \Tr\left(I_{r}\right) - 2\Tr\left(\Theta_{r}\Theta_{r}^{T}\right) = 0
\end{align*}
Using the cyclic property of the trace operator,
\begin{align*}
    \Tr\left(\Theta_{r}^{T}\Theta_{r}\Theta_{r}^{T}\Theta_{r}\right) + \Tr\left(I_{r}\right) - 2\Tr\left(\Theta_{r}^{T}\Theta_{r}\right) = 0\\
    \Tr\left((\Theta_{r}^{T}\Theta_{r} - I_{r})(\Theta_{r}^{T}\Theta_{r} - I_{r})^{T}\right) = 0\\
    \|\Theta_{r}^{T}\Theta_{r} - I_{r}\|_{F}^{2}=0
\end{align*}
This implies $\Theta_{r}^{T}\Theta_{r} = I_{r}$. Since $\Theta_{r}^{T}\Theta_{r} = \Theta_{r}\Theta_{r}^{T} = I_{r}$ this implies that $\Theta_{r}$ is a real-valued unitary matrix in the $r$-dimensional Euclidean space.

\begin{theorem}
Given a rank-r matrix $XV_{*} \in \mathbb{R}^{m\times r}$, such that $V_{*}V_{*}^{T} = V_{r}V_{r}^{T}$ where $V_{r}$ is a matrix with column vectors as the top-$r$ right singular vectors of $X$, and a real-valued unitary matrix $\Theta_{r} \in \mathbb{R}^{r\times r}$ then $(XV_{*}\Theta_{r})^{T}(XV_{*}\Theta_{r})$ is a diagonal matrix $\Sigma_{r}^{2}$ where $\Sigma_{r}^{2} = \mathrm{diag}(\sigma^{2}_{1},\sigma^{2}_{2},\cdots,\sigma^{2}_{r})$ and $\sigma_{i}$s are the top-$r$ singular values of $X$ if and only if $V_{*}\Theta_{r} = V_{r}$.
\label{thm:st22}
\end{theorem}

\begin{align*}
    (XV_{*}\Theta_{r})^{T}(XV_{*}\Theta_{r}) = \Sigma_{r}^{2}\\
    \Theta_{r}^{T}V_{*}^{T}(X^{T}X)V_{*}\Theta_{r} = \Sigma_{r}^{2}\\
    V_{*}^{T}X^{T}XV_{*} = \Theta_{r}\Sigma_{r}^{2}\Theta_{r}^{T}\\
    V_{*}V_{*}^{T}X^{T}XV_{*}V_{*}^{T} = V_{*}\Theta_{r}\Sigma_{r}^{2}\Theta_{r}^{T}V_{*}^{T}
\end{align*}
Using $V_{*}V_{*}^{T} = V_{r}V_{r}^{T}$ from \textbf{Theorem \ref{thm:eym_mod}},
\begin{align*}
    V_{r}V_{r}^{T}X^{T}XV_{r}V_{r}^{T} = V_{*}\Theta_{r}\Sigma_{r}^{2}\Theta_{r}^{T}V_{*}^{T}\\
    V_{r}\Sigma_{r}^{2}V_{r}^{T} = (V_{*}\Theta_{r})\Sigma_{r}^{2}(V_{*}\Theta_{r})^{T}\\
    (V_{r}-V_{*}\Theta_{r})\Sigma_{r}^{2}(V_{r}-V_{*}\Theta_{r})^{T} = 0\\
    \|(V_{r}-V_{*}\Theta_{r})\Sigma_{r}^{2}(V_{r}-V_{*}\Theta_{r})^{T}\|_{F} = 0\\
\end{align*}
Using Frobenius norm to bound the matrix product,
\begin{align*}
    \|\Sigma_{r}^{2}\|_{F}\|V_{r}-V_{*}\Theta_{r}\|^{2}_{F} \geq 0
\end{align*}
Since $\|\Sigma_{r}^{2}\|_{F}>0$, equality is achieved if and only if $V_{*}\Theta_{r}=V_{r}$.

\begin{remark}
Note that $\Theta_{r}$ is a rank-r unitary matrix wherein both rotation ($det(\Theta_{r}) = 1$) and reflection ($det(\Theta_{r}) = -1$) are valid since the order of the orthonormal vectors in the matrix $V_{*}\Theta_{r} = V_{r}$ do not alter $\|V_{r}-V_{*}\Theta_{r}\|_{F}$. In practice, $\Theta_{r}$ manifests itself predominantly as a rotation matrix during the iterative minimization using gradient descent.
\end{remark}

\section{Energy Minimization, Loss Surface Geometry, and Convergence} \label{app:converge}

\begin{figure}[th]
	\centering
	\includegraphics[width=0.3\linewidth]{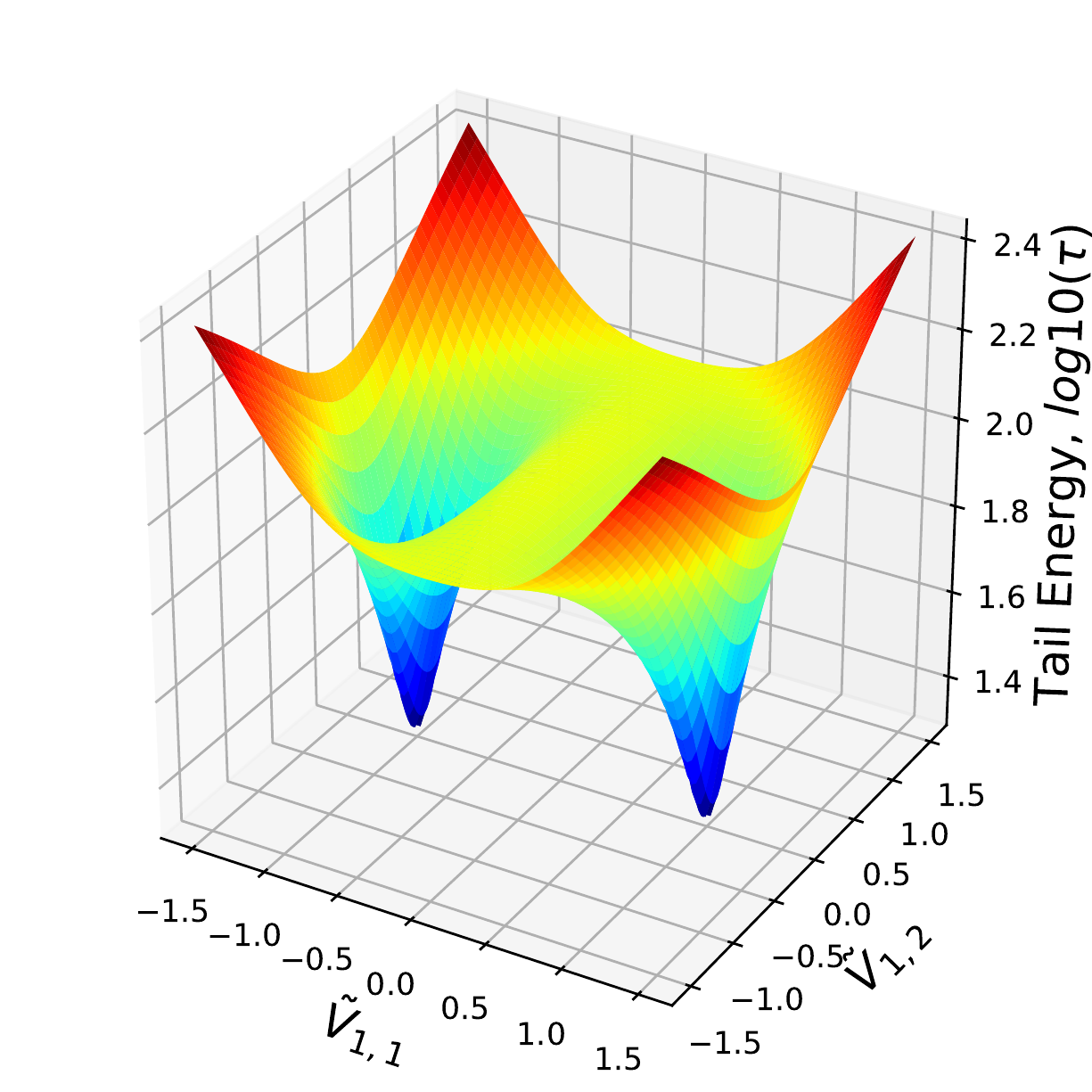}
	\caption{Surface plot for the bi-quadratic tail energy (log-scale) with 2 minima, 2 saddle points, and one maxima for $X = \mathrm{diag}(5,1)$.}
	% \Description{figure description}
	\label{fig:energy}
\end{figure}

In this section, we consider the energy minimization problem that constructs the projection space spanning the rank-$r$ sub-space of a given data matrix. For ease of visualization, we consider a $2\times2$ matrix $X = \mathrm{diag}(5,1)$ with singular values 5 and 1 corresponding to right singular vectors $v_{1} = [1,0]^{T}$ and $v_{2} = [0,1]^{T}$, respectively. Our objective here is to extract a rank $r=1$ approximation of this rank $f=2$ matrix $X$. Certainly, this corresponds to identifying the right singular vector $v_{1}$ with singular value $5$. The tail-energy surface (log-scale) corresponding to the $\| X\tilde{v}\tilde{v}^{T}-X\|_{F}$ is shown in \textbf{Fig. \ref{fig:energy}}. Here, $\tilde{v}$ is the test vector for a rank 1 approximation of X. The tail-energy is a bi-quadratic function in $\tilde{v}$ with 1 maximum, $2^{r}$ minima and $2^{r}$ saddle points, where $r$ is the desired low-rank approximation of a given data matrix. Furthermore, all minima have the same tail energy: a property of bi-quadratic functions. For the current specific example, the two equal tail-energy minima correspond to $v_{1}$ and $-v_{1}$, respectively.

Although the minimization problem is non-convex, convergence is guaranteed since any perturbed-gradient descent approach converges to either of the stable fixed points (minima). In other words, any test vector $\tilde{v}$ other than $v_{1}$ or $-v_{1}$ will increase the tail-energy and hence will not be the solution. The same argument applies for a high-dimensional dataset $X$ where a low rank ($r$) approximation is desired with the number of equal tail-energy minima corresponding to $2^{r}$ for all possible negative and positive combinations of the $r$ right singular vectors $v_{i}, \, i=1,\cdots r$.  In effect, the stage 1 minimization problem constructs a right projection space $\tilde{V} = span\{v_{1}, v_{2},\ldots,v_{r}\}$ that spans the top rank-r subspace of a given dataset $X$. A similar line of argument then applies to our stage 2 minimization problem as well. A mild limitation, that will be addressed in our future work, occurs when $X = \mathrm{diag}(5,5+\epsilon)$, $0\leq\epsilon<<1$, wherein the two right singular values cannot be resolved accurately (still better than Randomized SVD methods) without further considerations. This latter case, with near algebraic multiplicity in singular values is a special case for conventional SVD as well.

\section{Additional Experiments}

\subsection{Low Rank Approximation: Sandy Big Data}\label{app:sandy_low}

\begin{figure}[h]
    \centering
    \begin{subfigure}{.3\linewidth}
      \centering
      \includegraphics[width=\linewidth]{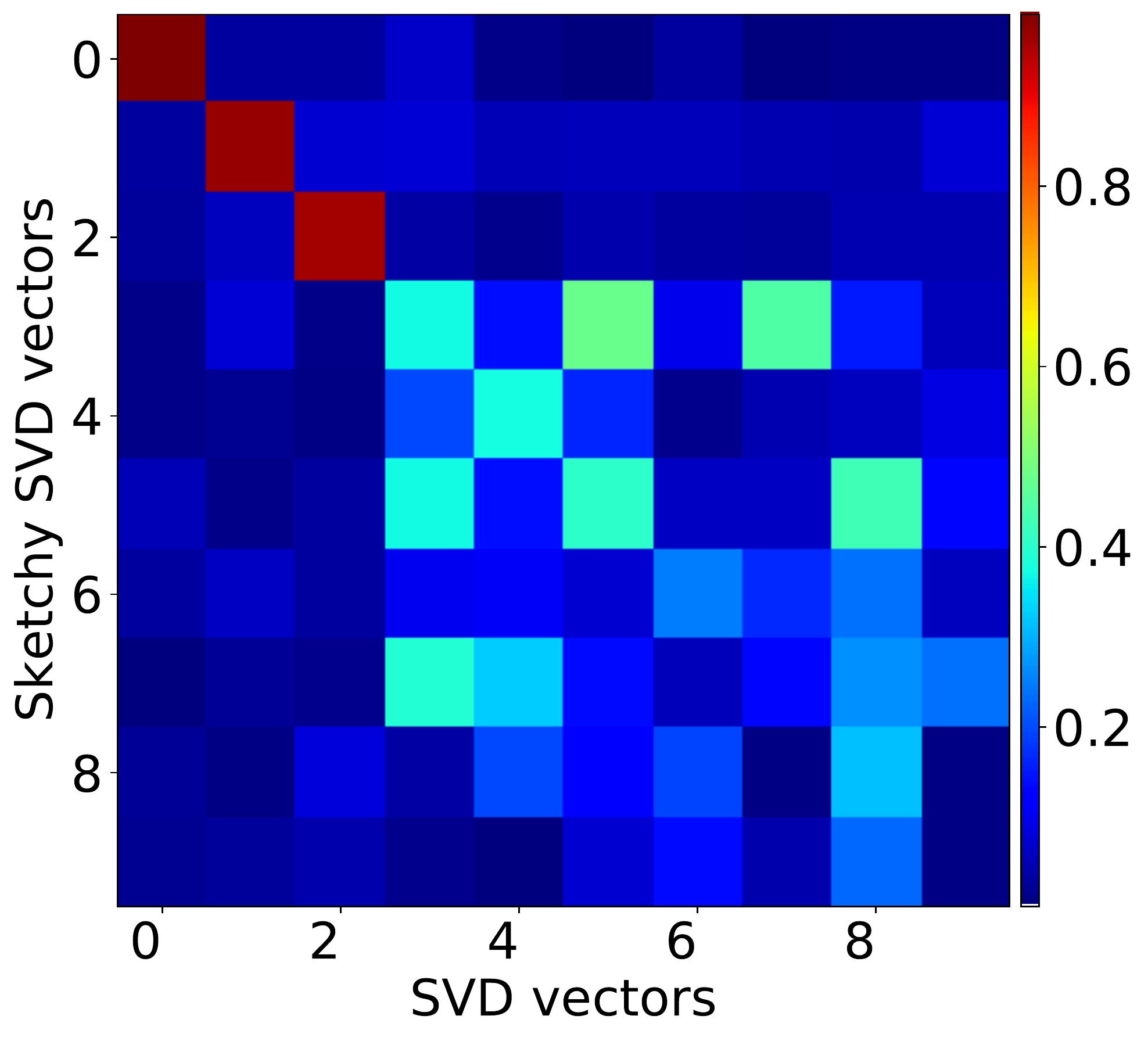}
      \caption{SketchySVD}
    \end{subfigure}
    \qquad
    \begin{subfigure}{.3\linewidth}
      \centering
      \includegraphics[width=\linewidth]{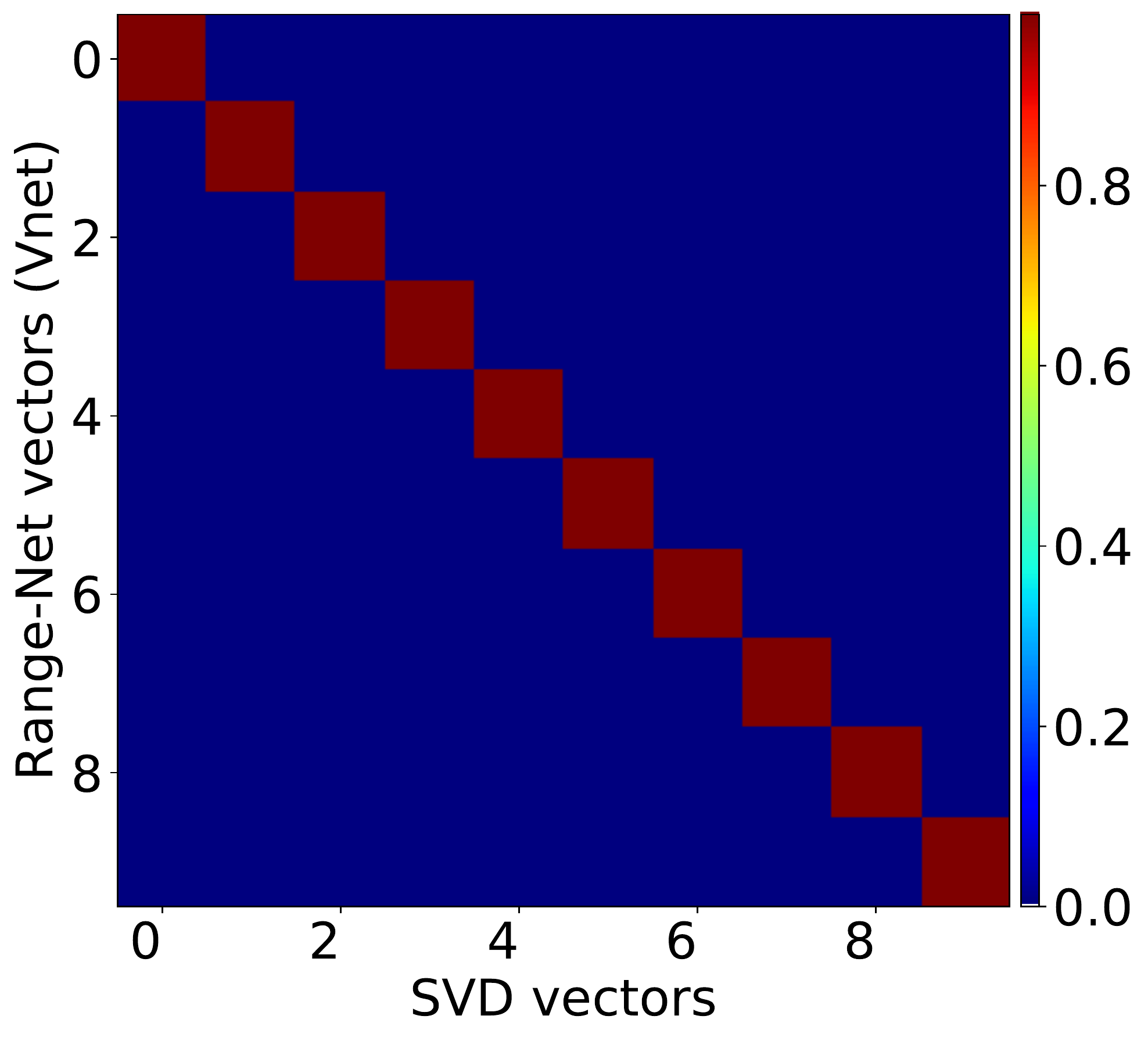}
      \caption{Range-Net}
    \end{subfigure}
    \caption{Cross-correlation between extracted and true (conventional SVD) right singular vectors for (a) SketchySVD and (b) Range-Net for a rank $r = 10$ approximation. SketchySVD deviates substantially after index 3 while Range-Net is in good agreement for all the 10 indices.}
    \label{fig:sandy_corr_low}
\end{figure}

In this experiment, we extract the rank-$10$ SVD factors using SketchySVD and Range-Net for the Sandy dataset. For this case, the oversampled ranks for sketchy SVD are $k = 4r+1 = 41$ and $s = 2k+1 = 83$ where $k,s \ll \min(m,n)$. As before,  \textbf{Fig. \ref{fig:sandy_corr_low}} show the cross-correlation between the extracted and true (conventional SVD) right singular vectors using SketchySVD and Range-Net. 

\begin{figure}[h]
    \centering
    \begin{subfigure}{.4\linewidth}
      \centering
      \includegraphics[width=\linewidth]{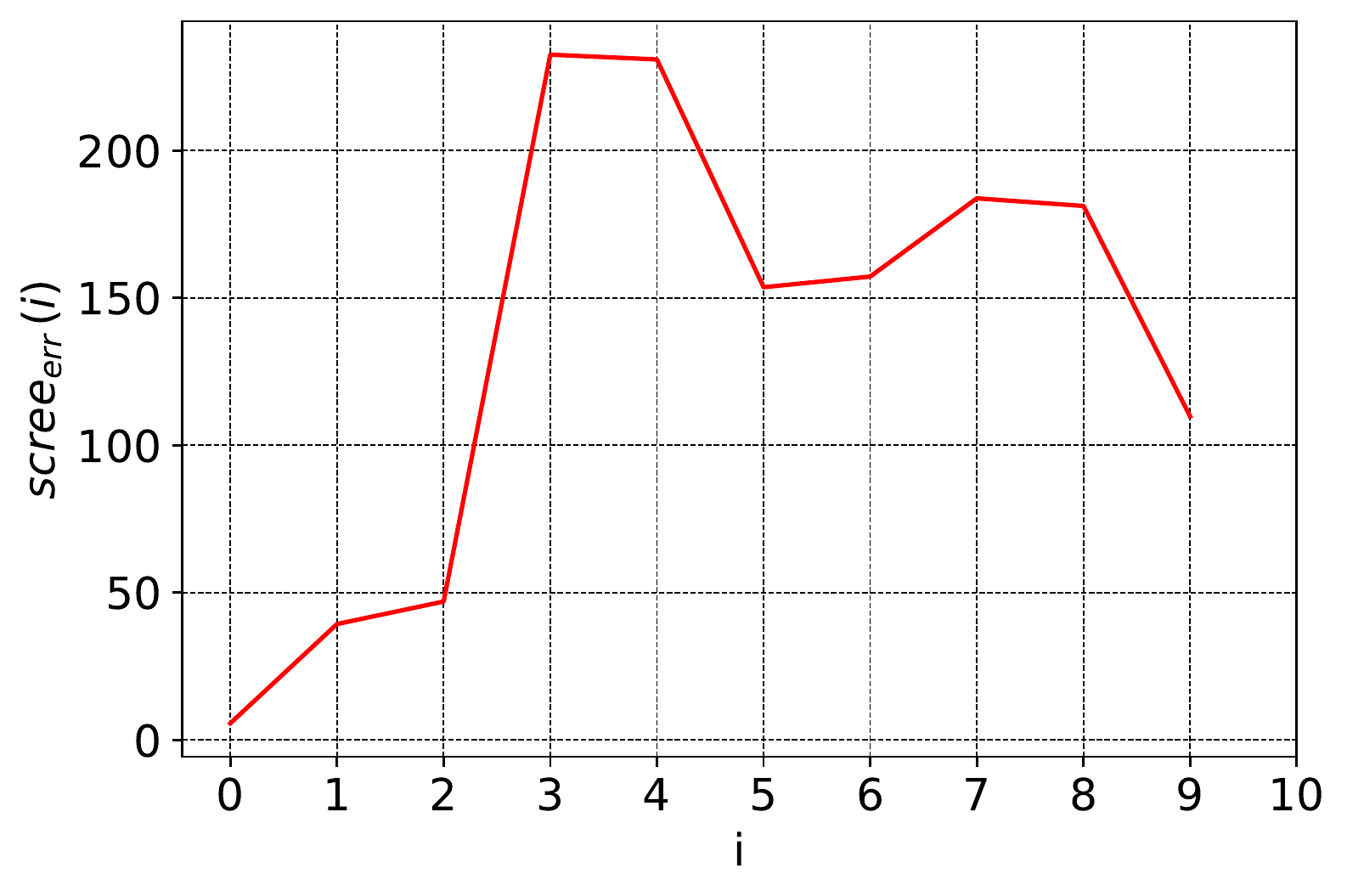}
      \caption{SketchySVD}
    \end{subfigure}
    \qquad
    \begin{subfigure}{.4\linewidth}
      \centering
      \includegraphics[width=\linewidth]{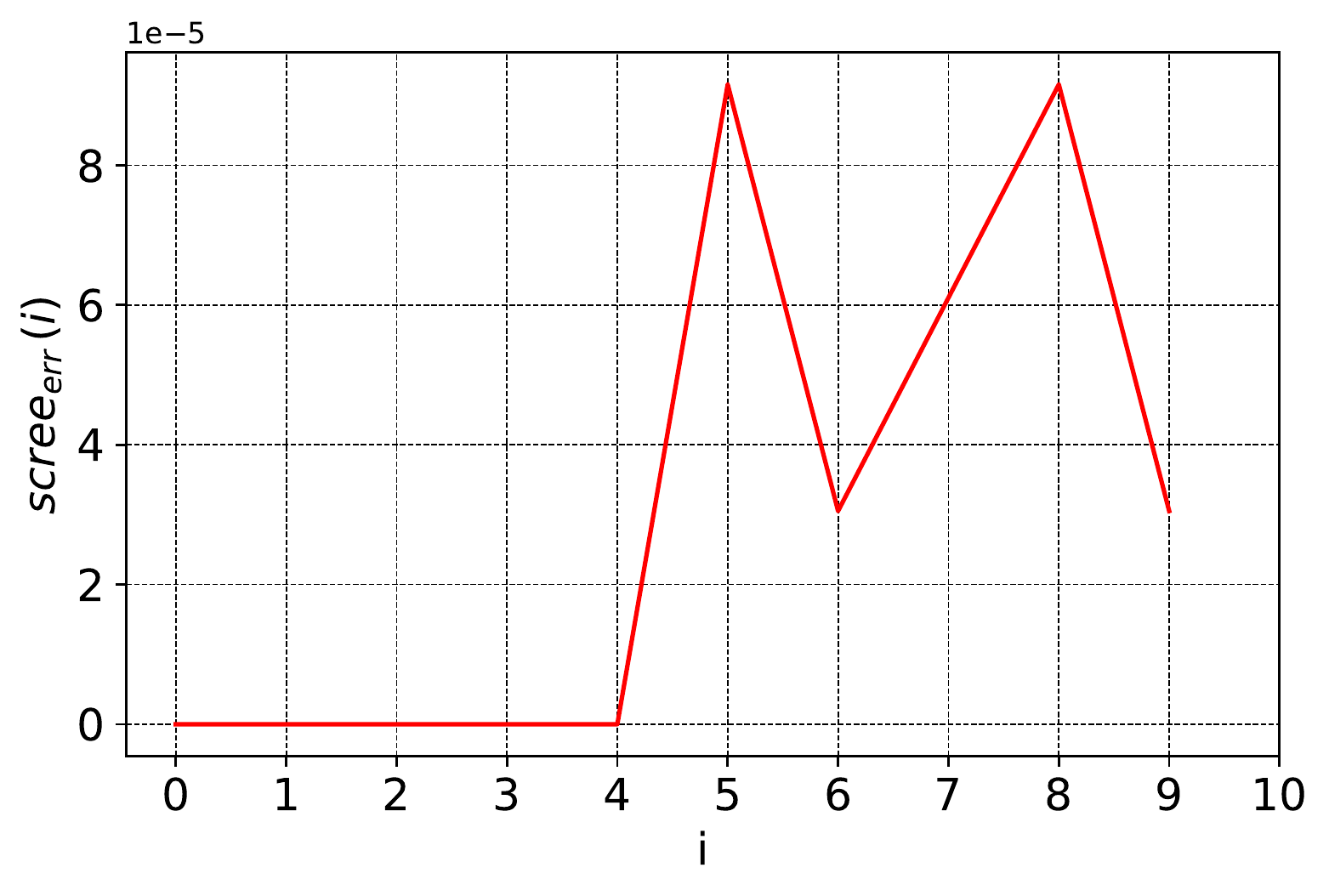}
      \caption{Range-Net}
    \end{subfigure}
    \caption{Scree-error in singular values for (a) SketchySVD and (b) Range-Net where a conventional SVD is used as the baseline in scree-error metric. Note that for Range-Net the error is at a scale of $10^{-5}$, $7$ orders of magnitude apart from SketchySVD ($10^2$).}
    \label{fig:sandy_scree_low}
\end{figure}

\textbf{Fig. \ref{fig:sandy_scree_low}} shows the scree error in the extracted singular values for the two methods with singular values from conventional SVD as the baseline.

\begin{figure}[h]
    \centering
    \begin{subfigure}{.25\linewidth}
      \centering
      \includegraphics[width=\linewidth]{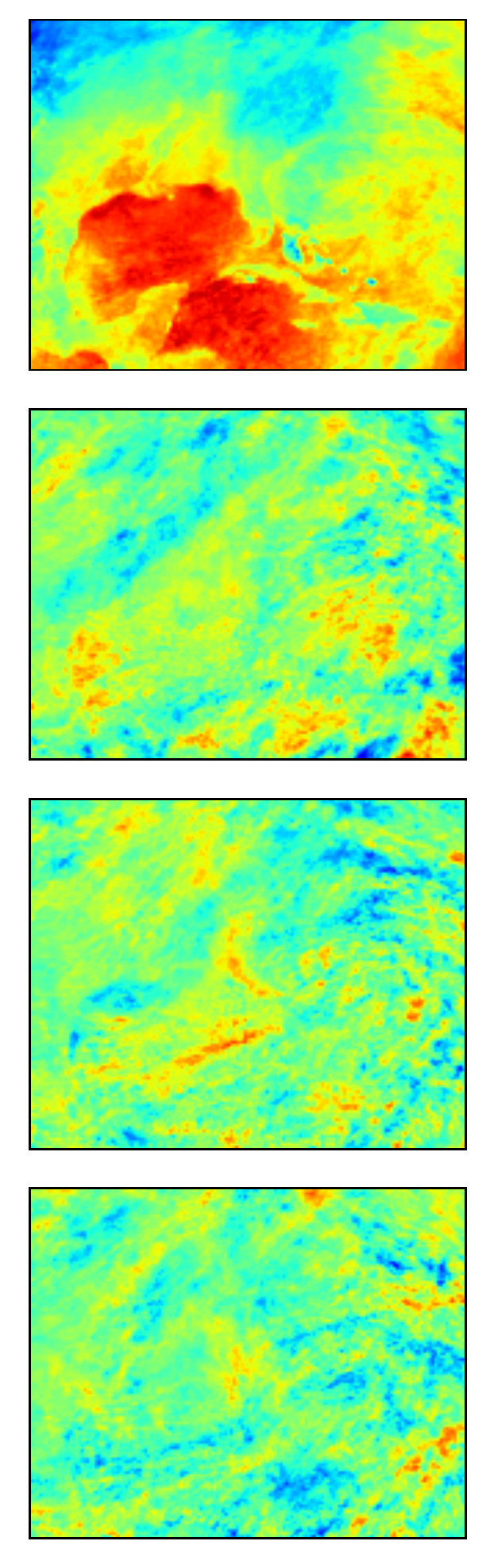}
      \caption{SketchySVD}
    \end{subfigure}
    \begin{subfigure}{.25\linewidth}
      \centering
      \includegraphics[width=\linewidth]{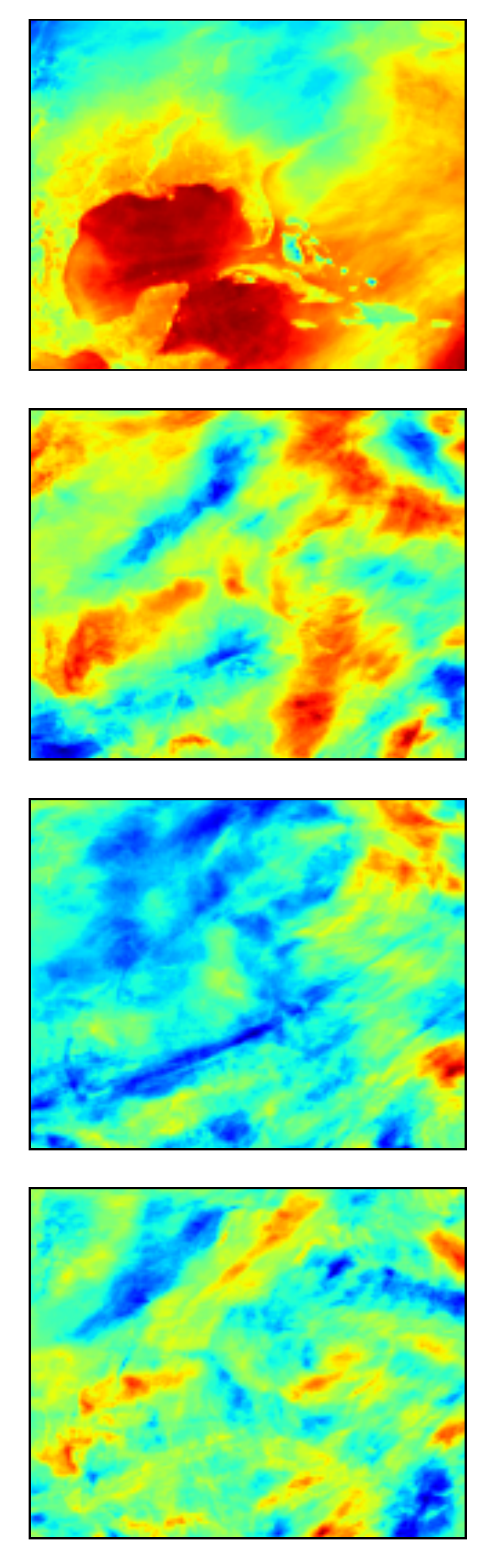}
      \caption{Truth}
    \end{subfigure}
    \begin{subfigure}{.25\linewidth}
      \centering
      \includegraphics[width=\linewidth]{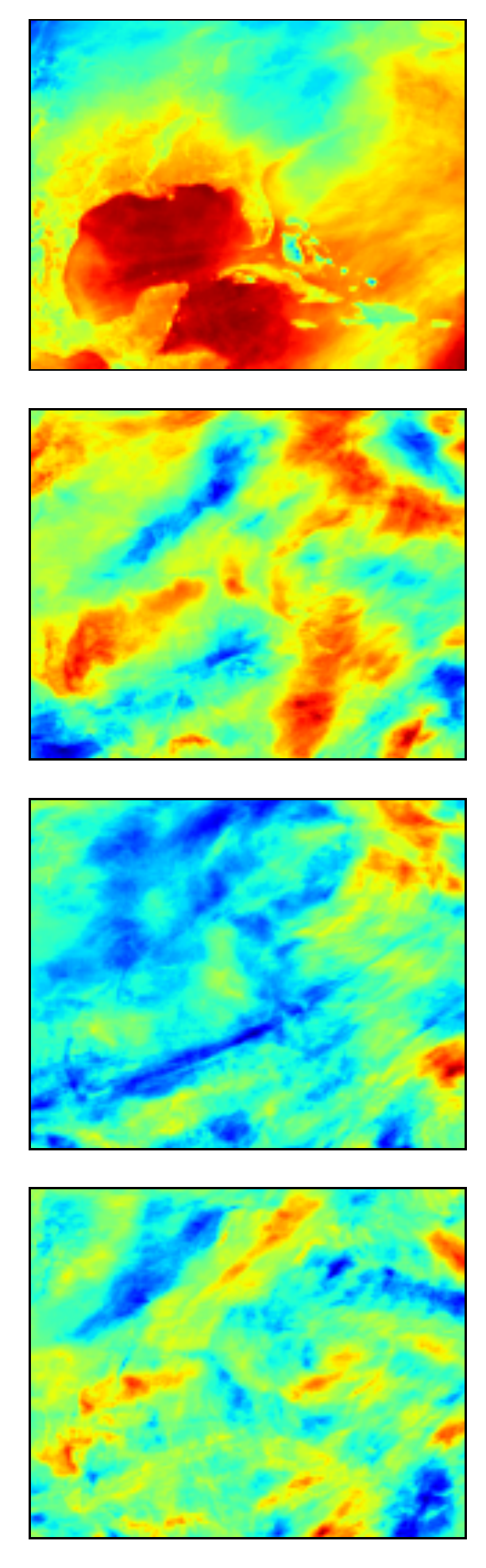}
      \caption{Range-Net}
    \end{subfigure}
    \caption{Reshaped $U_i$ indicative of dynamic modes, corresponding to $i=1,4,7,10$ for $r=10$. The dynamic mode approximation error stand out visually for SketchySVD for indices $50,100$. Our method does not have such artifacts.}
    \label{fig:sandy_sing_low}
\end{figure}

Finally, \textbf{Fig. \ref{fig:sandy_sing_low}} shows a comparison between extracted dynamic modes corresponding to indices $i = [1,4,7,10]$ from SketchySVD, conventional SVD, and Range-Net. One can easily see that Sketchy SVD extracted dynamic modes/right singular vectors deviate quite substantially for $i = [4,7,10]$.

\subsection{Feature Extraction: Graph (Eigen)} \label{app:graph}

\begin{table}[ht]
    \centering
    \caption{Description and Metrics of the Network Graphs}
    % \resizebox{\linewidth}{!}{
    \begin{tabular}{c|c|c|c|ccc} \toprule
        Dataset & Nodes & Edges & rank & $err_{fr}$ & $err_{sp}$ & $\chi^2_{err}$ \\ \midrule
        Airlines \cite{airlines} & 235 & 2101 & 200 & 0 & 0 & 0.011 \\
        Twitter \cite{airlines} & 3556 & 188712 & 200 & 0 & 0 & 0.014 \\
        Wikivote \cite{snapnets} & 8297 & 103689 & 200 & 0 & 0 & 0.027 \\
        Wikipedia \cite{levmuch} & 49728 & 941425 & 100 & 4.27e-6 & 1.23e-7 & 0.034 \\
        Slashdot \cite{snapnets} & 82168 & 948464 & 100 & 8.56e-6 & 6.92e-7 & 0.045 \\ \bottomrule
    \end{tabular}%}
    \label{tab:graph_data}
\end{table}

Large scale networks occur in many applications where SVD is primarily used to identify the most important nodes or as a pre-processing step to community detection. For these kind of graph based datasets, we either perform SVD or Eigen decomposition on the graph, depending on the format in which the data arrives. We demonstrate results on the following graphs of varying size, tabulated in \textbf{Table \ref{tab:graph_data}}. If the data arrives directly in the form of an adjacency matrix, we can perform SVD or Eigen decomposition on it directly. For cases, where an adjacency list is provided, a pre-processing step is required to convert the list representation in a sparse vector. Since Keras can handle sparse input data and sparse matrix operations, our method is trivially scalable to large sparse graphs. Since an Eigen decomposition problem is a special case of SVD, where the data matrix is symmetric positive semi-definite, our neural SVD solver is directly applicable.

% \begin{figure}[h]
%     \centering
%       \includegraphics[width=0.65\linewidth]{figs/graph_scr_net.pdf}  
%     \caption{Scree Error for Slashdot Graph}
%     \label{fig:scree_graph}
% \end{figure}

The benchmark was generated for smaller graphs using a conventional SVD solver. For larger graphs, a similar benchmark was constructed using the \textit{irlba} routine by \cite{baglama2005augmented}. \textbf{Table. \ref{tab:graph_data}} shows the error metrics for all the graphs, where consistently low values are observed for Frobenius and Spectral error metrics.

\section{Sketchy SVD Implementation}

A brief outline of the single-pass Sketch-SVD algorithm from \cite{tropp2019streaming}. Note that for the numbers reported in terms of storage, we implemented this code with additional memory optimization and sparse matrices.
\begin{algorithm}[ht]
	\caption{Sketchy SVD}
	\noindent \textbf{Input:} $X \in \mathbb{R}^{m \times n}$, $r:$ expected rank \\
	\noindent \textbf{Output:} $\tilde{X} \in \mathbb{R}^{m \times k}$ the approximated rank $k$-dim data
	\begin{algorithmic}[1]
		\State Initialize $k=4r+1,s=2k+1$ \Comment{Oversampling parameters}
		\State Projection maps: $\Upsilon \in \mathbb{R}^{k \times m}, \Omega \in \mathbb{R}^{k \times n}, \Phi \in \mathbb{R}^{s \times m}, \Psi \in \mathbb{R}^{s \times n}$
		\State Projection matrices: $A \in \mathbb{R}^{k \times n}, B \in \mathbb{R}^{m \times k}, Z \in \mathbb{R}^{s \times s}$ as empty
		\For{$i = 1:n$} \Comment{Streaming Update}
		\State Form $H \in \mathbb{R}^{m \times n}$ as a sparse empty matrix
		\State $H(i,:) = X(i,:)$ \Comment{Streamed columns}
		\State $A \leftarrow A + \Upsilon H$ \Comment{Update Co-Range}
		\State $B \leftarrow B + H \Omega^T$ \Comment{Update Range}
		\State $Z \leftarrow Z + \Phi H \Psi^T$ \Comment{Update Core Sketch}
		\EndFor
		\State $Q \in \mathbb{R}^{m \times k} \leftarrow qr\_econ(B)$ \Comment{Basis for Range}
		\State $P \in \mathbb{R}^{n \times k} \leftarrow qr\_econ(A^T)$ \Comment{Basis for Co-Range}
		\State $C \in \mathbb{R}^{s \times s} \leftarrow ((\Phi Q) \setminus Z) / (\Psi P)$ \Comment{Core Matrix}
		\State $[U,\Sigma,V^T] \leftarrow svd(C)$ \Comment{Full SVD of Core Matrix}
		\State $\Sigma \in \mathbb{R}^{r \times r} \leftarrow \Sigma[1:r,1:r]$ \Comment{Pick top $r$}
		\State $U \in \mathbb{R}^{k \times r} \leftarrow U[:,1:r]$ \Comment{Pick top $r$}
		\State $V^T \in \mathbb{R}^{r \times k} \leftarrow V^T[1:r,:]$ \Comment{Pick top $r$}
		\State $U \in \mathbb{R}^{m \times r} \leftarrow QU$ \Comment{Project to Row Space}
		\State $V^T \in \mathbb{R}^{r \times n} \leftarrow PV^T$ \Comment{Project to Column Space}
		\BState $\tilde{X} \in \mathbb{R}^{m \times n} \leftarrow U \Sigma V^T$ \Comment{Approximation}
	\end{algorithmic} \label{alg:sketchy}
\end{algorithm}

\end{document}